\documentclass{amsart}%
\usepackage{graphicx}
\usepackage{amscd}
\usepackage{amsmath}%
\usepackage{amsfonts}%
\usepackage{amssymb}
\newtheorem{theorem}{Theorem}[section]
\newtheorem{corollary}[theorem]{Corollary}
\newtheorem{lemma}[theorem]{Lemma}
\newtheorem{proposition}[theorem]{Proposition}
\theoremstyle{definition}
\newtheorem{definition}[theorem]{Definition}
\newtheorem{remark}[theorem]{Remark}

\newtheorem{observation}[theorem]{Observation}
\newtheorem{example}[theorem]{Example}

\theoremstyle{remark}
\newtheorem*{acknowledgements}{Acknowledgements}
\newlength{\submatrixwidth}
\newlength{\qedskip}
\newlength{\baselineadjust}
\newcounter{enumlink}
\renewcommand{\theenumi}{\roman{enumi}}
\renewcommand{\labelenumi}{(\theenumi)}
\numberwithin{equation}{section}
\def\hooklongrightarrow{\mathrel{\lhook\joinrel\longrightarrow}}

\newbox\frogdown
\newlength\frogdrop
\def\hookdownarrow{\setlength{\unitlength}{0.4pt}\setbox\frogdown=\hbox to 0pt{\hss $\displaystyle \downarrow $\hss }\setlength{\frogdrop}{2.5\ht\frogdown}\raisebox{0pt}[5\unitlength][\frogdrop]
{\begin{picture}(0,5)(10,0)
\put(5,0){\oval(10,10)[t]}
\end{picture}\lower\ht\frogdown\box\frogdown}}
\def\openone
{\mathchoice
{\hbox{\upshape \small1\kern-3.3pt\normalsize1}}
{\hbox{\upshape \small1\kern-3.3pt\normalsize1}}
{\hbox{\upshape \tiny1\kern-2.3pt\SMALL1}}
{\hbox{\upshape \Tiny1\kern-2pt\tiny1}}}
\makeatletter
\newbox\ipbox
\newcommand{\ip}[2]{\left\langle #1\mathrel{\mathchoice
{\setbox\ipbox=\hbox{$\displaystyle \left\langle\mathstrut #1#2\right\rangle$}
\vrule height\ht\ipbox width0.25pt depth\dp\ipbox}
{\setbox\ipbox=\hbox{$\textstyle \left\langle\mathstrut #1#2\right\rangle$}
\vrule height\ht\ipbox width0.25pt depth\dp\ipbox}
{\setbox\ipbox=\hbox{$\scriptstyle \left\langle\mathstrut #1#2\right\rangle$}
\vrule height\ht\ipbox width0.25pt depth\dp\ipbox}
{\setbox\ipbox=\hbox{$\scriptscriptstyle \left\langle\mathstrut #1#2\right\rangle$}
\vrule height\ht\ipbox width0.25pt depth\dp\ipbox}
} #2\right\rangle}
\newcommand{\diracb}[1]{\left\langle #1\mathrel{\mathchoice
{\setbox\ipbox=\hbox{$\displaystyle \left\langle\mathstrut #1\right.$}
\vrule height\ht\ipbox width0.25pt depth\dp\ipbox}
{\setbox\ipbox=\hbox{$\textstyle \left\langle\mathstrut #1\right.$}
\vrule height\ht\ipbox width0.25pt depth\dp\ipbox}
{\setbox\ipbox=\hbox{$\scriptstyle \left\langle\mathstrut #1\right.$}
\vrule height\ht\ipbox width0.25pt depth\dp\ipbox}
{\setbox\ipbox=\hbox{$\scriptscriptstyle \left\langle\mathstrut #1\right.$}
\vrule height\ht\ipbox width0.25pt depth\dp\ipbox}
}\right. }
\newcommand{\dirack}[1]{\left. \mathrel{\mathchoice
{\setbox\ipbox=\hbox{$\displaystyle \left.\mathstrut #1\right\rangle$}
\vrule height\ht\ipbox width0.25pt depth\dp\ipbox}
{\setbox\ipbox=\hbox{$\textstyle \left.\mathstrut #1\right\rangle$}
\vrule height\ht\ipbox width0.25pt depth\dp\ipbox}
{\setbox\ipbox=\hbox{$\scriptstyle \left.\mathstrut #1\right\rangle$}
\vrule height\ht\ipbox width0.25pt depth\dp\ipbox}
{\setbox\ipbox=\hbox{$\scriptscriptstyle \left.\mathstrut #1\right\rangle$}
\vrule height\ht\ipbox width0.25pt depth\dp\ipbox}
} #1\right\rangle}
\newcommand{\rip}[2]{\left( #1\mathrel{\mathchoice
{\setbox\ipbox=\hbox{$\displaystyle \left(\mathstrut #1#2\right)$}
\vrule height\ht\ipbox width0.25pt depth\dp\ipbox}
{\setbox\ipbox=\hbox{$\textstyle \left(\mathstrut #1#2\right)$}
\vrule height\ht\ipbox width0.25pt depth\dp\ipbox}
{\setbox\ipbox=\hbox{$\scriptstyle \left(\mathstrut #1#2\right)$}
\vrule height\ht\ipbox width0.25pt depth\dp\ipbox}
{\setbox\ipbox=\hbox{$\scriptscriptstyle \left(\mathstrut #1#2\right)$}
\vrule height\ht\ipbox width0.25pt depth\dp\ipbox}
} #2\right)}
\let\subsubsubsectionname\@empty
\newcounter{subsubsubsection}[subsubsection]

\def\l@subsubsubsection{\@tocline{4}{0pt}{1pc}{9pc}{}}

\def\subsubsubsection{\@startsection{subsubsubsection}{4}%
  \z@{.5\linespacing\@plus.7\linespacing}{-.5em}%
  {\normalfont\itshape}}
\AtBeginDocument{%
  \@for\@tempa:=-1,0,1,2,3,4\do{%
    \@ifundefined{r@tocindent\@tempa}{%
      \@xp\gdef\csname r@tocindent\@tempa\endcsname{0pt}}{}%
  }%
}
\def\@writetocindents{%
  \begingroup
  \@for\@tempa:=-1,0,1,2,3,4\do{%
    \immediate\write\@auxout{%
      \string\newlabel{tocindent\@tempa}{%
        \csname r@tocindent\@tempa\endcsname}}%
  }%
  \endgroup}
\makeatother
\def\LaTeXparent#1{}%
\def\ChildStyles#1{}%

\input cyracc.def

\begin{document}
\title[The isomorphism problem for stationary AF-algebras]{Decidability of the isomorphism problem for stationary AF-algebras and the
associated ordered simple dimension groups}
\author{Ola~Bratteli}
\address[O.~Bratteli]{Department of Mathematics\\
University of Oslo\\
PB 1053 -- Blindern\\
N-0316 Oslo\\
Norway}
\email{bratteli@math.uio.no}
\author{Palle~E.~T.~ Jorgensen}
\address[P.E.~T.~ Jorgensen]{Department of Mathematics\\
The University of Iowa\\
14 MacLean Hall\\
Iowa City, IA 52242-1419\\
U.S.A.}
\email{jorgen@math.uiowa.edu}
\author{Ki Hang Kim}
\address[K.H. Kim and F. Roush]{Mathematics Research Group\\
Alabama State University\\
Montgomery AL~36101--0271\\
U.S.A.}
\email{kkim@asunet.alasu.edu, froush@asunet.alasu.edu}
\author{Fred Roush}
\thanks{Work supported by U.S. N.S.F. grants DMS9900265 (K.H.K., F.R.) and DMS9700130
(P.E.T.J.), and by the Norwegian N.F.R. (Norges Forskningsr\aa d) and the
University of Oslo, Norway (O.B., P.E.T.J.)}

\begin{abstract}
The notion of isomorphism of stable AF-$C^{\ast}$-algebras is considered in
this paper in the case when the corresponding Bratteli diagram is stationary,
i.e., is associated with a single square primitive incidence matrix.
\mbox{$C^{\ast}$-isomorphism} induces an equivalence relation on these
matrices, called \mbox{$C^{\ast}$-equivalence}. We show that the associated
isomorphism equivalence problem is decidable, i.e., there is an algorithm that
can be used to check in a finite number of steps whether two given primitive
matrices are \mbox{$C^{\ast}$-equivalent} or not.

\end{abstract}
\maketitle

\section*{Introduction}

%
%
%
In \cite{BJKR98} we studied isomorphism of the stable AF-algebras associated
with constant square primitive nonsingular incidence matrices. This
isomorphism is called $C^{\ast}$-equivalence of the matrices in \cite{BJKR98}
and weak equivalence of the (transposed) matrices in \cite{SwVo00}. In this
paper we prove that the isomorphism problem in this setting is decidable, even
when the assumption of nonsingularity is removed. The decision procedure is
spelled out explicitly in Section \ref{SX}. This result was announced in
\cite{BJKR98}, and it is interesting in view of the fact that the
corresponding problem for non-constant incidence matrices is undecidable
\cite{MuPa98}. That isomorphism is decidable means that there is an algorithm
that can be used to decide, in a finite number of steps, whether two given
primitive matrices are $C^{\ast}$-equivalent or not. (See below.)

The significance of this result goes well beyond the theory of AF-algebras,
since the result may be viewed as a decision procedure for isomorphism of the
ordered simple dimension groups associated to the AF-algebras, and this class
of groups is important for a variety of other problems, especially in symbolic
and topological dynamics, see \cite{PaTa95}, \cite{Han81}, \cite{BMT87} and
\cite{Kit98}. The decision result is a fundamental and nontrivial fact one
wants in all these applications.

Bratteli diagrams were introduced in \cite{Bra72} with a view to understanding
the structure and the classification of those $C^{\ast}$-algebras which arise
as inductive limits of finite-dimensional $C^{\ast}$-algebras, the so-called
AF-algebras. In fact, the equivalence relation on Bratteli diagrams which is
generated by the operation of telescoping is a complete $C^{\ast}$-isomorphism
invariant for the AF-algebras; see \cite[Remark 5.6]{BJO99}. It is the
decidability of this isomorphism problem in the case of stationary Bratteli
diagrams which is our main result here. The diagrams are called
\emph{stationary} if the incidence matrix is constant; in the general case it
is not constant, but varies from one level to the next. However, it was the
stationary class of AF-algebras which came from the problem addressed in
\cite{BJO99}, and while special, this subfamily is still general enough for
the study of substitution dynamical systems, as noted in \cite{DHS99}.
Consider, for example, a substitution dynamical system $\sigma$ (letters to
words) derived from a given alphabet $S$ of size $N$. For $i,j\in S$, let
$a_{ij}$ count the number of occurrences of $i$ in the word $\sigma(j)$,
resulting from the substitution $\sigma$, and let $A$ be the corresponding
matrix with dimension group $G(A)$ (see (\ref{eqInt.10})). In \cite{DHS99},
the co-authors use $G(A)$ in their classification of these systems, which may
also be realized as shift dynamical systems on the paths in the corresponding
Bratteli diagrams. These systems have significance in formal languages,
quasi-crystals, aperiodic tilings of the plane \cite{Rad99}, and
$p$-recognizable sets of numbers. Hence the classification we address here has
some bearing not only on the original setting of AF-algebras, but also on
recent developments in dynamical systems. For a survey of other dynamical
system classifications related to more standard shifts than those considered
in \cite{DHS99}, and the relation of our present classification to these, see
\cite{BJKR98}. In particular, it is explained in \cite{BJKR98} that the notion
of $C^{\ast}$-equivalence of two primitive nonsingular matrices is strictly
weaker than shift equivalence, strong shift equivalence, or elementary shift
equivalence. Specifically, formula (\ref{eqInt.2}) below shows that $C^{\ast}%
$-equivalence may be expressed also as a certain system of matrix
factorizations, but these conditions for $C^{\ast}$-equivalence are less
restrictive than those which define shift equivalence \cite[Proposition
2]{BJKR98}. This means that some techniques which are common in the study of
shift equivalence, see, e.g., \cite{BMT87}, are also common in the study of
isomorphism of $C^{\ast}$-algebras. The dimension group is one such tool, see
\cite{Ell76}, \cite{Eff81}.

Our approach is based on studying isomorphism of ordered dimension groups (the
order is essential!). We introduce those groups in (\ref{eqInt.6}%
)--(\ref{eqInt.11}), and we formulate the associated isomorphism problem. We
then go on to prove that this problem is decidable, in Theorem \ref{ThmDec.7}.
A general algorithm which can be used to decide whether or not two primitive
matrices $A$, $B$ are $C^{\ast}$-equivalent is spelled out point by point in
Section \ref{SX}.

After decidability, the next question is a presentation of the answer in terms
of numerical invariants. We take this up in Sections \ref{Cas}--\ref{GZN},
which are a continuation of \cite{BJO99}. Here the answers are not yet
complete, so we present in Section \ref{Cas} (Proposition \ref{Prop1} and
Corollary \ref{ProFred.2.7.99}) a subclass of incidence matrices for which the
$C^{\ast}$-equivalence question is decided by the value of a numerical
invariant. The matrices $A$ in the subclass allow a direct-sum decomposition,
$A=A_{0}\oplus\left(  \lambda\right)  $, such that $A_{0}$ is unimodular up to
sign, and $\left(  \lambda\right)  $ is multiplication by the
Perron--Frobenius eigenvalue $\lambda$ on the one-dimensional subspace spanned
by the right Perron--Frobenius eigenvector. This property is equivalent to
$\left|  \det A\right|  =\lambda$.

Section \ref{Som} and Section \ref{Str} address symmetry properties, pointing
out that there are nonsymmetric primitive incidence matrices $A$ which are
$C^{\ast}$-equivalent to $A^{\operatorname*{tr}}$, the transposed matrix. But
even in the $2$-by-$2$ case, there are also examples where $A$ and
$A^{\operatorname*{tr}}$ are not $C^{\ast}$-equivalent. The related symmetry
question for shift equivalence comes from the issue of reversibility for
topological Markov chains, which was studied in \cite{PaTu82} and
\cite{CuKr80}.

While the dimension group $G\left(  A\right)  $ associated with an incidence
matrix $A$ is torsion-free, it has a certain torsion group quotient $G\left(
A\right)  /L$ by a lattice $L$ in $G\left(  A\right)  $. We show in
Proposition \ref{Proposition1} that this quotient is natural in the sense that
it is an invariant. It is well known that abelian torsion groups have
explicitly computable and complete numerical invariants, and these invariants
are thus also invariants for the dimension group (but not complete because
they do not reflect order and some of the group structure). In this case they
take an especially simple form, and they can be read off from the
characteristic polynomial. This is proved in Section \ref{GZN}. Section
\ref{Str} presents a formulation of $C^{\ast}$-equivalence for matrices $A$,
$B$ in terms of a certain explicit matrix factorization $B=CAD$, where the two
factors $C$, $D$ are specified in the statement of the result, Theorem
\ref{Theorem5Sep15}.

By decidability of a class of problems, we will here mean that there is an
algorithm (which could be converted into a computer program) to solve the
problem \cite{Her69,Her78,Knu81}. There may be no simple way to tell how many
steps the algorithm will use, but it must eventually terminate in all cases.
This is equivalent to saying that there is a Turing machine, which given the
necessary inputs (two matrices here), will give an output which here is zero
or one accordingly as the problem has an answer ``No'' or ``Yes''. The theory
of algorithmic decidability begins essentially with the proof that the halting
problem, the problem of whether an arbitrary Turing machine on a given input
will halt, is algorithmically undecidable; a result equivalent to this was
proved by G\"{o}del, though the theory was cast into different forms by
Church, Kleene, and Turing. Its high-water mark was the proof by Davis,
Matijasevi\v{c}, Putnam, and Robinson \cite{DMR76} that diophantine equations
over the integers are algorithmically unsolvable (Hilbert's Tenth Problem).
Since then many other problems have been proved undecidable (such as the
result of \cite{MuPa98} on a different class of $C^{\ast}$-algebras), though
others like the diophantine problem over $\mathbb{Q}$, have resisted all
efforts. On the other hand, major decidability results have appeared too, such
as the proof by Ax and Kochen \cite{AxKo65a,AxKo65b,AxKo66} that it is
decidable whether a given system of diophantine equations is solvable
simultaneously over every $p$-adic field or ring, results on power series
rings, Rabin's result on the theory of the Cantor set, and the
Grunewald--Segal result \cite{GrSe80a,GrSe80b} that isomorphism of forms over
algebraic number rings is decidable (which is related to our work here as well
as the proof in \cite{KiRo79} that shift equivalence is decidable). Many
conjectured decidability results remain open, such as the question of whether
abelian and hyperbolic algebraic varieties have a rational point \cite[Parts C
and F4]{HiSi00}.

The classical treatise on decidability in the context of algorithmic algebraic
number theory is the book \cite{PoZa97}, and we will, sharing the view of
those authors, not try to give a definition of algorithm in terms of
mathematical logic.

There are some general blanket references which we will use throughout the
paper: \cite{Wei98} and \cite{BoSh66} on algebraic number theory,
\cite{PoZa97} on algorithms of algebraic number theory, \cite{New72} on
integral matrices and their factorizations, and \cite{Kit98}, \cite{Wag99} on
symbolic dynamics. Especially \cite{PoZa97}, \cite{Wei98}, and \cite{New72}
are used frequently in the proofs to follow, each one containing algorithmic
constructions which we cite as they are needed. The proofs involve diverse
areas of mathematics which are not always thought to be directly related. They
fall at the interface of techniques from these different subjects. For that
reason, we include a bit more detail and discussion than is customary in a
paper which does not cut across boundaries between fields.

\section{\label{Int}Equivalent isomorphism conditions}

Recall from \cite{BJKR98} that two matrices $A$, $B$ with nonnegative integer
matrix entries are said to be $C^{\ast}$-equivalent if there exist two
sequences $n_{1},n_{2},\dots$ and $m_{1},m_{2},\dots$ of natural numbers and
two sequences of matrices $J\left(  1\right)  ,J\left(  2\right)  ,\dots$ and
$K\left(  1\right)  ,K\left(  2\right)  ,\dots$ with nonnegative integer
matrix entries such that the diagram (\ref{eqInt.1}) below commutes.%

\begin{equation}
\setlength{\unitlength}{78bp}%
\begin{minipage} {1.67\unitlength}\begin{picture}(1.67,5)(-0.4,-4.895) \put
(-0.093,-4.3565){\includegraphics
[bb=133 0 215 348,height=348bp,width=82bp]{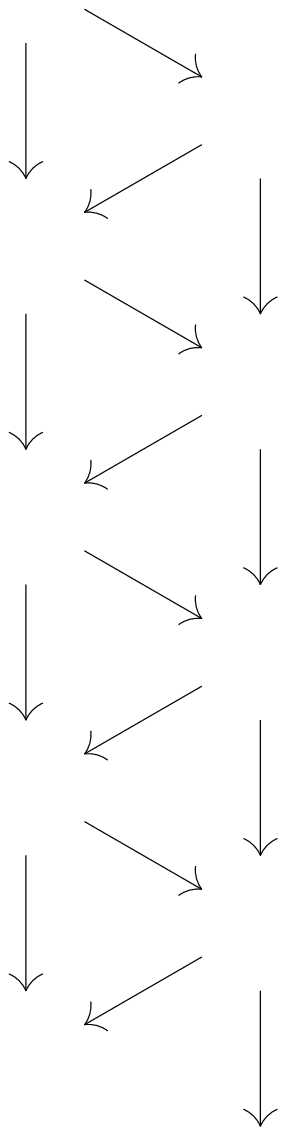}} \put(0,0){\makebox
(0,0){$\bullet$}} \put(0.866,-0.5){\makebox(0,0){$\bullet$}} \put
(0,-1){\makebox(0,0){$\bullet$}} \put(0.866,-1.5){\makebox(0,0){$\bullet$}%
}\put(0,-2){\makebox(0,0){$\bullet$}} \put(0.866,-2.5){\makebox(0,0){$\bullet
$}} \put(0,-3){\makebox(0,0){$\bullet$}} \put(0.866,-3.5){\makebox
(0,0){$\bullet$}} \put(0,-4){\makebox(0,0)[t]{$\vdots$}} \put
(0.866,-4.5){\makebox(0,0)[t]{$\vdots$}} \put(0.433,-0.225){\makebox
(0,0)[bl]{$J\left(1\right)$}} \put(0.433,-1.225){\makebox(0,0)[bl]{$J\left
(2\right)$}} \put(0.433,-2.225){\makebox(0,0)[bl]{$J\left(3\right)$}}%
\put(0.433,-3.225){\makebox(0,0)[bl]{$J\left(4\right)$}} \put
(0.433,-0.725){\makebox(0,0)[br]{$K\left(1\right)$}} \put
(0.433,-1.725){\makebox(0,0)[br]{$K\left(2\right)$}} \put
(0.433,-2.725){\makebox(0,0)[br]{$K\left(3\right)$}} \put
(-0.025,-0.5){\makebox(0,0)[r]{$A^{n_{1}}$}} \put(-0.025,-1.5){\makebox
(0,0)[r]{$A^{n_{2}}$}} \put(-0.025,-2.5){\makebox(0,0)[r]{$A^{n_{3}}$}}%
\put(0.891,-1){\makebox(0,0)[l]{$B^{m_{1}}$}} \put(0.891,-2){\makebox
(0,0)[l]{$B^{m_{2}}$}} \put(0.891,-3){\makebox(0,0)[l]{$B^{m_{3}}$}%
}\end{picture} \end{minipage}
\label{eqInt.1}%
\end{equation}

The diagram expresses the following two identities:%
\begin{equation}
A^{n_{k}}=K\left(  k\right)  J\left(  k\right)  ,\qquad B^{m_{k}}=J\left(
k+1\right)  K\left(  k\right)  , \label{eqInt.2}%
\end{equation}
for $k=1,2,\dots$. This corresponds to isomorphism of the associated stable
AF-algebras \cite{BJKR98,Bra72}, and it corresponds to homeomorphism of
one-di\-men\-sion\-al connected orientable hyperbolic attractors of
diffeomorphisms of manifolds by \cite{Jac97}; see also \cite{SwVo00}. We will
assume throughout that $A$ and $B$ are primitive square matrices (i.e.,
sufficiently high powers have only strictly positive matrix entries). For the
rest of this section we will also assume that $A$ and $B$ are nonsingular, but
this extra condition can be dispensed with by a remedy described in Section
\ref{SY}. So assume that $A$ and $B$ are nonsingular, and hence $C^{\ast}%
$-equivalence implies that they have the same dimension $N$, because $N$ is
the rank of the associated dimension group \cite{BJO99}. (We will argue in
Section \ref{Red} that the class of AF-algebras we obtain in this manner will
no longer be the same if $A$ and $B$ are merely required to be primitive but
not necessarily nonsingular. This does not contradict the results in Section
\ref{SY}, because the matrices replacing $A$, $B$ there no longer have
positive matrix entries, and the order is defined in a different manner.) In
this case we note that $J\left(  1\right)  $ and the sequences $n_{1},\dots$
and $m_{1},\dots$ determine all other $K\left(  k\right)  $ and $J\left(
j\right)  $ from (\ref{eqInt.1}), i.e.,
\begin{equation}%
\begin{aligned}
K\left( 1\right) &=A^{n_{1}}J\left( 1\right) ^{-1},  \\
J\left( 2\right) &=B^{m_{1}}J\left( 1\right) A^{-n_{1}},  \\
K\left( 2\right) &=A^{n_{1}+n_{2}}J\left( 1\right) ^{-1}B^{-m_{1}},  \\
J\left( 3\right) &=B^{m_{1}+m_{2}}J\left( 1\right) A^{-n_{1}-n_{2}},  \\
&\mathrel{\mkern4.5mu\vdots\mkern4.5mu}
\end{aligned}
\label{eqInt.3}%
\end{equation}
etc. If $n$ is a nonzero integer, let $\operatorname*{Prim}\left(  n\right)  $
denote the set of prime factors of $n$. Then (\ref{eqInt.2}) implies%
\begin{equation}
\operatorname*{Prim}\left(  \det\left(  A\right)  \right)
=\operatorname*{Prim}\left(  \det\left(  B\right)  \right)  , \label{eqInt.4}%
\end{equation}
and thus (\ref{eqInt.3}) implies%
\begin{equation}
\operatorname*{Prim}\left(  \det\left(  J\left(  1\right)  \right)  \right)
\subseteq\operatorname*{Prim}\left(  \det\left(  A\right)  \right)
=\operatorname*{Prim}\left(  \det\left(  B\right)  \right)  . \label{eqInt.5}%
\end{equation}

Thus a necessary and sufficient condition for $C^{\ast}$-equivalence of two
primitive, nonsingular $N\times N$ matrices $A$, $B$ with nonnegative integer
matrix entries, is the existence of a (necessarily nonsingular) matrix
$J\left(  1\right)  $ with nonnegative integer matrix entries and sequences
$n_{1},n_{2},\dots$ and $m_{1},m_{2},\dots$ of natural numbers such that the
matrices $K\left(  1\right)  ,J\left(  2\right)  ,\dots$ defined by
(\ref{eqInt.3}) have positive integer matrix entries.

Another way of formulating this is in terms of dimension groups (see
\cite{Bla86}, \cite{Eff81}, and \cite{BMT87} for details). Let $G\left(
A\right)  $ be the inductive limit of the sequence%
\begin{equation}
\mathbb{Z}^{N}\overset{A}{\longrightarrow}\mathbb{Z}^{N}\overset
{A}{\longrightarrow}\mathbb{Z}^{N}\longrightarrow\cdots\label{eqInt.6}%
\end{equation}
of free abelian groups with order generated by the order defined on each
$\mathbb{Z}^{N}$ by%
\begin{equation}
\left(  m_{1},\dots,m_{N}\right)  \geq0\iff m_{i}\geq0,\qquad i=1,\dots,N.
\label{eqInt.7}%
\end{equation}
Since we assume $\det A\neq0$, we may realize $G\left(  A\right)  $ concretely
as a subgroup of $\mathbb{Q}^{N}$ as follows: Put%
\begin{equation}
G_{n}\left(  A\right)  =A^{-n}\left(  \mathbb{Z}^{N}\right)  , \label{eqInt.8}%
\end{equation}
and equip $G_{n}\left(  A\right)  $ with the order%
\begin{equation}
\left(  G_{n}\right)  _{+}\left(  A\right)  =A^{-n}\left(  \left(
\mathbb{Z}_{+}\right)  ^{N}\right)  . \label{eqInt.9}%
\end{equation}
(Here and through the rest of the paper we use the dynamical-systems
convention that $\mathbb{Z}^{+}$ means the strictly positive integers and
$\mathbb{Z}_{+}$ the nonnegative integers, and correspondingly, $G^{+}$ means
the nonzero positive elements of $G$ and $G_{+}=G^{+}\cup\left\{  0\right\}
$.) Then $G_{0}\subseteq G_{1}\subseteq G_{2}\subseteq\cdots$ and we define%
\begin{equation}
G\left(  A\right)  =\bigcup_{n=0}^{\infty}G_{n} \label{eqInt.10}%
\end{equation}
with the order defined by%
\begin{equation}
g\geq0\text{\qquad if and only if }g\geq0\text{ in some }G_{n}\left(
A\right)  . \label{eqInt.11}%
\end{equation}
Then one fundamental characterization of $C^{\ast}$-equivalence is that there
exists a (necessarily nonsingular) matrix $J\left(  1\right)  $ in
$M_{N}\left(  \mathbb{Q}\right)  $ such that%
\begin{equation}
J\left(  1\right)  G\left(  A\right)  =G\left(  B\right)  \label{eqInt.12}%
\end{equation}
and%
\begin{equation}
J\left(  1\right)  G^{+}\left(  A\right)  =G^{+}\left(  B\right)  ;
\label{eqInt.13}%
\end{equation}
see \cite[Proposition 11.7]{BJO99}. The $1$--$1$ correspondence between group
isomorphism $\theta$ and matrix $J$ referred to in \cite{BJO99} is as follows:
If a matrix $J=J\left(  1\right)  $ is specified as above, then $\theta\colon
G\left(  A\right)  \rightarrow G\left(  B\right)  $, given by $\theta\left(
g\right)  =Jg$, $g\in G\left(  A\right)  $, will be an isomorphism. Here the
product $Jg$ is matrix multiplication, and each $g$ is viewed as a column
vector. Conversely, the observation in \cite{BJO99} is that every isomorphism
arises this way. This can also be formulated in other ways, as we shall
presently do.

If $A$ is a given primitive $N\times N$ matrix, let $\lambda_{\left(
A\right)  }$ denote its Perron--Frobenius eigenvalue, and let $v\left(
A\right)  $ denote a corresponding left (row) eigenvector with strictly
positive components and $w\left(  A\right)  $ a corresponding right (column)
eigenvector with strictly positive components, and in both cases use a
normalization such that the components are contained in the field
$\mathbb{Q}\left[  \lambda_{\left(  A\right)  }\right]  $. Define $V\left(
A\right)  $ as the orthogonal complement of $v\left(  A\right)  $, i.e.,
$V\left(  A\right)  =\left\{  x\in\mathbb{Q}\left[  \lambda_{\left(  A\right)
}\right]  ^{N}\mid%
\ip{v\left( A\right) }{x}%
=0\right\}  =v\left(  A\right)  ^{\perp}$. Then $V\left(  A\right)  $ is an
$\left(  N-1\right)  $-dimensional vector space of column vectors which will
sometimes be referred to, somewhat informally, as the linear span of the
nonmaximal generalized eigenvectors of $A$; see (\ref{eqIntNew.28}). Thus%
\begin{multline}
v\left(  A\right)  A=\lambda_{\left(  A\right)  }v\left(  A\right)
,\;Aw\left(  A\right)  =\lambda_{\left(  A\right)  }w\left(  A\right)  ,\text{
and}\label{eqInt.14}\\%
\ip{v\left( A\right) }{V\left( A\right) }%
=\left\{  0\right\}  ,\;%
\ip{v\left( A\right) }{w\left( A\right) }%
\in\mathbb{Q}\left[  \lambda_{\left(  A\right)  }\right]  \cap\left(
0,\infty\right)  .
\end{multline}
In particular, $A$ leaves $V\left(  A\right)  $ invariant, for if $u\in
V\left(  A\right)  $, then%
\begin{equation}%
\ip{v\left( A\right) }{Au}%
=%
\ip{v\left( A\right) A}{u}%
=\lambda_{\left(  A\right)  }%
\ip{v\left( A\right) }{u}%
=0, \label{eqIntNew.15}%
\end{equation}
and it follows that $Au\in V\left(  A\right)  $. The same argument applies to
the matrix $J$ from (\ref{eqInt.17}) below. It shows that any $J$ satisfying
(\ref{eqInt.17}) must map $V\left(  A\right)  $ onto $V\left(  B\right)  $;
i.e., $JV\left(  A\right)  =V\left(  B\right)  $. The number $%
\ip{v\left( A\right) }{w\left( A\right) }%
$ from (\ref{eqInt.14}) plays an important role in the discussion of the
isomorphism problem here (Section \ref{Cas}) and in \cite{BJO99}.

Let us mention an alternative form of the isomorphism criterion
(\ref{eqInt.12})--(\ref{eqInt.13}), formulated in \cite[Proposition
11.7]{BJO99}. Two primitive nonsingular $N\times N$ matrices $A$, $B$ with
positive integer matrix entries are $C^{\ast}$-equivalent if and only if there
is a nonsingular $N\times N$ matrix $J=J\left(  1\right)  $ in $M_{N}\left(
\mathbb{Z}\right)  $ satisfying the two conditions:%
\begin{equation}
v\left(  B\right)  J=\mu v\left(  A\right)  \text{\qquad for some }\mu
\in\left(  0,\infty\right)  , \label{eqInt.17}%
\end{equation}%
\begin{multline}
\text{for all }n\in\mathbb{Z}_{+}\text{, there is an }m\in\mathbb{Z}_{+}\text{
such that}\label{eqInt.18}\\
B^{m}JA^{-n}\text{ and }A^{m}J^{-1}B^{-n}\text{ both have integer matrix
entries;}%
\end{multline}
and then $J^{-1}$ has matrix entries in $\mathbb{Z}\left[  1/\det\left(
A\right)  \right]  =\mathbb{Z}\left[  1/\det\left(  B\right)  \right]  $. It
suffices to assume that $J\in\mathrm{\operatorname*{GL}}\left(  N,\mathbb{R}%
\right)  $, but then (\ref{eqInt.18}) forces $J$, $J^{-1}$ to lie in
$M_{N}\left(  \mathbb{Z}\left[  1/\det A\right]  \right)  =M_{N}\left(
\mathbb{Z}\left[  1/\det B\right]  \right)  $. So $J$ is not unique: one may,
for example, replace the given $J$ with $B^{m}JA^{-n}$ for any $m,n\in
\mathbb{Z}_{+}$. By choosing $m$ large enough, one may assure that $B^{m}J$
has integer matrix entries, and choosing it even larger one may also assure
that these entries are positive, and in fact (\ref{eqInt.17}) may be replaced
by the condition%
\begin{equation}
J\text{ has positive matrix entries.} \tag*{(\ref{eqInt.17})$^\prime$}%
\end{equation}
(But again, a given $J$ may satisfy (\ref{eqInt.17})--(\ref{eqInt.18}) without
having positive or integer matrix entries.) The combined two conditions
(\ref{eqInt.17}), (\ref{eqInt.18}) are equivalent to the two conditions
(\ref{eqInt.17})$^{\prime}$, (\ref{eqInt.18}), and to (\ref{eqInt.12}),
(\ref{eqInt.13}). For this one uses Perron--Frobenius theory (see, e.g.,
\cite{New72}): asymptotically when $m\rightarrow\infty$, $B^{m}$ behaves like
$\lambda_{\left(  B\right)  }^{m}$ times the projection onto $w\left(
B\right)  $, and $w\left(  B\right)  $ has strictly positive components.

In the two conditions (\ref{eqInt.17})--(\ref{eqInt.18}) on $J$, positivity of
the matrix entries is just hidden away in the first of the subconditions.
However, from (\ref{eqInt.1}), one may merge the two conditions into the joint
condition: There is a nonsingular $N\times N$ matrix $J=J\left(  1\right)  $
in $M_{N}\left(  \mathbb{Z}\right)  $ such that,%
\begin{multline}
\text{for all }n\in\mathbb{Z}_{+}\text{, there is an }m\in\mathbb{Z}_{+}\text{
such that}\label{eqInt.19}\\
B^{m}JA^{-n}\text{ and }A^{m}J^{-1}B^{-n}\text{ both have positive integer
matrix entries.}%
\end{multline}
Thus the single condition (\ref{eqInt.19}) is equivalent to each of the three
pairs of conditions (\ref{eqInt.12})--(\ref{eqInt.13}), (\ref{eqInt.17}%
)$^{\prime}$--(\ref{eqInt.18}), and (\ref{eqInt.17})--(\ref{eqInt.18}).

Let us record a fact which was not mentioned in \cite{BJO99}, namely that the
$m$ in (\ref{eqInt.19}) can be taken to depend linearly on $n$:

\begin{proposition}
\label{ProInt.1}Let $A$, $B$ be nonsingular primitive $N\times N$ matrices
with positive integer matrix entries, and assume that there is a nonsingular
matrix $J\in\mathrm{\operatorname*{GL}}\left(  N,\mathbb{R}\right)  $ such
that \textup{(\ref{eqInt.19})} holds. It follows that there exists a positive
integer $k$ and an integer $l$ such that%
\begin{multline}
\text{for all positive integers }n\text{, the matrices}\label{eqInt.20}\\
B^{kn+l}JA^{-n}\text{ and }A^{kn+l}J^{-1}B^{-n}\text{ both have positive
integer matrix entries.}%
\end{multline}
\end{proposition}

\begin{proof}
To show the existence of $k$, $l$ giving positivity we may modify the proof of
Theorem 6 in \cite{BJKR98} so as to make some specific estimates, i.e., we
show that if a solution to (\ref{eqInt.1}) exists, then the sequences $n_{i}$,
$m_{i}$ may be taken to grow at most linearly. Let $\lambda_{1}$, $\lambda
_{2}$ be the maximum eigenvalues of $B$, $A$. Let $\lambda_{3}<\lambda
_{1},\lambda_{2}$ exceed the largest absolute value of any other eigenvalue,
and let $\lambda_{4}$ be the largest absolute value of the reciprocal of any
eigenvalue. Consider $B^{m}JA^{-n}$. Using the above-mentioned (see
(\ref{eqInt.14})) two invariant complex vector-space (column vectors)
decompositions%
\begin{equation}
\mathbb{C}^{N}=V\left(  A\right)  \oplus\mathbb{C}w\left(  A\right)
\text{\quad and\quad}\mathbb{C}^{N}=V\left(  B\right)  \oplus\mathbb{C}%
w\left(  B\right)  , \label{eqIntNew.21}%
\end{equation}
we note that the contribution of the maximum eigenvector in $(JA^{-n})$ will
be at least $C\lambda_{2}^{-n}$ for some positive $C$. When we multiply it by
$B^{m}$ we get $\lambda_{1}^{m}C\lambda_{2}^{-n}$. The largest magnitude of
any other term will be some $\lambda_{3}^{m}C_{1}^{{}}\lambda_{4}^{n}$. We
want the former terms to dominate the sum of all the others, say to be $N^{2}$
times the largest, where $N$ is the dimension of the matrices. Take
logarithms, and we want
\begin{equation}
m\log\lambda_{1}+\log C-n\log\lambda_{2}>m\log\lambda_{3}+\log\left(
C_{1}N^{2}\right)  +n\log\lambda_{4} \label{ineq.1}%
\end{equation}
or rearranged equivalently as
\begin{equation}
m\left(  \log\lambda_{1}-\log\lambda_{3}\right)  >-\log C+\log\left(
C_{1}N^{2}\right)  +n\left(  \log\lambda_{2}+\log\lambda_{4}\right)  .
\label{ineq.2}%
\end{equation}
Then some arithmetic progression where the ratio of $m$ to $n$ exceeds
\begin{equation}
\frac{\log\lambda_{2}+\log\lambda_{4}}{\log\lambda_{1}-\log\lambda_{3}}
\label{ineq.3}%
\end{equation}
will give the domination.

Consider denominators in the matrix entries which have as divisor some
algebraic prime $p$. The prime $p$ is fixed, but we will do this for all prime
divisors in $\det\left(  A\right)  $. (For the definition of ``algebraic
prime'', see the end of Section \ref{Sub}.) For simplicity extend the
coefficient field and assume we can diagonalize the matrices (the case of a
standard Jordan form can be treated similarly). The maximum denominator in
$A^{-n}$ is $p^{-kn}$ for some constant $k$, which for instance can be worked
out from the determinant. Then consider the matrix entries in $B^{m}JA^{-n}$.
They will be sums of constants from the diagonalizing matrices times $m$
powers of the eigenvalues $\mu_{i}$ of $B$, i.e., $\sum_{i}c_{i}^{{}}\mu
_{i}^{m}$. The eigenvalues, when factored, only involve nonnegative powers of
$p$, since they are algebraic integers.

The terms in this sum for eigenvalues not divisible by $p$ must add up to be
an integer at the prime $p$: otherwise, no very large powers $m$ could make
the total an integer. For the other terms, as soon as $m$ exceeds $nk$ plus
the degrees of constants arising from diagonalization process, we will have
algebraic integers.
\end{proof}

There is another general observation about solving for $J$ and $K$ in
(\ref{eqInt.2}), with $A$ and $B$ given, which motivates the $p$-adic analysis
to follow and is a key point in the decidability argument. The identities
(\ref{eqInt.2}) are quadratic. Since the matrix entries on the left are all
integral, solving for $J$ and $K$ is therefore a quadratic diophantine problem
in the sense of \cite[Ch.~1]{BoSh66}: We thus have a system of quadratic
equations in the respective matrix entries of $J$ and $K$, and \cite[Theorem 1
on p.~61, Ch.~1, Section 7.1]{BoSh66} amounts to the assertion that the
solution to a quadratic diophantine problem is equivalent to instead solving a
finite system of related $p$-adic congruences, but for all $p$. (See the next
paragraph.) Hence, in the following, we will be stating criteria for $C^{\ast
}$-equivalence in terms of $p$-adic conditions. We will specify for which $p$
we need the conditions, for example in Corollary \ref{Corpad.2}, and we will
show that there are finite algorithms for deciding the problem.

The simplest case of solving an equation by congruences is to replace a
diophantine equation%
\begin{equation}
F\left(  x\right)  =\sum_{k_{1}\dots k_{n}}a_{k_{1}\dots k_{n}}x_{1}^{k_{1}%
}\cdots x_{n}^{k_{n}}=0 \label{eqIntNew.star}%
\end{equation}
where
\[
a_{k_{1}\dots k_{n}}\in\mathbb{Z}\text{,\qquad the sum is finite and }k_{i}%
\in\left\{  0,1,2,\dots\right\}  ,
\]
by the corresponding equation over the ring $\mathbb{Z}_{p^{k}}=\mathbb{Z}%
/p^{k}\mathbb{Z}$ of residues modulo $p^{k}$ for each prime $p$ and each
positive integer $k$. The latter problem amounts to checking only a finite
number of cases, since $\mathbb{Z}_{p^{k}}$ is obviously finite. The point
made in \cite{BoSh66} under the name of Hasse--Minkowski's Theorem is that
this is possible when (\ref{eqIntNew.star}) is quadratic, but not in general.
It is, for example, noted in \cite[p.~3]{BoSh66} that the congruences
\begin{equation}
\left(  x^{2}-13\right)  \left(  x^{2}-17\right)  \left(  x^{2}-221\right)
\equiv0\qquad\left(  \operatorname{mod}p^{k}\right)  \label{eqIntNew.25}%
\end{equation}
are solvable for all $p^{k}$, i.e., with solution $x\in\mathbb{Z}_{p^{k}}$,
while
\begin{equation}
\left(  x^{2}-13\right)  \left(  x^{2}-17\right)  \left(  x^{2}-221\right)  =0
\label{eqIntNew.26}%
\end{equation}
clearly has no solution $x\in\mathbb{Z}$. When matrices $A$ and $B$ are given
in $M_{N}\left(  \mathbb{Z}\right)  $, then solving equation (\ref{eqInt.2})
for $K$ and $J$ (in $M_{N}\left(  \mathbb{Z}\right)  $) is a quadratic
diophantine problem in the matrix entries of $K$ and $J$.

It seems to be difficult to convert Proposition \ref{ProInt.1} directly into
an effective decision procedure for isomorphism, since $J$ is not unique, and
hence it is difficult to obtain \emph{a priori} estimates on the norm of $J$
and on the coefficients $k$ and $l$. Instead we will turn to the completely
different method developed in \cite{KiRo88}, which is described in the
previous paragraph and in Section \ref{Dec}. Instead of starting with an
explicit norm estimate on $J$, we reduce the problem to a collection of
congruences and norm restrictions which are decidable by Lemma \ref{LemDec.1}.

The simple-minded way of trying to determine the dimension group from
(\ref{eqInt.8})--(\ref{eqInt.11}) is to take the algebraic extension of
$\mathbb{Q}$ determined by all the roots of the characteristic equation of
$A$, write $A$ in generalized Jordan form \cite{New72}, and then compute
$\bigcup_{m}A^{-m}\left(  \mathbb{Z}^{N}\right)  $ in the new basis. Vestiges
of this approach appear in our argument, but instead of using the complete
Jordan form we merely use a reduction to block-diagonal form where the blocks
correspond to generalized eigenspaces, first when we determine the subspaces
which a rational matrix $J_{0}$ has to preserve, and then in studying the
matrix giving the difference of the actual matrix $J$ from $J_{0}$.

If $A$ has the eigenvalues $\lambda_{1},\dots,\lambda_{n}$ in $\mathbb{C}$ and
we view $A$ as a map on the module $V=\mathbb{Q}\left[  \lambda_{1}%
,\dots,\lambda_{n}\right]  ^{N}$, then we have the standard direct-sum
decomposition%
\begin{equation}
V=\sum_{\lambda\in\operatorname*{Sp}\left(  A\right)  }V_{\lambda},
\label{eqIntNew.27}%
\end{equation}
where $V_{\lambda}$ is the \emph{generalized eigenspace}%
\begin{equation}
V_{\lambda}=\left\{  x\in V\mid\exists\,k\in\mathbb{Z}^{+}\ni\joinrel
\relbar\left(  \lambda-A\right)  ^{k}x=0\right\}  . \label{eqIntNew.28}%
\end{equation}
The elements of $V_{\lambda}$ are called the \emph{generalized eigenvectors}
corresponding to $\lambda$.

\section{\label{SY}Reduction to nonsingular matrices}

The first step in our decision procedure is to reduce the problem to the
corresponding problem for two matrices $A$, $B$ with integer coefficients
which are no longer positive, but nonsingular, and with a different definition
of positivity in $G\left(  A\right)  $ and $G\left(  B\right)  $. The
discussion of the general case when $A$ is not assumed nonsingular is resumed
in Section \ref{Red}.

The following result is well-known \cite{BMT87}. We include the proof since it
introduces terminology and details that we will need later.

\begin{lemma}
\label{LemSY.1}Every $N\times N$ matrix $A$ over $\mathbb{Z}$ is shift
equivalent to a nonsingular matrix over $\mathbb{Z}$. Specifically, let
$\mathcal{W}(A)=A^{N}\mathbb{Q}^{N}$. This is an $A$-invariant subspace, and
$A$ restricted to it is nonsingular. Let $\mathcal{W}_{0}=\mathcal{W}%
(A)\cap\mathbb{Z}^{N}$. This is also $A$-invariant. Choose a basis for this
free abelian group and express the restriction of $A$ to $\mathcal{W}_{0}$ as
a matrix $C$. Then $A$ is shift equivalent to $C$.
\end{lemma}

\begin{proof}
There is an $M\leq N$ such that $\mathcal{W}_{0}\cong\mathbb{Z}^{M}$, and let
$x_{1},\dots,x_{M}$ be a column vector basis for the free abelian group
$\mathcal{W}_{0}\subseteq\mathbb{Z}^{N}$. Define a matrix $R$ by%
\[
R=\left[  x_{1},\dots,x_{M}\right]  \colon\mathbb{Z}^{M}\longrightarrow
\mathcal{W}_{0}\subseteq\mathbb{Z}^{N}.
\]
The $N\times M$ matrix $R$ maps $\mathbb{Z}^{M}$ bijectively onto
$\mathcal{W}_{0}\subseteq\mathbb{Z}^{N}$. Let $R^{-1}$ be an $M\times N$
matrix extending the inverse of this bijective map. Then
\[
C=R^{-1}AR\colon\mathbb{Z}^{M}\longrightarrow\mathbb{Z}^{M}%
\]
is an $M\times M$ matrix. The mapping $A^{N}$ on $\mathbb{Z}^{N}$ has image in
$\mathcal{W}_{0}$, and hence we may define an $M\times N$ matrix map $S$ by
\[
S=R^{-1}A^{N}\colon\mathbb{Z}^{N}\longrightarrow\mathbb{Z}^{M}.
\]
Now one immediately verifies that $AR=RC$, $SA=CS$, $SR=C^{N}$, $RS=A^{N}$,
the equations of shift equivalence. Over an extension field, we may
triangularize $A$, and find that $\mathcal{W}_{0}$ over this field becomes the
sum of all nonzero generalized eigenspaces, since on them $A^{N}$ is an
isomorphism, but on the zero generalized eigenspace of dimension at most $N$,
$A^{N}$ is zero. Therefore $C$ is nonsingular.
\end{proof}

This mapping $A\mapsto C$ preserves unordered dimension groups. The lemma also
shows that every unordered dimension group of an integer matrix (see
(\ref{eqRedNew.1}) below) is the unordered dimension group of a nonsingular
integer matrix. Moreover, the order structure is given by evaluation on a
Perron-Frobenius row eigenvector. To see this, note that with the notation of
the proof of Lemma \ref{LemSY.1}, if $v$ is a row vector in $\mathbb{R}^{N}$
such that%
\[
vA=\lambda v,
\]
then it follows that%
\[
vRC=vAR=\lambda vR,
\]
i.e., $vR$ is a row vector in $\mathbb{R}^{M}$ which is an eigenvector of $C$
with the same eigenvalue $\lambda$. Applying this on the Perron--Frobenius
eigenvalue $\lambda$ and corresponding eigenvector $v$, we see that the order
structure of the dimension group, realized in $\mathbb{Z}^{M}$, is given by
evaluation on $vR$. This evaluation can be considered as the projection on
column vectors over an extension field which is the identity on the
Perron-Frobenius column eigenvector and is zero on all other generalized
eigenspaces. Note that the nonsingular matrix $C$ has the property that it has
a positive eigenvalue $\lambda$ which has strictly larger modulus than any
other eigenvalue, and if $G\left(  C\right)  =\bigcup_{n}C^{-n}\mathbb{Z}^{M}$
is the corresponding dimension group, then $G\left(  C\right)  $ is isomorphic
to $G\left(  A\right)  $ as a group by an isomorphism taking strictly positive
elements in $G\left(  A\right)  $ into elements $g^{\prime}\in\bigcup
_{n}C^{-n}\mathbb{Z}^{N}$ such that $%
\ip{vR}{g^{\prime}}%
>0$. In this way we consider the ordered dimension groups of a singular matrix
as the unordered dimension group of a nonsingular matrix together with an
order structure which amounts to preservation of the sum of
non-Perron--Frobenius generalized eigenspaces. If we start out working with a
singular matrix, then we will replace it with this nonsingular matrix as we
continue with the decision procedure in Section \ref{Dec}.

We will not need Theorems 6 and 7 from \cite{BJKR98} directly in Section
\ref{Dec}, as we may work directly with the dimension groups defined from $C$
as above, but note that these theorems could also be generalized to the
setting of primitive singular matrices as follows.

The argument of Theorem 7 of our first paper applies over $\mathbb{Z}$ to
characterize maps giving isomorphism of unordered dimension groups. The
argument of Theorem 6 will also apply over $\mathbb{Z}$; however, there is a
problem with the inverses which are used in the proof. Let $M^{\langle
-1\rangle}$ denote the Drazin inverse of an $N\times N$ matrix $M$. The Drazin
inverse is a matrix having the row and column space of $M^{N}$ and is such
that $MM^{\langle-1\rangle}=M^{\langle-1\rangle}M$ is a projection to this
column space (is the identity on it). The Drazin inverse is multiplicative
over matrices having the same eventual row and column space and is unique, and
agrees with the ordinary inverse for nonsingular matrices. In effect, it is
the ordinary inverse restricted to the nonzero generalized eigenspaces. Let
$D$ be the determinant of $M$ restricted to the nonzero generalized
eigenspaces, the product of its nonzero eigenvalues. Then the Drazin inverse
has denominators which divide $D^{N}$. Use $M^{\left\langle -n\right\rangle }$
to denote the Drazin inverse of $M^{n}$.

\begin{lemma}
\label{LemSY.2}For singular matrices $A$, $B$ over $\mathbb{Z}$, a necessary
and sufficient condition for $C^{\ast}$-equivalence of matrices $A$, $B$ is
the existence of a nonnegative matrix $J(1)$ such that for all $d\in
\mathbb{Z}_{+}$ there is a $c\in\mathbb{Z}_{+}$ such that%
\[
B^{c}J(1)A^{\left\langle -d\right\rangle }%
\]
and%
\[
A^{c}J(1)^{-1}B^{\left\langle -d\right\rangle }%
\]
are nonnegative integer matrices.
\end{lemma}

\begin{proof}
Form the matrices $J\left(  i\right)  $, $K\left(  i\right)  $, using the
formulas (\ref{eqInt.3}), but now with the respective Drazin inverses in place
of inverses. In doing so, we can multiply by suitably chosen higher powers of
$A$ or $B$, and thereby force the $J(i)$'s and $K(i)$'s to have their row and
column spaces contained in those of $A^{N}$, and then we must have the same
inverse formulas for them as in the nonsingular case, provided that we use the
Drazin inverse.
\end{proof}

With this change, the proofs of Theorems 6 and 7 in \cite{BJKR98} will also go
over to the singular case, and characterize $C^{\ast}$-isomorphism of matrices.

\section{\label{Sub}Subspace structure and localization}

In the proof of the main theorem, Theorem \ref{ThmDec.7}, the structure of
subspaces of $\mathbb{C}^{N}$ which are mapped into each other by a possible
intertwiner matrix $J\in M_{N}\left(  \mathbb{Z}\right)  $ will be important.
One general idea is the following: Consider a certain subset $D_{A}$ of
$G\left(  A\right)  $ which is defined by a property which is invariant under
group isomorphism. Then%
\begin{align}
\tilde{D}_{A}  &  =\left\{  g\in G\mid\exists\,n,n_{1},\dots,n_{m}%
\in\mathbb{Z},\;g_{1},\dots,g_{m}\in D_{A}\ni\joinrel\relbar ng=%
{\textstyle\sum\nolimits_{i=1}^{m}}
n_{i}g_{i}\right\} \label{eqSub.1}\\
&  =\text{``the subgroup of }G\text{ linearly spanned by }D_{A}\text{''}%
\nonumber
\end{align}
is a subgroup of $G$. If $G\left(  B\right)  $ is another subgroup of
$\mathbb{C}^{N}$, and $G\left(  B\right)  =JG\left(  A\right)  $, and $D_{B}$
and $\tilde{D}_{B}$ are defined as above for $G\left(  B\right)  $, we must
have%
\begin{equation}
D_{B}=JD_{A},\qquad\tilde{D}_{B}=J\tilde{D}_{A}, \label{eqSub.2}%
\end{equation}
and hence%
\begin{equation}
J\mathbb{R}D_{A}=\mathbb{R}D_{B}. \label{eqSub.3}%
\end{equation}
This idea was much exploited in \cite{BJO99} on the subgroups%
\begin{align}
D_{m}\left(  G\left(  A\right)  \right)   &  =\bigcap\nolimits_{i}%
m^{i}G\left(  A\right) \label{eqSub.4}\\
&  =\text{the set of elements of }G\left(  A\right)  \text{ which are
infinitely divisible by }m,\nonumber
\end{align}
and we will soon give an example of this in a more general setting than in
\cite{BJO99}. Note in particular that if $m$ is a rational eigenvalue of $A$,
then $m$ is an integer since the characteristic equation of $A$ is monic, and
hence $D_{m}\left(  G\left(  A\right)  \right)  $ is nonzero and gives
nontrivial information about $J$. We would like to exploit this idea also when
$\lambda$ is an irrational eigenvalue of $A$, but since $G\left(  A\right)
\subset\mathbb{Z}\left[  1/\left|  \det A\right|  \right]  ^{N}$, $G\left(
A\right)  $ then clearly does not contain eigenvectors of $A$. To remedy this
situation, we may augment or localize $G\left(  A\right)  $ and $G\left(
B\right)  $ by equipping them with coefficients outside $\mathbb{Z}$, i.e., by
considering tensor products%
\begin{equation}
\tilde{G}\left(  A\right)  =E\otimes G\left(  A\right)  ,\qquad\tilde
{G}\left(  B\right)  =E\otimes G\left(  B\right)  , \label{eqSub.5}%
\end{equation}
where $E$ is any $\mathbb{Z}$-module, and then $J$ still defines an
isomorphism between $\tilde{G}\left(  A\right)  $ and $\tilde{G}\left(
B\right)  $. One then tries to choose $E$ to optimize the information about
subspaces. In \cite{BJO99} this remedy was used with $E$ finite cyclic groups,
but one may use $p$-adic numbers, or, as we will also do, various finite
algebraic extensions of $\mathbb{Z}$. Which extension is used has to be
fine-tuned to the problem. For example, if $E=\mathbb{Z}\left[  1/\left|  \det
A\right|  \right]  $, then $\tilde{G}\left(  A\right)  =\mathbb{Z}\left[
1/\left|  \det A\right|  \right]  ^{N}$, and all information about $G\left(
A\right)  $ disappears (except for its rank and the prime factors of $\left|
\det A\right|  $, which both are invariants). Similarly, if $\lambda$ is an
algebraic integer which is a unit, i.e., is such that the constant term in its
minimal polynomial is $\pm1$, then $\lambda^{-1}\in\mathbb{Z}\left[
\lambda\right]  $, and hence \emph{all} elements of $\mathbb{Z}\left[
\lambda\right]  \otimes G\left(  A\right)  $ are divisible by $\lambda$, and
no information on the subspace structure is obtained. One useful choice of $E$
is based on Theorem 10 in \cite{BJKR98}: If $G\left(  A\right)  $ and
$G\left(  B\right)  $ are isomorphic and $\lambda_{\left(  A\right)  }$ and
$\lambda_{\left(  B\right)  }$ are the Perron--Frobenius eigenvalues of $A$
and $B$, then the fields $\mathbb{Q}\left[  \lambda_{\left(  A\right)
}\right]  $ and $\mathbb{Q}\left[  \lambda_{\left(  B\right)  }\right]  $ are
the same, and $\lambda_{\left(  A\right)  }$ and $\lambda_{\left(  B\right)
}$ are the products of the same primes over this field. A prime in this
context means a prime ideal in the associated subring $\mathbb{O}\left[
\lambda\right]  $ of algebraic integers, i.e., $\mathbb{O}\left[
\lambda\right]  $ is the ring of all elements of $\mathbb{Q}\left[
\lambda\right]  $ which satisfy equations in monic polynomials over
$\mathbb{Z}$, so that%
\begin{equation}
\mathbb{Z}\left[  \lambda\right]  \subset\mathbb{O}\left[  \lambda\right]
\subset\mathbb{Q}\left[  \lambda\right]  . \label{eqSub.6}%
\end{equation}
Recall that an ideal $\mathcal{I}$ in a ring is a prime ideal if whenever
$\mathcal{I}=\mathcal{I}_{1}\mathcal{I}_{2}$ for two ideals $\mathcal{I}_{1}$,
$\mathcal{I}_{2}$, then $\mathcal{I}=\mathcal{I}_{1}$ or $\mathcal{I}%
=\mathcal{I}_{2}$. One useful choice for $E$ is thus $\mathbb{O}\left[
\lambda_{\left(  A\right)  }\right]  =\mathbb{O}\left[  \lambda_{\left(
B\right)  }\right]  $. One other choice we shall use is%
\begin{equation}
\Omega=\mathbb{O}\left[  \lambda_{1},\dots,\lambda_{N},\mu_{1},\dots,\mu
_{N}\right]  , \label{eqSub.7}%
\end{equation}
where $\lambda_{1},\dots,\lambda_{N},\mu_{1},\dots,\mu_{N}$ are the respective
roots in $\mathbb{C}$ of the characteristic equations of $A$ and $B$:
\begin{equation}
\det\left(  \lambda\openone-A\right)  =0,\qquad\det\left(  \mu\openone
-B\right)  =0. \label{eqSub.8}%
\end{equation}

Our decision procedure will involve even other rings and fields. For example,
in Section \ref{pad} below, we will consider $\mathbb{Z}_{\left(  p\right)
}\subset\mathbb{Q}_{\left(  p\right)  }$, i.e., the $p$-adic integers and the
$p$-adic numbers. More generally, the setting is $\mathcal{O}\subset
\mathcal{F}$ where $\mathcal{F}$ is an algebraic number field, and
$\mathcal{O}$ denotes the algebraic integers in $\mathcal{F}$. This is
explained in \cite[Chapter 1]{Wei98}. There prime divisors are defined as
equivalence classes of valuations, and the terminology is calibrated in such a
way that a compactness argument shows that prime ideals are of the form
$\pi\mathcal{O}$ for suitable elements $\pi\in\mathcal{O}$. The element $\pi$
is associated with a given valuation $\varphi$ by requiring that
$\varphi\left(  \pi\right)  $ assumes the maximal value $<1$ taken on by
$\varphi$.

Another example of the setup $\mathcal{O}\subset\mathcal{F}$ is $\mathcal{F}%
=\mathbb{Q}\left(  \sqrt{-5}\right)  $ and $\mathcal{O}=\mathbb{Z}%
+\mathbb{Z}\sqrt{-5}$. In general we have unique factorization in terms of
prime ideals, but the examples%
\begin{equation}
6=2\cdot3=\left(  1+\sqrt{-5}\right)  \left(  1-\sqrt{-5}\right)
\label{eqSubNew.9}%
\end{equation}
and%
\begin{equation}
14=7\cdot2=\left(  3+\sqrt{-5}\right)  \left(  3-\sqrt{-5}\right)
\label{eqSubNew.10}%
\end{equation}
show that we may have nonunique factorization in terms of irreducibles.

In general, the various field extensions are independent of one another, but
there are still some embeddings (perhaps unexpected) which will be used in our
analysis of $p$-adic eventual row spaces in Section \ref{pad}. (If $p$ is a
prime number, let $\mathbb{Z}_{\left(  p\right)  }$ and $\mathbb{Q}_{\left(
p\right)  }$ denote the $p$-adic integers and $p$-adic numbers, respectively;
see the first paragraph of Section \ref{pad} for definitions.) These field
extensions may be identified by use of a Newton approximation scheme; see,
e.g., \cite[Chapter 1, Section 5]{BoSh66} and \cite[Chapter 5]{Wei98}. For
example, the field $\mathbb{Q}\left(  \sqrt{-5}\right)  $ of (\ref{eqSubNew.9}%
)--(\ref{eqSubNew.10}) is embedded in $\mathbb{Q}_{\left(  3\right)  }$ and in
$\mathbb{Q}_{\left(  7\right)  }$, but not in $\mathbb{Q}_{\left(  11\right)
}$. This is because the equation $x^{2}+5=0$ has solutions in $\mathbb{Z}%
_{\left(  3\right)  }$ and $\mathbb{Z}_{\left(  7\right)  }$, but not in
$\mathbb{Z}_{\left(  11\right)  }$. (The polynomial $x^{2}+5$ is irreducible
in $\mathbb{Z}_{\left(  11\right)  }\left[  x\right]  $.) All the extensions
$\mathbb{Q}_{\left(  3\right)  }$, $\mathbb{Q}_{\left(  7\right)  }$, and
$\mathbb{Q}_{\left(  11\right)  }$ are, however, mutually independent; see
\cite[Section 1-2]{Wei98}.

To make it clear when our primes refer to those in the standard setup
$\mathbb{Z}\subset\mathbb{Q}$, i.e., when the primes are just $2,3,5,7,\dots$,
we refer to the latter as ``rational primes''; but if there is no danger of
confusion, we will simply refer to them as primes. Recall that ``algebraic
prime'' means a prime in the associated subring of algebraic integers.

We will also work with Galois field extensions $\mathbb{Q}\subset\mathcal{F}$
(see \cite{Rot98}); for example $\mathcal{F}$ may be obtained by adjoining
roots to $\mathbb{Q}$. As usual, the Galois group is defined as the group of
automorphisms of $\mathcal{F}$ leaving $\mathbb{Q}$ pointwise fixed; thus,
elements in the Galois group permute the roots and are uniquely determined by
this permutation. The Galois group will act on vectors over $\mathcal{F}$ by
$\left(  x_{i}\right)  \rightarrow\left(  x_{i}^{g}\right)  $, where $x^{g}$
for $x\in\mathcal{F}$ and $g\in\mathcal{G}$ denotes the action of $g$ on $x$.
Hence $\mathcal{G}$ also acts on matrices over $\mathcal{F}$; and, either way,
the respective actions will be used in defining Galois conjugacy.

\section{\label{pad}$p$-adic characterization of $J$}

We have already given several characterizations of the intertwiner $J$ more or
less in terms of the dimension groups $G\left(  A\right)  $, $G\left(
B\right)  $, i.e., (\ref{eqInt.1}), ((\ref{eqInt.17}) \& (\ref{eqInt.18})),
((\ref{eqInt.17})$^{\prime}$ \& (\ref{eqInt.18})), and (\ref{eqInt.19}). Here
$G\left(  A\right)  $ and $G\left(  B\right)  $ are defined in terms of
asymptotic properties of $A^{-n}$ and $B^{-n}$ as $n\rightarrow\infty$. We
will now give an exposition of another property of $J$ given in terms of
asymptotic properties of the \emph{positive} powers $A^{n}$ and $B^{n}$ as
$n\rightarrow\infty$. Since $n\rightarrow A^{n}\mathbb{Z}^{N}$ is decreasing,
and
\begin{equation}
\bigcap_{n}A^{n}\mathbb{Z}^{N}=\left\{  m\in\mathbb{Z}^{N}\mid q\left(
A\right)  m=0\right\}  \label{eqpad.1}%
\end{equation}
by \cite[Proposition 12.1]{BJO99}, where $q\left(  t\right)  $ is the product
of those irreducible (over $\mathbb{Z}$) factors of $\det\left(
t\openone-A\right)  $ which have constant term $\pm1$, the lattices
$\bigcap_{n}A^{n}\mathbb{Z}^{N}$ give very little information except that $J$
has to map $\bigcap_{n}A^{n}\mathbb{Z}^{N}$ onto $\bigcap_{n}B^{n}%
\mathbb{Z}^{N}$. However, if one replaces these intersections by $p$-adic
limits, one can say much more. Recall that if $p\in\left\{  2,3,5,7,11,\dots
\right\}  $ is an ordinary prime, the ring of $p$-adic integers $\mathbb{Z}%
_{\left(  p\right)  }$ is the projective limit%
\begin{equation}
0\overset{p}{\longleftarrow}\mathbb{Z}_{p}\overset{p}{\longleftarrow
}\mathbb{Z}_{p^{2}}\overset{p}{\longleftarrow}\mathbb{Z}_{p^{3}}%
\longleftarrow\dots\longleftarrow\mathbb{Z}_{\left(  p\right)  },
\label{eqpad.2}%
\end{equation}
where the left maps are multiplication by $p$. It can be equipped with a
topology making it into a compact totally disconnected ring. This is in fact
the topology the additive group $\mathbb{Z}_{\left(  p\right)  }$ has as a
dual group to $\mathbb{Z}_{p^{\infty}}$ viewed as the inductive limit of the
discrete groups%
\begin{equation}
0\hooklongrightarrow\mathbb{Z}_{p}\hooklongrightarrow\mathbb{Z}_{p^{2}%
}\hooklongrightarrow\mathbb{Z}_{p^{3}}\hooklongrightarrow\dots
\hooklongrightarrow\mathbb{Z}_{p^{\infty}}, \label{eqpad.3}%
\end{equation}
where the injections come from the standard realization of $\mathbb{Z}%
_{p^{\infty}}=\mathbb{Z}\left[  1/p\right]  /\mathbb{Z}$ as a subgroup of the
circle group $\mathbb{T}$; see \cite{Kob84}, \cite{Ser79}, \cite{Ser98}.
Koblitz uses the terminology $\mathbb{Z}_{p}$ for the $p$-adic integers, our
$\mathbb{Z}_{\left(  p\right)  }$, while we reserve $\mathbb{Z}_{p}$ for
$\mathbb{Z}/p\mathbb{Z}$. Other authors, e.g., \cite{BoSh66}, use $O_{p}$ for
the $p$-adic integers. In the duality consideration of the two groups
$\mathbb{Z}_{\left(  p\right)  }$ and $\mathbb{Z}_{p^{\infty}}$ of
(\ref{eqpad.2})--(\ref{eqpad.3}), we use the duality notion of locally compact
abelian \emph{groups,} e.g., $\mathbb{Z}_{p^{\infty}}$ is realized as the
group of continuous characters on $\mathbb{Z}_{\left(  p\right)  }$, and
conversely, $\mathbb{Z}_{\left(  p\right)  }$ acts as the group of all
characters on $\mathbb{Z}_{p^{\infty}}$. Now $\mathbb{Z}_{\left(  p\right)  }$
is a ring and thus a $\mathbb{Z}$-module, but it is not a field: If $q$ is an
integer, then $1/q\in\mathbb{Z}_{\left(  p\right)  }$ if and only if $q$ is
mutually prime with $p$. However, $\mathbb{Z}\left[  1/p\right]
\otimes\mathbb{Z}_{\left(  p\right)  }=\mathbb{Q}_{\left(  p\right)  }$ is a
field called the $p$-adic numbers.

Now if $A\in M_{N}\left(  \mathbb{Z}\right)  $ is a matrix, we may view $A$ as
a matrix with matrix entries in $\mathbb{Z}_{\left(  p\right)  }$, and we may
associate a unique idempotent%
\begin{equation}
E_{\left(  p\right)  }\left(  A\right)  =E\left(  A\right)  \in M_{N}\left(
\mathbb{Z}_{\left(  p\right)  }\right)  \label{eqpad.4}%
\end{equation}
with $A$ by using the following presumably known lemma (we did not find a reference).

\begin{lemma}
\label{Lempad.1}If $A\in M_{N}\left(  \mathbb{Z}/q\mathbb{Z}\right)  $ for a
$q\in\mathbb{Z}$, then the semigroup $\left\{  A,A^{2},A^{3},\dots\right\}  $
contains an idempotent $E$. This idempotent is unique, and $\left\{  n\mid
A^{n}=E\right\}  $ is a subsemigroup of $\mathbb{Z}^{+}$.
\end{lemma}

\begin{proof}
Since $M_{N}\left(  \mathbb{Z}/q\mathbb{Z}\right)  $ is finite, there is an
$m_{0}\in\mathbb{Z}_{+}$ and an $n_{0}\in\mathbb{Z}^{+}$ such that
$A^{n_{0}+m_{0}}=A^{m_{0}}$.But then $A^{n_{0}+m}=A^{m}$ for all $m\geq m_{0}$
and thus $A^{kn_{0}+m}=A^{m}$ for all $k\in\mathbb{Z}^{+}$. Choose $k$ such
that $kn_{0}\geq m_{0}$ and put $m=kn_{0}$. This gives $\left(  A^{kn_{0}%
}\right)  ^{2}=A^{kn_{0}}$ so $A^{kn_{0}}$ is idempotent.

If $A^{n}$ and $A^{m}$ are idempotents, then $A^{n}=\left(  A^{n}\right)
^{m}=\left(  A^{m}\right)  ^{n}=A^{m}$, so the idempotent is unique. If it is
called $E$, then if $A^{n}=A^{m}=E$, then $A^{n+m}=E\cdot E=E$, so $\left\{
n\mid A^{n}=E\right\}  $ is a semigroup.
\end{proof}

We now turn to part of the construction of the idempotent $E_{\left(
p\right)  }\left(  A\right)  $ in (\ref{eqpad.4}).

Fix a prime $p$, and let $e\left(  m\right)  $ be an increasing sequence of
integers such that $A^{e\left(  m\right)  }$ is an idempotent modulo $p^{m}$
in $M_{N}$,%
\begin{equation}
\left(  A^{e\left(  m\right)  }\right)  ^{2}=A^{e\left(  m\right)  }%
\mod{p^mM_N\left(  \mathbb{Z}\right)  }. \label{eqpad.5}%
\end{equation}
This sequence exists because of Lemma \ref{Lempad.1}, and by thinning out the
sequence, and using Lemma \ref{Lempad.1} again, we may also assume%
\begin{equation}
\left(  B^{e\left(  m\right)  }\right)  ^{2}=B^{e\left(  m\right)  }%
\mod{p^mM_N\left(  \mathbb{Z}\right)  }. \label{eqpad.6}%
\end{equation}
But by the uniqueness of the idempotent, it follows that%
\begin{equation}
m^{\prime}>m\Longrightarrow A^{e\left(  m^{\prime}\right)  }=A^{e\left(
m\right)  }\mod{p^m}, \label{eqpad.6bis}%
\end{equation}
and hence, by passing to yet another subsequence,%
\begin{equation}
E_{\left(  p\right)  }\left(  A\right)  =\lim_{m\rightarrow\infty}A^{e\left(
m\right)  } \label{eqpad.7}%
\end{equation}
exists in $M_{N}\left(  \mathbb{Z}_{\left(  p\right)  }\right)  $, and
$E_{\left(  p\right)  }\left(  A\right)  $ is an idempotent matrix.
Correspondingly, $E_{\left(  p\right)  }\left(  B\right)  $ is an idempotent
matrix. Now, if $A$ and $B$ define isomorphic dimension groups $G\left(
A\right)  $ and $G\left(  B\right)  $, it follows from (\ref{eqInt.18}) that
there exist for each $n\in\mathbb{Z}_{+}$ integer matrices $K_{n},L_{n}\in
M_{N}\left(  \mathbb{Z}\right)  $ and positive integers $m_{n}$such that%
\begin{align}
B^{m_{n}}J  &  =K_{n}A^{n},\label{eqpad.8}\\
A^{m_{n}}  &  =L_{n}B^{n}J. \label{eqpad.9}%
\end{align}
We may replace the powers $m_{n}$ by a new sequence (and thus $K_{n}$, $L_{n}$
by new integer matrices) to ensure that $A^{m_{n}}$, $B^{m_{n}}$ have
subsequences converging $p$-adically to the idempotents $E_{\left(  p\right)
}\left(  A\right)  $ and $E_{\left(  p\right)  }\left(  B\right)  $. Since
$\mathbb{Z}_{\left(  p\right)  }$ is compact (and metrizable), it follows that
there is a subsequence of $n\rightarrow\infty$ such that $\lim_{n}K_{n}=K$ and
$\lim_{n}L_{n}=L$ exist in $M_{N}\left(  \mathbb{Z}_{\left(  p\right)
}\right)  $, and we get from the relations above that%
\begin{align}
E_{\left(  p\right)  }\left(  B\right)  J  &  =KE_{\left(  p\right)  }\left(
A\right)  ,\label{eqpad.10}\\
E_{\left(  p\right)  }\left(  A\right)   &  =LE_{\left(  p\right)  }\left(
B\right)  J. \label{eqpad.11}%
\end{align}
Now define the $\mathbb{Z}_{\left(  p\right)  }$-eventual row space
$G_{\left(  p\right)  }\left(  A\right)  $ of $A$ as the linear combinations
over $\mathbb{Z}_{\left(  p\right)  }$ of the row-vectors of $E_{\left(
p\right)  }\left(  A\right)  $, and similarly for $E_{\left(  p\right)
}\left(  B\right)  $. Then (\ref{eqpad.10}) and (\ref{eqpad.11}) together say
that%
\begin{equation}
G_{\left(  p\right)  }\left(  B\right)  J=G_{\left(  p\right)  }\left(
A\right)  . \label{eqpad.12}%
\end{equation}
Thus (\ref{eqpad.12}) holds for any prime $p$. But conversely, by taking
$p$-adic limits as in the proof of Theorem 7 in \cite{BJKR98}, if
(\ref{eqpad.12}) holds for all primes $p$ in the set $\operatorname*{Prim}%
\left(  \det\left(  A\right)  \right)  =\operatorname*{Prim}\left(
\det\left(  B\right)  \right)  $, then we can recover (\ref{eqInt.18}). Thus
\begin{equation}
G_{\left(  p\right)  }\left(  B\right)  J=G_{\left(  p\right)  }\left(
A\right)  \text{\qquad for all }p\in\operatorname*{Prim}\left(  \det\left(
A\right)  \right)  =\operatorname*{Prim}\left(  \det\left(  B\right)  \right)
\tag*{(\ref{eqpad.12})$^\prime$}%
\end{equation}
is equivalent to (\ref{eqInt.18}) (the equivalence of (\ref{eqInt.18}) and
(\ref{eqpad.12})$^{\prime}$ is Theorem 7 in \cite{BJKR98}). The details
supplied above expand on the arguments from \cite{BJKR98}, which were somewhat
terse. Let us cast Theorem 7 in \cite{BJKR98} in a somewhat different, but
equivalent, form:

\begin{corollary}
\label{Corpad.2}In order that the unordered dimension groups $\bigcup
_{n}A^{-n}\mathbb{Z}^{N}$ and \linebreak $\bigcup_{n}B^{-n}\mathbb{Z}^{N}$
associated with a pair of nonsingular matrices $A$, $B$ be isomorphic, it is
necessary and sufficient that $\operatorname*{Prim}\left(  \det\left(
A\right)  \right)  =\operatorname*{Prim}\left(  \det\left(  B\right)  \right)
$, and that there exists a nonsingular matrix $J\in\mathrm{\operatorname*{GL}%
}\left(  N,\mathbb{Z}\left[  1/\det\left(  A\right)  \right]  \right)  $
\textup{(}i.e., the matrix entries of $J$ are in $\mathbb{Z}\left[
1/\det\left(  A\right)  \right]  $ and $\det\left(  J\right)  $ is invertible
in the ring $\mathbb{Z}\left[  1/\det\left(  A\right)  \right]  $\textup{)}
such that%
\begin{equation}
G_{\left(  p\right)  }\left(  B\right)  J=G_{\left(  p\right)  }\left(
A\right)  \label{eqCorpad.2}%
\end{equation}
for each prime $p\in\operatorname*{Prim}\left(  \det\left(  A\right)  \right)
$.
\end{corollary}

What makes this particularly useful for the decidability problem is that any
countably generated torsion-free module over the $p$-adic integers has a
trivial structure: such a module is merely a direct sum of replicas of the
$p$-adic numbers or the $p$-adic integers (\cite{Pru25}; see also
\cite{KaMa51}). The total number of direct summands in $G_{\left(  p\right)
}\left(  B\right)  $ and $G_{\left(  p\right)  }\left(  A\right)  $ is bounded
by the rank $N$ of $A$ or $B$. This makes it possible to decide whether or not
$J$ exists with the property (\ref{eqpad.12}) for each $p$, but the remaining
problem is to find a joint $J$ for all $p$ in $\operatorname*{Prim}\left(
A\right)  $ and to ensure the positivity property (\ref{eqInt.17}). Note that
in our setting we have $G_{\left(  p\right)  }\left(  A\right)  \subseteq
\mathbb{Z}_{\left(  p\right)  }^{N}$ by construction as $p$-adic limits of
integer vectors, and hence $G_{\left(  p\right)  }\left(  A\right)  $ cannot
contain any element which is infinitely divisible by $p$, and thus $G_{\left(
p\right)  }\left(  A\right)  $ as a $\mathbb{Z}_{\left(  p\right)  }$-module
is just a direct sum of at most $N$ copies of $\mathbb{Z}_{\left(  p\right)
}$ (no direct summand $\mathbb{Q}_{\left(  p\right)  }$ can occur). However,
be warned, since $\mathbb{Z}_{\left(  p\right)  }$ is not a field, this is not
as useful as knowing that a vector space (over a field) has a certain
dimension, since the usual operations of change of basis, etc., cannot be
performed within the ring $\mathbb{Z}_{\left(  p\right)  }$. In particular,
(\ref{eqCorpad.2}) says much more than that the $p$-adic row spaces have the
same rank.

\begin{remark}
\label{RempadNew.3}To see that the $p$-adic idempotents and row spaces are
independent of the chosen subsequences, note more generally that, when an
algebraic prime $\pi$ is given, we may determine which eigenvalues of $A$ are
divisible by $\pi$. The Newton polygon \cite[Section 3-1, pp.~73--78]{Wei98}
for the characteristic polynomial helps to tell which eigenvalues can be taken
as $\pi$-divisible for algebraic primes $\pi$. Then diagonalize $A$, and
replace the $\pi$-divisible eigenvalues by $0$ and other eigenvalues by $1$,
to get the projection operator $E_{\left(  \pi\right)  }\left(  A\right)  $
onto the eventual $\pi$-adic row space. If $\pi= p$ is a rational prime,
$E_{(\pi)}(A) = E _{(p)}(A)$ is the projection defined by
\textup{(\ref{eqpad.7}).} In the case that $\pi$ is a nonrational algebraic
prime, the procedure above gives a working man's definition of $E_{(\pi)}(A)$.
Strictly speaking, the idempotents $E_{\left(  \pi\right)  }\left(  A\right)
$, and the eventual ranges $G_{\left(  \pi\right)  }\left(  A\right)  $, were
constructed in \textup{(\ref{eqpad.7})} only in the case when the algebraic
prime $\pi$ is in the smaller set of rational primes, i.e., $2,3,5,\dots$; but
the construction in \textup{(\ref{eqpad.7})} goes over \emph{mutatis mutandis}
to the general case, see, e.g., \cite[Sections 4-4--4-5]{Wei98}. In view of
this, it is perhaps surprising that isomorphism of dimension groups in
Corollary \textup{\ref{Corpad.2}} is decided only by the much smaller set
$\operatorname*{Prim}\left(  \det\left(  A\right)  \right)  $.

In the case when $p\in\operatorname*{Prim}\left(  \det\left(  A\right)
\right)  $, then we saw that $G_{\left(  p\right)  }\left(  A\right)  $ is
derived from the space%
\begin{equation}
\sum_{\mu}\left\{  V_{\mu}\left(  \operatorname*{row}\right)  \mid\mu
\in\operatorname*{spec}\left(  A\right)  ,\;p\nmid\mu\right\}  ,
\label{eqRempadNew.3(1)}%
\end{equation}
where $V_{\mu}\left(  \operatorname*{row}\right)  $ is defined analogously as
in \textup{(\ref{eqIntNew.28})} by
\begin{equation}
V_{\mu}\left(  \operatorname*{row}\right)  =\left\{  x\in
V^{\operatorname*{tr}}\mid\exists\,k\in\mathbb{Z}^{+}\ni\joinrel\relbar
x\left(  \mu-A\right)  ^{k}=0\right\}  . \label{eqRempadNew.3(2)}%
\end{equation}
As noted, this sum space is initially computed in $F^{N}$, for a finite-index
field extension $F$. But, since $G_{\left(  p\right)  }\left(  A\right)
\subset\left(  \mathbb{Z}_{\left(  p\right)  }\right)  ^{N}$ as a
$\mathbb{Z}_{\left(  p\right)  }$-module of row vectors, we conclude that the
field $F$ in question must in fact be embedded in $\mathbb{Q}_{\left(
p\right)  }$.

When the extension field $F$ in \textup{(\ref{eqRempadNew.3(1)})} is computed
for a given $p\in\operatorname*{Prim}\left(  \det\left(  A\right)  \right)  $,
then the existence of this embedding of fields $F\subset\mathbb{Q}_{\left(
p\right)  }$ is a nontrivial consequence of Corollary \textup{\ref{Corpad.2}.}
Such an embedding amounts to the conclusion that the equation $f_{A}\left(
x\right)  =0$ has its roots $\mu$ from \textup{(\ref{eqRempadNew.3(1)})} in
$\mathbb{Q}_{\left(  p\right)  }$, where $f_{A}$ is the characteristic
polynomial of $A$. The roots in question must then in fact be in
$\mathbb{Z}_{\left(  p\right)  }$, since $f_{A}$ is monic. Hence, this
solvability of the characteristic equation in $\mathbb{Z}_{\left(  p\right)
}$ is a subtle consequence of Corollary \textup{\ref{Corpad.2},} since we show
that $C^{\ast}$-isomorphism is decided by the $\mathbb{Z}_{\left(  p\right)
}$-modules $G_{\left(  p\right)  }\left(  A\right)  $; and, in particular,
that the latter are nonzero as submodules in $\left(  \mathbb{Z}_{\left(
p\right)  }\right)  ^{N}$. Of course, after knowing existence, there is the
practical issue of having algorithms for finding the solutions.

This issue of field embeddings is addressed algorithmically in \cite[Chapter
1, Section 5]{BoSh66}. Our Example \textup{\ref{ExamplebisJul28}} below
further illustrates the point: The equation $x^{2}=2$ is solvable in
$\mathbb{Z}_{\left(  7\right)  }$, and so we get a natural field embedding of
$\mathbb{Q}\left[  \sqrt{2}\right]  $ into $\mathbb{Q}_{\left(  7\right)  }$,
but not into, for example, $\mathbb{Q}_{\left(  5\right)  }$. \textup{(}The
polynomial $x^{2}-2$ is irreducible in $\mathbb{Z}_{\left(  5\right)  }\left[
x\right]  $.\textup{)}

Similarly, for the complex case, the equation $x^{2}+1=0$ is solvable in
$\mathbb{Z}_{\left(  5\right)  }$ and in $\mathbb{Z}_{\left(  13\right)  }$,
but not in $\mathbb{Z}_{\left(  2\right)  }$ nor in $\mathbb{Z}_{\left(
7\right)  }$. And so we have field embeddings $\mathbb{Q}\left[  i\right]
\hookrightarrow\mathbb{Q}_{\left(  5\right)  }$ and $\mathbb{Q}\left[
i\right]  \hookrightarrow\mathbb{Q}_{\left(  13\right)  }$, but not an
analogous embedding of $\mathbb{Q}\left[  i\right]  $ into $\mathbb{Q}%
_{\left(  2\right)  }$, nor into $\mathbb{Q}_{\left(  7\right)  }$. More
generally, if $p$ is odd, then $x^{2}+1$ is irreducible over $\mathbb{Q}%
_{\left(  p\right)  }$ if and only if $p\equiv3\pmod{4}$, while for
$p\equiv1\pmod{4}$, $x^{2}+1$ has two distinct roots in $\mathbb{Z}_{\left(
p\right)  }$, by Hensel's theorem; see \cite[Section 2-4-7, p.~62]{Wei98}. The
solutions in the respective $\mathbb{Z}_{\left(  p\right)  }$'s may be found
by the standard $p$-adic algorithms, e.g., the Newton scheme \cite[Chapter 1,
Section 6]{BoSh66}.
\end{remark}

\begin{example}
\label{ExapadNew.3}A very simple example illustrating Corollary
\textup{\ref{Corpad.2}} is the pair%
\[
A=%
\begin{pmatrix}
1 & 1\\
2 & 0
\end{pmatrix}
,\qquad B=%
\begin{pmatrix}
0 & 1\\
2 & 1
\end{pmatrix}
.
\]
Then the matrix%
\[
J=%
\begin{pmatrix}
1 & 0\\
1 & 1
\end{pmatrix}
\]
defines an isomorphism of $G\left(  A\right)  $ onto $G\left(  B\right)  $.
Since $2^{n}\underset{n\rightarrow\infty}{\longrightarrow}0$ in $\mathbb{Z}%
_{\left(  2\right)  }$, the respective eigenspaces for the eigenvalue $2$ do
not contribute to the $2$-adic row spaces, and only the $-1$ eigenspaces
contribute. A simple computation shows%
\[
\left(  -2,1\right)  \cdot%
\begin{pmatrix}
1 & 0\\
1 & 1
\end{pmatrix}
=\left(  -1,1\right)  ,
\]%
\[
G_{\left(  2\right)  }\left(  B\right)  =\mathbb{Z}_{\left(  2\right)
}\left(  -2,1\right)  ,\qquad G_{\left(  2\right)  }\left(  A\right)
=\mathbb{Z}_{\left(  2\right)  }\left(  -1,1\right)  \text{\qquad(using Remark
\ref{RempadNew.3}),}%
\]
and
\[
\operatorname*{Prim}\left(  \det\left(  A\right)  \right)
=\operatorname*{Prim}\left(  \det\left(  B\right)  \right)  =\left\{
2\right\}  ,
\]
so \textup{(\ref{eqCorpad.2})} holds.
\end{example}

\section{\label{Dec}Decidability of $C^{\ast}$-equivalence}

In order to digest the steps taken in this central section of the paper, the
reader might find it useful to read this section in conjunction with the road
map in the following section, Section \ref{SX}.

In this section we will prove that the problem of finding an integer matrix
$J=J\left(  1\right)  $, satisfying any of the equivalent conditions
(\ref{eqInt.12})--(\ref{eqInt.13}), (\ref{eqInt.17})--(\ref{eqInt.18}),
(\ref{eqInt.17})$^{\prime}$--(\ref{eqInt.18}), (\ref{eqInt.19}),
(\ref{eqInt.20}), (\ref{eqpad.12})$^{\prime}$ together with positivity, is
decidable.\textbf{ }In these considerations, positivity and singularity will
be dispensed with as in Section \ref{SY}, i.e., we will henceforth assume in
this section that $A$ and $B$ are nonsingular matrices with integer matrix
entries, with the property that $A$, $B$ have positive eigenvalues
$\lambda_{\left(  A\right)  }$, $\lambda_{\left(  B\right)  }$ dominating
strictly all other eigenvalues in absolute value, and such that the
corresponding left eigenvectors $v\left(  A\right)  $ and $v\left(  B\right)
$ are unique up to a scalar multiple, and for a suitable choice of this
scalar, a $g\in G\left(  A\right)  =\bigcup_{n}A^{-n}\mathbb{Z}^{N}$ is
positive if and only if $%
\ip{v\left( A\right) }{g}%
>0$ or $g=0$. We will, as partially explained in Section \ref{Sub}, work in
various algebraic extensions $R$ of $\mathbb{Z}$. The idea is roughly that if
$J$ satisfies (\ref{eqInt.12}):%
\begin{equation}
J\left(  1\right)  G\left(  A\right)  =G\left(  B\right)  , \label{eqDec.1}%
\end{equation}
then $J\left(  1\right)  $ also satisfies%
\begin{equation}
J\left(  1\right)  \left(  R\otimes_{\mathbb{Z}}G\left(  A\right)  \right)
=R\otimes_{\mathbb{Z}}\left(  G\left(  B\right)  \right)  , \label{eqDec.2}%
\end{equation}
and, conversely, if (\ref{eqDec.2}) has no solution $J\left(  1\right)  \in
M_{N}\left(  R\right)  $, then (\ref{eqDec.1}) certainly has no solution, and
this can be used to decide absence of $C^{\ast}$-equivalence.

The operator $J$, as a mapping of the column vectors in $R\otimes_{\mathbb{Z}%
}G\left(  A\right)  $, must preserve Galois conjugation (see the end of
Section \ref{Sub}). We will see in Proposition \ref{ProDec.2} that the
conditions (\ref{eqDec.1})--(\ref{eqDec.2}) amount to having a linear mapping
which preserves a lattice of subspaces defined by a lattice of basis elements
over an extension field, having only specified primes in its determinant, and
satisfying congruences. The lattices of subspaces are sums of generalized
eigenspaces (see (\ref{eqIntNew.28}) for the definition of \emph{generalized
eigenspace\/}). The summands are determined by conditions of divisibility of
eigenvalues by algebraic primes and by the Perron--Frobenius eigenvector. In
addition we can multiply the matrix $J\left(  1\right)  $ by powers of $A$,
$B$ which can automatically make it divisible by any power of an algebraic
prime $\pi$ at the $\pi$-eigenspace. We will show that these conditions are decidable.

By congruences, we mean that a finite set of vectors over a ring $R$ has its
image modulo some ideal $I$ to lie in a specified finite set in $R/I$. In
particular, any Boolean or logical combination of congruences is a set of
congruences. We can test congruences by testing each element of this set of
residue classes.

Over the integers, a matrix which preserves a sequence of rational subspaces
in a direct sum decomposition can be conjugated into a block-triangular form,
by taking bases over the integers corresponding to the sequence of subspaces
\cite{New72}. Every subgroup of a free abelian group is free, and a finitely
generated subgroup is a summand if and only if it has no elements which are
not divisible by a prime $p$ in it but are divisible in the total group
\cite{Kap69}. However, an integer matrix which preserves a sequence of
rational subspaces in a direct sum decomposition cannot always be conjugated
further to be block-diagonal over the integers without introducing fractions.

In an algebraic number ring, some finite, computable power of any ideal (the
order of the class group \cite{Ser79}) will be principal. It basically follows
from ideas in \cite[Section 5-3]{Wei98} that we can get a finite list of
representatives for the class group, and then we need only have a procedure to
test whether an ideal is principal. If it were principal then we can bound the
norm of some field element giving the equivalence and test all possibilities.
A general algorithm is embodied in the free number-theory software PARI; a
theoretical treatment of this problem is in \cite[Section 6.5]{PoZa97},
\cite{Buc86}. A general algorithm for determining the class group is given in
an appendix of \cite{KiRo79}. This means that congruences to a modulus which
is an ideal, or fractions whose denominators lie in an ideal, can be restated
as congruences to a modulus which is an element, or fractions whose
denominators divide a power of some element. Thus we need only to consider
ideals $\left(  m\right)  $ generated by a single element $m\in\Omega$ in the
following lemma, which will be used in the last step of the decision procedure.

\begin{lemma}
\label{LemDec.1}Let $\Omega$ be an algebraic number ring with quotient field
$F$, and let $m_{1}$, $m_{2}$ be relatively prime elements of $\Omega$, i.e.,
$\left(  m_{1}\right)  +\left(  m_{2}\right)  =\Omega$. Let $f\in F$ be
relatively prime to $m_{1}$ also. Let $\operatorname*{CC}\left[  m_{1}%
,m_{2},f\right]  $ be the following set of congruence classes of matrices $M$:%
\begin{multline}
\operatorname*{CC}\left[  m_{1},m_{2},f\right]  =\{M\pmod{m_1}\mid M\in
M_{N}(\Omega\left[  1/m_{2}\right]  )\text{ and }\label{eqLemDec.1}\\
\text{there exists an }x\in\Omega\left[  1/m_{2}\right]  \text{ such that
}1/x\in\Omega\left[  1/m_{2}\right]  \text{ and }\det\left(  M\right)  =fx\}.
\end{multline}
In words, $\operatorname*{CC}\left[  m_{1},m_{2},f\right]  $ is the set of
modulo-$m_{1}$ reductions of matrices $M$ over $F$ whose entries $m_{ij}$ can
be expressed as fractions of elements of $\Omega$ whose denominators divide a
power of $m_{2}$ and such that the determinant of $M$ is a product of $f$ and
units and powers of primes dividing $m_{2}$.

It follows that there is a finite algorithm to list the finite set
$\operatorname*{CC}\left[  m_{1},m_{2},f\right]  $.
\end{lemma}

\begin{remark}
The set $\operatorname*{CC}\left[  m_{1},m_{2},f\right]  $ is finite since
$m_{2}$ is invertible modulo $m_{1}$. It is a subset of $M_{N}\left(
\Omega/\left(  m_{1}\right)  \right)  $, and the quotient ring $\Omega/\left(
m_{1}\right)  $ is finite. That an algorithm determines something means that
the algorithm always gives the correct answer in a finite number of steps.
\end{remark}

\begin{proof}
[Proof of Lemma \textup{\ref{LemDec.1}}]We first show that to any given
modulus such as $m_{1}^{k}$, we can put a matrix modulo $m_{1}^{k}$ into
diagonal form using row and column operations (elementary matrices, having one
non-zero off-main diagonal entry) modulo $m_{1}$. Each such operation lifts to
a similar operation over $\Omega\lbrack1/m_{2}]$, so it preserves the given
set, and moreover, these lifted operations will preserve the norm of the
matrix. The reason that this works is that $\Omega\lbrack1/m_{2}%
]/(m_{1})=\Omega/(m_{1})$ is a principal ideal domain, even though $\Omega$
may not be. This means that every ideal is modulo $m_{1}$ generated by some
element---for this it suffices to factor $m_{1}$ into primes, use the fact
that finite extensions of the $p$-adic integers are principal ideal domains,
and then use the Chinese remainder theorem \cite[Remark 4-1-5 and Theorem
2-2-10]{Wei98} to assemble primes.

We can determine the group $G_{41}$ generated by row operations modulo $m_{1}%
$; it is a subsemigroup of the finite semigroup $M(N,\Omega/(m_{1}))$
generated by a given finite list of generators. The criterion for being in
$\operatorname*{CC}[m_{1},m_{2},f]$ is that a matrix is in $G_{41}D_{41}%
G_{41}$ where%
\begin{multline*}
D_{41}=\{D\pmod{m_1}\mid D\in M(N,\mathbb{Z}_{m_{1}^{N}})\text{ is
diagonal,}\\
\det\left(  D\right)  =xf\pmod{m_1^N}\text{ for an }x\in\Omega\left[
1/m_{2}\right]  \text{ such that }1/x\in\Omega\left[  1/m_{2}\right]  \}.
\end{multline*}
We can determine (list) the finite set of such reductions by determining a
finite list of generators for the group of units of $\Omega\lbrack1/m_{2}]$,
that is, units of $\Omega$ and combinations of prime factors of $m_{2}$, as
well as the class group, the prime factors of $f$, and their images in the
class group. Weiss \cite[Corollary 3-3-3, Proposition 4-4-8]{Wei98} gives
methods for finding an integral basis---we can take some rational basis,
compute its discriminant, and then find integral bases locally at primes which
divide the discriminant (in effect by finding extensions of the ring of
$p$-adic integers). Weiss \cite[Chapters 2, 3, and Section 4-9]{Wei98} further
gives $p$-adic methods for determining what the prime ideals are which lie
over given rational primes $p$. We can factor an element into primes, by
factoring its rational norm into primes \cite[p.~64]{Jac75}, and successively
attempting division by the algebraic primes over $p$. Also from Weiss
\cite[Chapter 5]{Wei98} we have algorithms for bounding the norms of
representatives of the class group, so that to compute the class group, it
remains to tell when two given elements generate the same ideal class, which
is discussed by Pohst and Zassenhaus \cite[Chapter 6]{PoZa97}. These authors
\cite[Chapter 5]{PoZa97} also give an algorithm for finding generators for the
group of units in algebraic number fields.

Necessity of this condition follows by the first paragraph. We now show
sufficiency, i.e., that every such diagonal matrix modulo $m_{1}$ actually
arises from a matrix of $\operatorname*{CC}[m_{1},m_{2},f]$. To do this, we
start with a diagonal matrix the product of whose entries modulo $m_{1}^{N}$
is some $fx$. Then we alter it by multiples of $m_{1}$ so as to modify the
determinant by an arbitrary multiple of $m_{1}^{N}$. To do this, make the
entries $(i,i+1)$ for each $i$ equal to $m_{1}$ and the entry $(N,1)$ equal to
any $xm_{1}\in m_{1}\Omega\lbrack1/m_{2}]$. This adds the single product
$xm_{1}^{N}$ to the determinant.
\end{proof}

We next return to the early steps in the decision procedure and describe the
algebra of endomorphisms which preserves a collection of subspaces (those in
the next definition).

\begin{definition}
\label{DefDecNew.2}Let $A$, $B$ be matrices with rational matrix entries.
Assume that $A$, $B$ are nonsingular, and that $A$, $B$ have positive
eigenvalues $\lambda_{\left(  A\right)  }$, $\lambda_{\left(  B\right)  }$
which are larger than the absolute value of all other eigenvalues, with
corresponding left eigenvectors $v\left(  A\right)  $, $v\left(  B\right)  $
unique up to scalars. Let $K$ denote the field generated by the eigenvalues of
$A$ and $B$. Assume that $A$ and $B$ act on vector spaces $V$, $W$ over
$\mathbb{Q}$. Let $R$ be a subring of $K$. Then $A$ and $B$ act in a natural
manner on $V\otimes_{\mathbb{Z}}R$ and $W\otimes_{\mathbb{Z}}R$ respectively.
Then $\operatorname*{DGI}(A,B,R)$ denotes the additive group of $R$%
-homomorphisms $J\left(  1\right)  $ from $V\otimes R$ to $W\otimes R$ such that

\begin{enumerate}
\item[$(\alpha)$] the direct sums $v\left(  A\right)  ^{\perp}$ and $v\left(
B\right)  ^{\perp}$ of all nonmaximal generalized eigenspaces are mapped into
each other, and, more generally, for each $g$ in the Galois group of $K$ over
$\mathbb{Q}$, the Galois conjugate $\left(  v\left(  A\right)  ^{\perp
}\right)  ^{g}$ is mapped into $\left(  v\left(  B\right)  ^{\perp}\right)
^{g}$,

\item[$(\beta)$] for each algebraic prime $\pi$ of $K$ which divides an
eigenvalue, $J\left(  1\right)  $ preserves the span $E\left(  \pi\right)  $
of the generalized eigenspaces whose eigenvalues are divisible by $\pi$.
\end{enumerate}

\noindent We shall use the abbreviation $\operatorname*{DGI}$, for ``dimension
group isomorphisms'', although ``dimension group pre-isomorphisms'' would be a
more accurate description.

If $\pi=p$ is an ordinary prime, then $E\left(  \pi\right)  \otimes
\mathbb{Z}_{\left(  p\right)  }$ is the orthogonal complement of the
$\mathbb{Z}_{\left(  p\right)  }$-eventual row space $G_{\left(  p\right)
}\left(  A\right)  $ of $A$ defined in {\upshape(\ref{eqpad.11}%
)--(\ref{eqpad.12}),} which is spanned by generalized row eigenvectors for
eigenvalues not divisible by $p$.
\end{definition}

Note that $\operatorname*{DGI}\left(  A,B,R\right)  $ really depends on $A$,
$B$ and not merely on $V$, $W$, because the generalized eigenspaces and
eigenvalues of $A$ and $B$ occur in these conditions. Then our criterion for a
dimension group isomorphism says that there is such a map $J(1)$ defined over
$\mathbb{Z}$ with the following properties for all algebraic primes $\pi$
dividing $\det\left(  A\right)  $ and thus $\det\left(  B\right)  $. (We
identify $J\left(  1\right)  $ with the map it defines on various sub- and quotient-modules.)

\begin{enumerate}
\item \label{DecDGIadditional.1}$J\left(  1\right)  $ is nonzero modulo the
non-Perron--Frobenius generalized eigenspaces (which can be ensured by
congruences relatively prime to $\pi$),

\item \label{DecDGIadditional.2}on the quotient $V/E\left(  \pi\right)  $ the
determinant of $J\left(  1\right)  $ is relatively prime to $\pi$,

\item \label{DecDGIadditional.3}the determinant of $J\left(  1\right)  $ is
divisible only by the primes $\pi$.
\end{enumerate}

\noindent Here (\ref{DecDGIadditional.1}), (\ref{DecDGIadditional.2}) are
congruence conditions and (\ref{DecDGIadditional.3}) is a determinant
condition; these will be transformed a little so that they become the basic
criteria whose satisfiability we must decide. By linear algebra, as outlined
in the next paragraph, we find a nonsingular map $J_{0}$ over the rational
numbers satisfying the first two conditions $(\alpha)$--$(\beta)$, if it
exists, from $V$ to $W$, and then a general hypothetical map must differ from
$J_{0}$ by a map $J_{a}$ in $\operatorname*{DGI}(A,A,\mathbb{Q})\colon
J(1)=J_{0}J_{a}$. Replace $J_{0}$ by some $c_{0}J_{0}$, $c_{0}\in\mathbb{Q}$,
so that $J_{0}^{-1}\in M_{N}(\mathbb{Z})$, where $M_{N}\left(  R\right)  $ is
the algebra of all $N\times N$ matrices over the ring $R$. Then $J_{a}\in
M_{N}(\mathbb{Z})$. Write $J_{0}=J_{c}/N_{c}$, $J_{c}\in M_{N}(\mathbb{Z})$,
$N_{c}\in\mathbb{Z}$. Then (\ref{DecDGIadditional.1}),
(\ref{DecDGIadditional.2}), (\ref{DecDGIadditional.3}) translate into
congruence conditions and norm conditions on $J_{a}$: \renewcommand{\theenumi
}{\roman{enumi}$_{\text{a}}$}

\begin{enumerate}
\item \label{DecDGIsuba.1}$J_{c}J_{a}$ on the chosen maximal eigenvector
$v\left(  A\right)  $ of $A$ is nonzero modulo $p_{a}$ (a fixed prime
relatively prime to $\det(A)$, $\det(J_{c})$, $N_{c}$);

\item \label{DecDGIsuba.2}$J_{c}J_{a}$ on the quotient $V/E(\pi)$ has
determinant a multiple of $N_{c}^{N}$ by a invertible number modulo $\pi$;

\item \label{DecDGIsuba.3}the determinant of $J_{a}$ is $N_{c}^{N}/\det
(J_{c})$ times a number dividing some power of $\det(A) \det(B)$;

\item \label{DecDGIsuba.4}
\begin{equation}
J_{c}J_{a}\equiv0 \pmod{N_c}. \label{eqDecDGIsuba.4.1}%
\end{equation}
\end{enumerate}

The vector space $\operatorname*{DGI}(A,A,\mathbb{Q})$ is in fact also an
algebra, which we next describe. Let $K$ now denote the field generated by the
eigenvalues of $A$. The next proposition is based on general principles of
Galois theory, see, e.g., \cite{Rot98} and \cite{Jac75}, as well as standard
facts about linear resolutions, see \cite{New72}, \cite{Ser77}, \cite{Ser98}.
First, we can find a basis for the set of linear mappings between two vector
spaces $V$, $W$ which map a finite list of subspaces $X_{i}$ into another list
$Y_{i}$. This can be done by writing these inclusion conditions as linear
equations in the entries of a matrix. Then write out the determinant of a
general matrix in this subspace in terms of variables; if this determinant is
not identically zero as a polynomial, then we can find a nonsingular mapping.
The next proposition also extends a more primitive variant which appeared
earlier in \cite[Corollary 9.5]{BJO99}. To understand the statement of the
proposition, recall the following standard terminology: If $K\supset
\mathbb{Q}$ is a number field, the Galois group $\Gamma=\operatorname*{Gal}%
\left(  K/\mathbb{Q}\right)  $ is the group of automorphisms $g$ of $K$ which
fix $\mathbb{Q}$ pointwise, i.e., $x^{g}=x$ for $g\in\Gamma$ and
$x\in\mathbb{Q}$. But we shall also consider $\Gamma$ as a group of
transformations of column vectors $K^{N}$. If $x=\left(  x_{i}\right)
_{i=1}^{N}\in K^{N}$, we set $x^{g}=\left(  \left(  x_{i}\right)  ^{g}\right)
_{i=1}^{N}$.

The submodules $V_{i}$ of the vector space $K\otimes\mathbb{Z}^{N}$ on which
$A$ acts, in the following proposition, are all direct sums of generalized
eigenspaces of $A$ (see (\ref{eqIntNew.27})--(\ref{eqIntNew.28})), and they
are defined as follows: Recall that, for each algebraic prime $\pi$, $E\left(
\pi\right)  $ is the linear span of the generalized eigenspaces $V_{\mu}$
where the eigenvalue $\mu$ has $\pi$ as a factor. Thus the Galois action
permutes the spaces $E\left(  \pi\right)  $ among themselves. Also throw in
$v\left(  A\right)  ^{\perp}$ and its Galois conjugates $\left(  v\left(
A\right)  ^{\perp}\right)  ^{g}$ in addition to the $E\left(  \pi\right)  $'s,
recalling that $v\left(  A\right)  ^{\perp}$ is also a sum of generalized
eigenspaces. Note that this implies that Galois conjugation by $g$ will send
the generalized eigenspace $V_{\mu}$ for any eigenvalue $\mu$ to the
generalized eigenspace $V_{\mu^{g}}$ for $\mu^{g}$ (since it will send, for
example, generating eigenvectors of one to those of the other, if we make the
first coordinate $1$, and will send $K$-linear combinations to possibly
different $K$-linear combinations). So the Galois conjugates of $v\left(
A\right)  ^{\perp}$ will still be sums of generalized eigenspaces. Thus each
finite intersection
\begin{multline*}
E\left(  \pi_{1}\right)  \cap E\left(  \pi_{2}\right)  \cap\dots\cap E\left(
\pi_{n}\right)  \cap\left(  v\left(  A\right)  ^{\perp}\right)  ^{g_{1}}%
\cap\left(  v\left(  A\right)  ^{\perp}\right)  ^{g_{2}}\cap\dots\cap\left(
v\left(  A\right)  ^{\perp}\right)  ^{g_{k}},\qquad\\
g_{1},g_{2},\dots,g_{k}\in\Gamma=\operatorname*{Gal}\left(  K/\mathbb{Q}%
\right)  ,
\end{multline*}
is a direct sum of generalized eigenspaces (if nonzero). The Galois group of
$K$ over $\mathbb{Q}$ must map each such finite intersection into another one.
By Definition \ref{DefDecNew.2}, all these finite intersections are preserved
by $\operatorname*{DGI}(A,A,K)$. Now, choose a linear ordering $I_{i}$,
$i=1,\dots,l$, of these intersections (the ordering is not unique), such
that\renewcommand{\theenumi}{\arabic{enumi}}

\begin{enumerate}
\item \label{PreProDec.2(1)}if $I_{i}\supseteq I_{j}$ then $j\leq i$,

\item \label{PreProDec.2(2)}if $I_{i}$ is a Galois conjugate of $I_{j}$ and
$i<j$, then for all $k$ with $i\leq k\leq j$, $I_{k}$ is also a Galois
conjugate of $I_{j}$.
\end{enumerate}

\noindent Define $V_{j}=\bigoplus_{i=j}^{l}I_{i}$. This gives a decreasing
filtration. Since $I_{i}$ are invariant subspaces of $\operatorname*{DGI}%
(A,A,K)$, all $V_{i}$ are also invariant subspaces. It is not true that
$\operatorname*{DGI}(A,A,K)$ is precisely the algebra fixing all $V_{i}$, but
Proposition \ref{ProDec.2} next gives a partial converse. Since the
construction of the filtration above is rather involved, we have fleshed it
out in a simple example (Example \ref{ExaDec} below) in order to highlight the idea.

\begin{proposition}
\label{ProDec.2}There is a filtration $V_{i}$, $i=0,\dots,l$, of the vector
space on which $A$ acts in which $\operatorname*{DGI}(A,A,K)$ has a
block-triangular structure. The ideal $\mathcal{J}=\{M\in\operatorname*{DGI}%
(A,A,K)\mid MV_{i}\subset V_{i+1},\;i=0,\dots,l-1\}$ is a nilpotent ideal and
$\operatorname*{DGI}(A,A,K)/\mathcal{J}$ has a natural embedding by the block
structure into $\bigoplus_{i}\operatorname*{GL}(V_{i}/V_{i+1})$. This
embedding is an isomorphism. There is a subfiltration $V_{s(i)}$ defined over
$\mathbb{Q}$ such that $V_{s(i)}/V_{s(i+1)}$ is a direct sum of Galois
conjugates of $V_{s(i+1)-1}/V_{s(i+1)}$. These structures can be finitely computed.
\end{proposition}

\begin{proof}
We find the eigenvalues of $A$, diagonalize $A$ over $K$, factoring ideals
into primes, using standard algorithms, e.g., \cite{PoZa97}. Define $I_{i}$
and $V_{i}$ as in the paragraph before the proposition. The effect of Galois
action and the families of intersections of these spaces can be considered by
taking Galois-invariant bases $B_{\pi}$ for $E\left(  \pi\right)  $. We have
ordered the intersections $I_{j}$ with bases $B_{j}$ by inclusion, and have
put Galois conjugates next to each other. Then the subspace generated by all
bases succeeding any given basis is preserved, and we have a block-triangular
structure corresponding to it, and a larger block-triangular structure, whose
blocks are the sets of Galois-conjugate blocks of those from the former
structure. The latter will be defined over $\mathbb{Q}$ as required. Since the
elements of $\mathcal{J}$ strictly increase filtration, any $l$-fold product
of elements of $\mathcal{J}$ is zero, where $l$ is the filtration length,
i.e., the elements of $\mathcal{J}$ are the matrices in the algebra which are
zero on the main-diagonal blocks, and so the quotient maps isomorphically into
the sum of the general linear groups on $V_{j}/V_{j+1}$ with basis
$B_{0j}=B_{j}\setminus\bigcup_{k>j}B_{k}$. Each $I_{j}$, by induction, and
thus each $V_{j}$ is spanned by the union of the $B_{i}$'s contained in it.
But we note that the general linear group on the span of $B_{0j}$ will
preserve all subspaces $E\left(  \pi\right)  $, and their images will span
each of the required summands, so together they will span the sum. Finally,
the larger filtration mentioned above gives the $s\left(  i\right)  $'s$.$
\end{proof}

It follows that all Galois-invariant linear maps on $V_{s(j)}/V_{s(j+1)}$ will
also lift to $\operatorname*{DGI}(A,A,\mathbb{Q})$.

\begin{proposition}
\label{ProDec.3}Suppose a vector space $V$ over $\mathbb{Q}$ is a direct sum
over an extension field $K\supset\mathbb{Q}$ of Galois-conjugate subspaces
$V_{i}$ \textup{(}with corresponding bases\/\textup{),} transitively permuted
by the Galois group of $K$. Then the algebra generated by automorphisms of $V$
over $\mathbb{Q}$ which preserve each space $V_{i}$ is isomorphic to the
general linear group of $V_{1}$ over the minimal field $K_{1}$ required to
define $V_{1}$, which corresponds to the subgroup $N$ of the Galois group that
sends $V_{1}$ to itself.
\end{proposition}

\begin{proof}
If $V_{1}$ can be defined over a subfield of $K$, then the Galois group of
that field must fix $V_{1}$; conversely if the Galois group fixes $V_{1}$, it
will also fix the complementary sum of generalized eigenspaces, hence it will
fix a projection operator to the subspace whose kernel is the complementary
sum of generalized eigenspaces, and from the columns of a matrix for this
operator, the subspace can be defined.

Given an endomorphism of $V_{1}$ over $K$ which arises from a mapping over
$\mathbb{Q}$, the endomorphisms of all other $V_{i}$ are uniquely determined
as its Galois conjugates. This means we have a one-to-one linear mapping from
endomorphisms of $V$ over $\mathbb{Q}$ fixing $V_{1}$ (and these by Galois
conjugacy fix every $V_{i}$), into the general linear group of $V_{1}$ over
$K$. In fact the image lies in the general linear group over $K_{1}$ since
over it, we can define a projection operator to $V_{1}$. This mapping is also
an epimorphism, since, given any $K_{1}$-linear mapping $h$ of $V_{1}$ to
itself, there are Galois conjugates defined on the other $V_{i}$ (the Galois
operator is unique up to the subgroup fixing $V_{1}$, which also fixes $h$).
We can take the sum of $h$ and its Galois conjugates on the other $V_{i}$, and
the sum will be a Galois-invariant mapping of $V$, and therefore defined over
$\mathbb{Q}$.
\end{proof}

\begin{example}
\label{ExaDec}\textup{(}Illustrating the construction from Proposition
\textup{\ref{ProDec.2}} of a covariant filtration.\textup{)} Consider some
matrix $A$ with three eigenvalues $p$, $q$, $pq$ with respective generalized
eigenspaces $E_{1}$, $E_{2}$, $E_{3}$, so that the two sum spaces $E_{1}\oplus
E_{3}$ and $E_{2}\oplus E_{3}$ are preserved under the Galois action, as is
their intersection $E_{3}$. Then the algebra of endomorphisms has a
block-triangular structure with three blocks and the main-diagonal blocks are
isomorphic to the respective endomorphism algebras $\operatorname*{End}%
(E_{1})$, $\operatorname*{End}(E_{2})$, $\operatorname*{End}(E_{3})$. Suppose
now that $p$ and $q$ are Galois conjugates so that the product $pq$ is
Galois-invariant. The larger block structure will then correspond to the two
spaces $E_{1}\oplus E_{2}$ and $E_{3}$. The group of endomorphisms of
$E_{1}\oplus E_{2}$ over the rational numbers will be isomorphic to the
automorphisms of $E_{1}$ over a quadratic extension field corresponding to the
Galois conjugation which interchanges $p$ and $q$.
\end{example}

We now describe how the elements in $\operatorname*{DGI}\left(  A,A,\mathbb{Q}%
\right)  $ may be put into block-triangular form.

As indicated in the paragraph before Proposition \ref{ProDec.2}, the $V_{j}$'s
arise by taking direct sums of intersections of the $E(\pi)$'s and $\left(
v\left(  A\right)  ^{\perp}\right)  ^{g}$'s, ordered in such a way as to
refine the partial order by inclusion of subspaces, and such that Galois
conjugates are adjacent. Now add all $V_{j}$'s which are smaller in the order,
to each given element, to make this a decreasing sequence of subspaces. Then
all $V_{j}$'s are invariant subspaces of $\operatorname*{DGI}(A,A,K)$. Choose
a base of algebraic integer vectors for each generalized eigenspace, so that
the Galois conjugate of a base is chosen as a base for the Galois conjugate
subspace. Choose also a second basis for the sum of all Galois conjugates of
each generalized eigenspace, which exists over the rational integers. Then
each $V_{j}$ is a sum of generalized eigenspaces, so it is spanned by the
union of bases which are in it; likewise each $V_{s(i)}$ is spanned over
$\mathbb{Q}$ by the union of the second basis elements which are in it. Make
the elements of the second basis, in order, the columns of a matrix $J_{f}$.
Then conjugation by $J_{f}$ will put the matrices in $\operatorname*{DGI}%
(A,A,\mathbb{Q})$ into block triangular form, because all columns
corresponding to each subspace $V_{s(i)}$, which has the form, all basis
vectors $v_{i}$, $i\geq n_{0}$, and will span an invariant subspace. In this
basis, Galois conjugation is expressed by permutation of basis elements. Then
substitution of $J_{a}=J_{f}J_{g}\det(J_{f})^{-1}J_{f}^{-1}$ applied to
(i$_{\text{a}}$), (ii$_{\text{a}}$), (iii$_{\text{a}}$), (iv$_{\text{a}}$)
gives the corresponding formulas (i$_{\text{g}}$), (ii$_{\text{g}}$),
(iii$_{\text{g}}$), (iv$_{\text{g}}$).\renewcommand{\theenumi}{\roman
{enumi}$_{\text{g}}$}

\begin{enumerate}
\item \label{DecDGIsubg.1}$J_{f}J_{g}\det(J_{f})J_{f}^{-1}$ on a chosen
maximal eigenvector of $A$ is nonzero modulo $p_{a}$ (which is a fixed prime
relatively prime to $\det(A)$, $\det(J_{c})$, $N_{c}$, $\det(J_{f}))$.

\item \label{DecDGIsubg.2}$J_{c}J_{f}J_{g}\det(J_{f})J_{f}^{-1}=J_{c}J_{a}%
\det(J_{f})^{2}$ on the quotient $V/E(\pi)$ has determinant a multiple of
$\det(J_{f})^{2N}N_{c}^{N}$ by a rational integer which is invertible modulo
$\pi$.

\item \label{DecDGIsubg.3}$\det(J_{g})=\det(J_{f})^{N}\det(J_{a})$ is
$\det(J_{f})^{N}N_{c}^{N} /\det(J_{c})$ times a number dividing some power of
$\det(A)\det(B)$.

\item \label{DecDGIsubg.4}
\begin{align}
J_{f}J_{g}\det(J_{f})J_{f}^{-1}  &  \equiv0\pmod{ \det(J_f)^2 }%
,\label{eqDecDGIsubg.4.1}\\
J_{c}J_{f}J_{g}\det(J_{f})J_{f}^{-1}  &  \equiv0\pmod{N_c det (J_f)^2 }
\label{eqDecDGIsubg.4.2}%
\end{align}
\end{enumerate}

\noindent The first of these says that $J_{a}$ is an integer matrix and the
second is the same as (\ref{eqDecDGIsuba.4.1}) multiplied by $\det(J_{f})^{2}%
$. (Any further multiples by constant matrices could be treated in similar
fashion; we are multiplying matrices by these quantities, so when we take
determinants we multiply by $N$th powers).\medskip

We now prove two general propositions about congruences which will be needed.

Recall some aspects of the theory of finite-dimensional algebras $\mathcal{A}$
with unit over a field. The Jacobson radical is the intersection of all
maximal proper ideals, equivalently the maximal nilpotent ideal, equivalently
in characteristic $0$, the kernel $\left\{  x\in\mathcal{A}\mid
\operatorname*{Trace}(xy)=0,\;\forall\,y\in\mathcal{A}\right\}  $, where
algebra elements are represented as matrices acting on a basis for the algebra
\cite[Ch.~I.14, p.~62]{Jac75}. Modulo the Jacobson radical, the algebra is
semisimple, which means it has no nilpotent ideals, and then that every
element $a$ is regular in the sense there exists $x$ such that $axa=a$. A
semisimple algebra is isomorphic to a direct sum of simple algebras
\cite{vdW91}; this decomposition is unique, and corresponds to the set of
central idempotents of the algebra. Simple algebras will be given as matrix
algebras over algebraic number rings. However, simple finite-dimensional
algebras must always be full matrix algebras over division rings.

We will apply the next proposition to integer matrices in $J_{f}%
^{-1}\operatorname*{DGI}(A,A,\mathbb{Q})J_{f}$ and the congruences
(\ref{DecDGIsubg.1}), (\ref{DecDGIsubg.2}), (\ref{DecDGIsubg.4}), and the
determinant condition (\ref{DecDGIsubg.3}).

Note that we can write any Boolean combination of congruences on a single
matrix variable $x$ to various moduli in the form
\begin{equation}
\exists\,s\in S\ni\joinrel\relbar x\equiv s\pmod{m} \label{eqDecNew.8}%
\end{equation}
for a finite computable set $S$, by \cite{Wei98}. In the application of
Proposition \ref{ProDec.5}, $m$ can be taken as, say, the product of the
$2N$th power of all denominators and determinants for $A$, $B$, $M_{c}$,
$M_{f}$, $p_{a}$.

The terminology in the following proposition, that we can solve a finite
system of congruences, means that there is an algorithm to determine whether
solutions exist, and to find some solution if it exists. The determinant
restrictions are those stated in the proof.

\begin{proposition}
\label{ProDec.5}Let $\mathcal{A}$ be a finite-dimensional algebra of matrices
over a commutative ring $R$ in block-triangular form, and let $J$ be its
Jacobson radical consisting of matrices which have zero main-diagonal blocks.
If we can solve any finite system of additive congruences on $\mathcal{A}/J$
subject to any restrictions on the determinant, then we can solve any finite
system of additive congruences on $\mathcal{A}$ subject to any restrictions on
the determinant. More generally we can restate the congruences on
$\mathcal{A}$ as congruences on $\mathcal{A}/J$ and use the same determinant conditions.
\end{proposition}

\begin{proof}
Note that for our matrix representation the norm conditions on $\mathcal{A}$
will give norm conditions on $\mathcal{A}/J$, since the latter gives the
main-diagonal blocks in a block-triangular representation, and the product of
their determinants is the determinant in $\mathcal{A}$. The condition that the
determinant is a fixed algebraic integer $f$ times products from a finite list
of primes and units will translate into a finite list of similar conditions at
each main diagonal block, based on the prime factorizations of $f$.
Additively, write an element which is to have determinant involving certain
primes, and satisfy congruences, as $x+j$ where $j$ is in the Jacobson
radical. The congruences will say, for some $j\in J$, a Boolean combination of
congruences $x+j\equiv c\pmod{m}$ hold. If we take all possibilities $j_{0}$
for $j\pmod{m}$, this will be a Boolean combination of congruences $x\equiv
c-j_{0}\pmod{m}$.
\end{proof}

Congruences on an element of an algebraic number ring $\Omega$ modulo $m$ will
not be changed if we pass to an extension field (but require the element to
belong in the original ring), and it will suffice to take congruences modulo
the prime power factors of $m$ in the new ring, that is, if $\Omega_{1}$ is
the algebraic number ring of a finite extension of the quotient field of
$\Omega$, and if $x,m\in\Omega$, then $x$ is divisible by $m$ in $\Omega$ if
and only if $x$ is divisible by $m$ in $\Omega_{1}$.

We make one further transformation of our congruences and determinant
conditions. Since it is of the same nature as the previous changes except that
we must use Proposition \ref{ProDec.3} and Proposition \ref{ProDec.5} in a way
which is difficult to predict, we will not state the formulas explicitly but
describe the changes. Using Proposition \ref{ProDec.5}, we pass to congruences
on the indecomposable blocks of the matrix representations. We use Proposition
\ref{ProDec.3} and a further conjugation to pass to congruences over an
algebraic number field on particular generalized eigenspaces. This will result
in congruences (i$_{\text{h}}$), (ii$_{\text{h}}$), (iv$_{\text{h}}$), and a
determinant condition (iii$_{\text{h}}$). The conditions (ii), and so on, will
bound the powers of all primes occurring in the determinant of $J(1)$, $J_{a}%
$, $J_{g}$ at that generalized eigenspace, except for those which divide the eigenvalue.

One way to determine the congruences is to find a basis for the space of
matrices satisfying (\ref{DecDGIsubg.1})--(\ref{DecDGIsubg.4}), compute their
images in the sum of main diagonal blocks, and then give a congruence
specifying the span of these images as a finite-index subgroup of a direct sum
of rings of the form $M_{N_{i}}\left( \Omega_{i}[1/\lambda_{i}]\right) $,
where the $\lambda_{i}$ are eigenvalues and $\Omega_{i}$ the algebraic number
rings of the fields in Proposition \ref{ProDec.3}.

\begin{proposition}
\label{ProDec.6}Given the congruences \textup{(i$_{\text{h}}$), (ii$_{\text{h}%
}$), (iv$_{\text{h}}$),} we may construct\linebreak \ equivalent congruences
of the same type in which the moduli for each generalized eigenspace are
relatively prime to the corresponding eigenvalue $\pi$. Moreover it suffices
to find matrices satisfying these conditions with matrix entries in
$\Omega\lbrack1/\pi]$. To do this, split the congruences into ordered tuples
of congruences on the respective indecomposable blocks, between each block and
a corresponding constant matrix, and replace each modulus with its quotient by
all powers of primes dividing the corresponding eigenvalue.
\end{proposition}

\begin{proof}
We can eliminate the other prime factors of moduli and the denominators by
multiplying by a power of the defining matrix $A$ large enough to cancel off
the denominators. That is, if we have a solution mapping $J_{c}$ at a
particular generalized eigenspace which satisfies congruences for all primes
except those which divide the eigenvalue $\lambda$, then $A^{n}J_{a}$ will
produce a solution at all the other primes, which is congruent to zero modulo
any set power of the primes in $\lambda$, and therefore exists over $\Omega$.
And if any solution does exist, multiplication by a large power of $A$ must
produce one which is congruent to zero modulo high powers of the primes in
$\lambda$; hence it is one that can be found in this way. The last statement
follows by invertibility of the matrix $A$ restricted to any eigenspace, at
all primes not dividing the eigenvalue, so that there will be arbitrarily
large powers of $A$ congruent to the identity.
\end{proof}

We are now ready to state and prove our main theorem. The proof is built up
from the previous results, and it is further spelled out in the next section.

As noted in the Introduction, instead of saying ``stationary stable
AF-algebras'' in the following theorem, we might of course say ``dimension
groups defined by direct limits using constant primitive integer matrices as
in (\ref{eqInt.1})''.

\begin{theorem}
\label{ThmDec.7}There is an algorithm to decide isomorphism of stationary
stable AF-algebras arising from primitive integer matrices.
\end{theorem}

\begin{proof}
This result is a consequence of the preliminary discussion and the
propositions above. First reduce the problem to the case of nonsingular
matrices by the method in Section \ref{SY}. This reduces the problem to one of
finding a matrix $J(1)$ which preserves certain subspaces, has certain primes
in its determinants, and satisfies congruences, going from $A$ to $B$. We find
such a matrix $J_{0}$ over the rational numbers; the proposed solution must
differ from it by multiplying with a matrix $J_{a}\in\operatorname*{DGI}\left(
A,A,\mathbb{Q}\right)  \cap M_{N}\left( \mathbb{Z}\right) $ meeting
corresponding conditions (we multiply by a constant $N_{c}$ to arrange that
$J_{a}$ have integer entries). We find the Jacobson radical of
$\operatorname*{DGI}\left( A,A,\mathbb{Q}\right) $ and the simple components
of the quotient by it, and restate the congruences in terms of those simple
components. They are determined in terms of certain combinations of
generalized eigenspaces, as general linear groups over algebraic number
fields. Again we conjugate, and obtain a new finite set of congruences on a
tuple of matrices over algebraic number rings (no longer necessarily fields,
because they are images of integer matrices) of the same general nature as the
originals, (i$_{\text{h}}$), (ii$_{\text{h}}$), (iv$_{\text{h}}$), and a
determinant condition (iii$_{\text{h}}$). We use Proposition \ref{ProDec.6} to
ensure that the congruences involve moduli relatively prime to the eigenvalues
and can allow these eigenvalues as denominators. By Lemma \ref{LemDec.1}, we
can solve them.
\end{proof}

\renewcommand{\theenumi}{\roman{enumi}}

\section{\label{SX}The explicit algorithm}

In this section we write out the algorithm hinted at in the proof of Theorem
\ref{ThmDec.7} in more detail. Given two square integer primitive matrices
$A$, $B$, the algorithm can be used to decide whether the associated (ordered)
dimension groups $G\left(  A\right)  $, $G\left(  B\right)  $ are isomorphic
or not.

The algorithm has seven steps, numbered \ref{Alg.1}--\ref{Alg.7}%
.\renewcommand{\theenumi}{\Roman{enumi}}\renewcommand{\labelenumi}{\theenumi.}

\subsection{Algorithm\medskip}

\begin{enumerate}
\item \label{Alg.1}Reduce to the nonsingular case as in Lemma \ref{LemSY.1}
and its proof:

\begin{enumerate}
\item Find the kernel of $A^{N}$, $B^{N}$ over $\mathbb{Q}$.

\item Find a free basis for the integer vectors in these subspaces.

\item Find the images of $A$, $B$ on those vectors, giving associated
nonsingular matrices which henceforth replace $A$, $B$. The rank of the two
nonsingular matrices must have the same value $N$, otherwise the algorithm stops.
\end{enumerate}

\item \label{Alg.2}Determine the eigenspace structure of $A$, $B$ (a general
reference for steps \ref{Alg.2} and \ref{Alg.3} is \cite[Chapter 1, Sections
1--7]{Jac75}):

\begin{enumerate}
\item Determine the Perron--Frobenius eigenvalue.

\item Determine all eigenvalues and corresponding generalized eigenspaces.

\item Determine the field $K$ generated by the eigenvalues and the action of
the Galois group on them.

\item Factor the eigenvalues into powers of algebraic primes, noting the norms
of primes and the action of the Galois group on them. A necessary condition
for isomorphism of dimension groups is that the algebraic primes dividing the
determinants of $A$, $B$ must be the same. Stop if not.
\end{enumerate}

\item \label{Alg.3}Determine $\operatorname*{DGI}(A,B,K)$ and
$\operatorname*{DGI}(A,B,\mathbb{Q})$:

\begin{enumerate}
\item From step \ref{Alg.2}, write out, for each algebraic prime $\pi$, the
sum $E\left(  \pi\right)  $ of generalized eigenspaces such that $\pi$ divides
the corresponding eigenvalue. Also write out $v\left(  A\right)  ^{\perp}$,
the sum of the non-Perron--Frobenius generalized eigenspaces, and its Galois
conjugates. Do this for both $A$, $B$. Call the results $E\left(  \pi\right)
\left(  A\right)  $, $E\left(  \pi\right)  \left(  B\right)  $, $\left(
v\left(  A\right)  ^{\perp}\right)  ^{g}$, $\left(  v\left(  B\right)
^{\perp}\right)  ^{g}$.

\item Write out the linear equations on a matrix $J_{01}$ over $K$, ensuring
that $J_{01}$ maps each $E\left(  \pi\right)  \left(  A\right)  $ onto
$E\left(  \pi\right)  \left(  B\right)  $ and each $\left(  v\left(  A\right)
^{\perp}\right)  ^{g}$ onto $\left(  v\left(  B\right)  ^{\perp}\right)  ^{g}$
when going from $A$ to $B$.

\item Find a $\mathbb{Q}$-basis for $K$ which is acted on by the Galois group,
and expand $J_{01}$ as a matrix $J_{02}$ using this basis to define a basis
for $K^{N}$ over the field $\mathbb{Q}$ (that is, the $K$-linear
transformation $J_{01}$ is a $\mathbb{Q}$-linear transformation $J_{02}$).
Find the linear equations for Galois invariance (hence definability) over
$\mathbb{Q}$ of $J_{02}.$

\item By linear algebra, find a basis for the space of all matrices $J_{02}.$
Write the polynomial for a generic matrix in it, and write out its
determinant. By algebra of polynomials over $K$, find whether this polynomial
is identically zero, and if it is not identically zero, find a matrix $J_{0}$
in it. If it is identically zero, there is no dimension group isomorphism, and
the algorithm stops.
\end{enumerate}

\item \label{Alg.4}Write out congruence conditions on a hypothetical matrix
giving an isomorphism $J(1)$ and associated matrices:

\begin{enumerate}
\item $J(1)$ is nonzero modulo the Perron--Frobenius eigenspace; write this in
terms of congruences relatively prime to the determinants of $A$, $B$. Also
write out the congruences that on $V/E(\pi)$ the determinant is relatively
prime to $\pi$ and the condition that its determinant is divisible only by
primes in $\det(A)$. This is \textup{(\ref{DecDGIadditional.1}),
(\ref{DecDGIadditional.2}), (\ref{DecDGIadditional.3})} of the text after
Definition \ref{DefDecNew.2}.

\item Write $J(1)=J_{0}J_{a}$ where $J_{a}\in\operatorname*{DGI}%
(A,A,\mathbb{Q})$ and restate \textup{(\ref{DecDGIadditional.1}),
(\ref{DecDGIadditional.2}), (\ref{DecDGIadditional.3})} as conditions on the
hypothetical $J_{a}$. This gives \textup{(\ref{DecDGIsuba.1}),
(\ref{DecDGIsuba.2}), (\ref{DecDGIsuba.3}), (\ref{DecDGIsuba.4})} of the text.
\end{enumerate}

\item \label{Alg.5}Compute the filtration $V_{I}$ and the associated
filtration of algebras:

\begin{enumerate}
\item This filtration is obtained as follows. Take all intersections of the
$E\left(  \pi\right)  $ and the Galois conjugates of $V_{0}$, which will be
various sums of generalized eigenspaces. Order them so that the numbers of
sets being intersected increases and that Galois conjugates are adjacent, and
add all previous spaces into the next ones so that the sequence of sets
increases. Take bases for these spaces using Galois permuted bases for the
generalized eigenspaces. (See the paragraph before Proposition \ref{ProDec.2}.)

\item Compute the filtration on the known algebra $\operatorname*{DGI}(A,A,K)$
which arises from this, the ideal $J$ of maps which send each $V_{i}%
\rightarrow V_{i+1}$. Compute the isomorphism $\operatorname*{DGI}%
(A,A,K)\rightarrow\sum_{i}\operatorname*{GL}(V_{i}/V_{i+1}).$ (This is
Proposition \ref{ProDec.2}).

\item Compute the isomorphism from the Galois invariant elements of
\linebreak $\operatorname*{DGI}(A,A,K)/J$ into a direct sum of general linear
groups $G_{s}$ over subfields of $K$. (This is Proposition \ref{ProDec.3}).
\end{enumerate}

\item \label{Alg.6}Further restate the congruences, as after Example
\ref{ExaDec}:

\begin{enumerate}
\item Restate the congruence and norm conditions on $J_{a}$ in terms of the
block triangular form on $\operatorname*{DGI}(A,A,\mathbb{Q})$ which arises
from step \ref{Alg.5}. The map is a conjugation by some matrix $J_{f}$ and we
have conditions \textup{(\ref{DecDGIsubg.1}), (\ref{DecDGIsubg.2}),
(\ref{DecDGIsubg.3}), (\ref{DecDGIsubg.4})} on $J_{g}=\det(J_{f})J_{f}%
^{-1}J_{a}J_{f}$.

\item Restate the congruences and norm conditions in terms of equivalent
conditions on the image element of $J_{g}$ in $G_{s}$, which we call
(i$_{\text{h}}$), (ii$_{\text{h}}$), (iii$_{\text{h}}$), (iv$_{\text{h}}$).

\item Restate the congruences so that we have only congruences on each
generalized eigenspace to moduli relatively prime to the corresponding
eigenvalue (Proposition \ref{ProDec.6}). This involves multiplying by some
power of $A$ and considering the resulting congruences.
\end{enumerate}

\item \label{Alg.7}Solve the congruences, that is determine whether any
solution exists, and if so, find a solution (Lemma \ref{LemDec.1}). We do this
by finding all matrices $N\times N$ modulo moduli $m_{1}$ over an algebraic
number ring $\Omega$ whose entries have denominators dividing $m_{2}$ and
whose determinant is a fixed number $f$ relatively prime to $m_{1}$ times
units in $\Omega\lbrack1/m_{2}]$. This is a finite set and we list it, with representatives:

\begin{enumerate}
\item We determine the basic theory of $\Omega$: integral basis,
multiplication, class group, units, primes dividing $m_{1}$, $m_{2}$, $f$.

\item Determine the group $G_{41}$ of products of elementary matrices modulo
$m_{1}$; we can compute it as a subset of a finite semigroup (all matrices
over $\Omega/(m_{1})$) having given generators.

\item Determine the possible determinants of matrices in this set each of
which is the reduction of a product of units in $\Omega\lbrack1/m_{2}]$ (a
finitely generated group) times products of primes dividing $f$ (a finite set)
which is the identity in the class group.

\item Determine the possible diagonal forms $F_{41}$ modulo $m_{1}$: list
diagonal elements modulo $m_{1}^{N}$ whose product is a possible determinant.
(These conditions are unchanged if we pass to the smaller modulus used in the proof).

\item The required list is $G_{41}F_{41}G_{41}.$
\end{enumerate}
\end{enumerate}

\section{\label{Cas}The case $\lambda=\left|  \det\left(  A\right)  \right|  $}

Theorem \ref{ThmDec.7} and Section \ref{SX} give a finite algorithm (but in
general a long one) to decide whether two square, nonsingular, integer
primitive matrices $A$, $B$ are $C^{\ast}$-equivalent or not. In special
cases, like those considered in \cite{BJO99}, this algorithm can be
substantially simplified. One nice feature of the algorithm is that it uses
only ``elementary'' algebraic results, and avoids using the deep results on
decidability from \cite{GrSe80a,GrSe80b}. Nevertheless, the implementation of
the algorithm for general pairs $A$, $B$ may of course be complicated. Let us
pick up and generalize one special case from \cite{BJO99}. In Theorem 17.18
and Corollary 17.21 there, it was proved that if $A$, $B$ had a special form,
and $\lambda_{\left(  A\right)  }=\left|  \det\left(  A\right)  \right|  $ and
$\lambda_{\left(  B\right)  }=\left|  \det\left(  B\right)  \right|  $, then
the ideal generated by $%
\ip{v\left( A\right) }{w\left( A\right) }%
$ in $\mathbb{Z}\left[  1/\det\left(  A\right)  \right]  $ is a complete
invariant, if the left and right Perron--Frobenius eigenvectors are taken to
have integer components, and $\gcd\left(  v\left(  A\right)  \right)  =1$,
$\gcd\left(  w\left(  A\right)  \right)  =1$, where $\gcd$ denotes the
greatest common divisor of the components. We now prove that this is also true
for more general matrices $A$, $B$.

In stating this more general result, there is a technical complication. In
picking extension fields $F$ and an associated ring $R$ of algebraic integers,
it is not automatically true that the ideals in $R$ are principal. But by a
result in \cite{Wei98} or \cite{Ser79}, there is always a finite extension $E$
of $F$ in which the associated ideals are automatically principal. We refer to
this in the statement of the proposition. To further simplify the terminology
in the statement of the proposition we denote the above-mentioned respective
Perron--Frobenius column vectors $w$, $w^{\prime}$, i.e., $Aw=\lambda w$ and
$Bw^{\prime}=\lambda^{\prime}w^{\prime}$, and similarly $v$, $v^{\prime}$ for
the two respective Perron--Frobenius row vectors.\renewcommand{\theenumi
}{\roman{enumi}}\renewcommand{\labelenumi}{(\theenumi)}

\begin{proposition}
\label{Prop1}Choose a finite extension $E$ of the algebraic number field $F$
of the eigenvalues of primitive nonsingular integer matrices $A$, $B$ in which
all ideals of $F$ become principal and consider primes in it.

\begin{enumerate}
\item \label{Prop1(1)}An isomorphism $J$ on ordered dimension groups from the
dimension group of $A$ to that of $B$ sends the row Perron--Frobenius
eigenvector $v^{\prime}$ \textup{(}normalized so all entries are algebraic
integers with $\gcd1$\textup{)} of $B$ to a multiple $c$ times the row
Perron--Frobenius eigenvector $v$ of $A$.

\item \label{Prop1(2)}The two Perron--Frobenius eigenvalues generate the same
algebraic number field and involve the same primes of that field.

\item \label{Prop1(3)}Assume that for each non-Perron--Frobenius eigenvalue
$\mu$ of $A$ or $B$, the Perron--Frobenius eigenvalue $\lambda$ is divisible
by some algebraic prime not dividing $\mu$. Then the Perron--Frobenius column
eigenvector $w$ is mapped to a multiple $\xi$ times the other
Perron--Frobenius eigenvector $w^{\prime}$.

\item \label{Prop1(4)}If $\lambda$ satisfies the hypothesis in
\textup{(\ref{Prop1(3)}),} then the latter coefficient $\xi$ factorizes into
the primes dividing $\lambda$.

\item \label{Prop1(5)}If $\lambda$ satisfies the hypothesis in
\textup{(\ref{Prop1(3)}),} then the former coefficient $c$ factorizes into
primes dividing $\lambda$. The inner products of left and right
Perron--Frobenius eigenvectors are equal modulo normalization: $v^{\prime
}Jw^{\prime}=cvw^{\prime}=\left(  c/\xi\right)  vJw$. Therefore the inner
product of the normalized Perron--Frobenius eigenvectors is an invariant up to
units in the algebraic number ring generated by $1/\lambda$.
\end{enumerate}
\end{proposition}

\begin{proof}
The first assertion is by \cite[Theorem 6]{BJKR98}, and the second is by
\cite[Theorem 10]{BJKR98}.

The third assertion follows because the space of vectors in the dimension
group $G\otimes E$ such that some fixed multiple is arbitrarily divisible by a
given algebraic prime is sent to the corresponding subspace of the other
dimension group, and this set is the sum of the generalized eigenspaces for
all eigenvalues divisible by the prime. If these spaces are intersected over
all primes dividing the Perron--Frobenius eigenvalue, we get, by our
hypothesis, only the Perron--Frobenius eigenspace.

The fourth assertion follows because the eigenspace of $w$ will consist
precisely of those vectors in the dimension group which are divisible by
arbitrary powers of primes occurring only in $\lambda$, so it must be
preserved by any isomorphism of dimension groups. In addition, vectors in this
$1$-dimensional space which are not divisible by primes other than those in
$\lambda$ will be unique up to multiplication by units and primes dividing
$\lambda$, so they will be preserved by any isomorphism, up to such multiplication.

Next note that the fourth statement and the first part of the fifth statement
are equivalent whenever we have an isomorphism from the dimension group to
itself induced by an integer matrix. The reason is that since row and column
Perron--Frobenius eigenvectors are preserved, this integer matrix in a basis
corresponding to generalized eigenvectors becomes block diagonal, and the
block for the Perron--Frobenius eigenvectors must be the same element for the
row eigenvectors as for the column eigenvectors, and by (\ref{Prop1(4)}) it
involves only primes dividing $\lambda$.

Now consider integer matrices $R$ and $S$ inducing mappings each way between
two different column dimension groups, with Perron--Frobenius eigenvectors $v$
and $v^{\prime}$ normalized over the algebraic number ring. We have
$v\rightarrow c_{1} v^{\prime}$ and $v^{\prime}\rightarrow c_{2} v$, where
$c_{1}$ and $c_{2}$ are algebraic integers since $R$ and $S$ are integer
matrices. But $c _{1}c_{2}=c$ arises from a map of the dimension group to
itself, so it divides a power of $\lambda$, hence so do $c_{1}$ and $c_{2}$.
\end{proof}

The following is a partial converse to Proposition \ref{Prop1}.

\begin{corollary}
\label{ProFred.2.7.99}Suppose $A$, $B$ are nonsingular primitive integer
matrices such that their Perron--Frobenius eigenvalues are integers and that
the inner products as above are equal, i.e., after normalization, that
$\ip{v\left(      A\right)               }{w\left(               A\right)
}=\ip{v\left(      B\right)       }{w\left(            B\right)
}$, that the primes dividing the Perron--Frobenius eigenvalues are equal, and
that the dimensions of the matrices are at least $\,3$. Suppose that the
Perron--Frobenius eigenvalues are the determinants of $A$, $B$ up to sign.
Then there exists an isomorphism between the ordered dimension groups of $A$
and $B$.
\end{corollary}

\begin{proof}
By Lemma 17.19 of \cite{BJO99}, there is a unimodular matrix $J$ sending the
Perron-Frobenius row eigenvector of $A$ to the Perron-Frobenius row
eigenvector of $B$ and the Perron-Frobenius column eigenvector of $A$ to the
Perron-Frobenius column eigenvector of $B$ (and we can choose signs for
positivity). By \cite[Theorem 6]{BJKR98} this gives a positive mapping on
dimension groups. Since the row eigenvectors are perpendicular to the sum
$V\left(  A\right)  =v\left(  A\right)  ^{\perp}$ of all non-Perron-Frobenius
generalized eigenspaces, $JV\left(  A\right)  =V\left(  B\right)  $, and also
$v\left(  B\right)  J\subseteq\mathbb{Q}v\left(  A\right)  $ as noted in
Section \ref{Int}. Write any vector $v$ as a direct sum according to
(\ref{eqIntNew.21}), $v=x+y$. This splitting can introduce certain fixed
primes $p$ in the denominator.

Note that the matrix $A$ is unimodular and integer restricted to the integer
vectors in $V\left(  A\right)  $ (and similarly for the matrix $B$), because
each determinant is the product of its determinant on this space and its
determinant on the Perron--Frobenius eigenspace, and because it is an integer
matrix preserving this subspace. Multiplication by $A$ is multiplication by
$\lambda$ on $x$ (see Figure \ref{FigCas.1}), and the same is true for $B$.%

\begin{figure}[h]
\setlength{\unitlength}{\baselineskip}\settowidth{\submatrixwidth}{$A\sim
\left( \vphantom{
\begin{tabular}{ccccc}
$*$ & $\cdots$ & $*$ & \vline& $0$  \\
$\vdots$
&
\begin{tabular}{c}
$A|_{V\left( A\right) }$ \\
with \\
$\left| \det\left( A|_{V\left( A\right) }\right) \right| =1$
\end{tabular}
& $\vdots$ & \vline& $\vdots$  \\
$*$ & $\cdots$ & $*$ & \vline& $0$  \\ \hline$0$ & $\cdots$ & $0$ & \vline
& \begin{tabular}{c}
$\lambda_{(A)}$ \\
$=\left| \det A \right|$
\end{tabular}
\end{tabular}
}\right. $}\makebox[\submatrixwidth]{}\begin{tabular}{cccccccl}
& & $V\left( A\right) $ &  &  & $\mathbb{C}w\left( A\right
) $ & &  \\[\medskipamount]
& $*$ & $\cdots$ & $*$ & \vline& $0$ & &  \\
\begin{picture}(0,0)
\put(0,-0.75){\makebox(0,0)[r]{$A\sim\left( \vphantom{
\begin{tabular}{ccccc}
$*$ & $\cdots$ & $*$ & \vline& $0$  \\
$\vdots$
&
\begin{tabular}{c}
$A|_{V\left( A\right) }$ \\
with \\
$\left| \det\left( A|_{V\left( A\right) }\right) \right| =1$
\end{tabular}
& $\vdots$ & \vline& $\vdots$  \\
$*$ & $\cdots$ & $*$ & \vline& $0$  \\ \hline$0$ & $\cdots$ & $0$ & \vline
& \begin{tabular}{c}
$\lambda_{(A)}$ \\
$=\left| \det A \right|$
\end{tabular}
\end{tabular}
}\right. $}}
\end{picture}
& $\vdots$
&
\begin{tabular}{c}
$A|_{V\left( A\right) }$ \\
with \\
$\left| \det\left( A|_{V\left( A\right) }\right) \right| =1$
\end{tabular}
& $\vdots$ & \vline& $\vdots$
&
\begin{picture}(0,0)
\put(0,-0.75){\makebox(0,0)[l]{$\left. \vphantom{
\begin{tabular}{ccccc}
$*$ & $\cdots$ & $*$ & \vline& $0$  \\
$\vdots$
&
\begin{tabular}{c}
$A|_{V\left( A\right) }$ \\
with \\
$\left| \det\left( A|_{V\left( A\right) }\right) \right| =1$
\end{tabular}
& $\vdots$ & \vline& $\vdots$  \\
$*$ & $\cdots$ & $*$ & \vline& $0$  \\ \hline$0$ & $\cdots$ & $0$ & \vline
& \begin{tabular}{c}
$\lambda_{(A)}$ \\
$=\left| \det A \right|$
\end{tabular}
\end{tabular}
}\right)$}}
\end{picture}
& $\;V\left( A\right) $ \\
& $*$ & $\cdots$ & $*$ & \vline& $0$ & & \\ \cline{2-6}
& $0$ & $\cdots$ & $0$ & \vline& \begin{tabular}{c}
$\lambda_{(A)}$ \\
$=\left| \det A \right|$
\end{tabular}
& & $\;\mathbb{C}w\left( A\right) $
\end{tabular}
\caption{The case $\left|  \det A \right|  =\lambda_{(A)}$: Decomposition
relative to (\ref{eqIntNew.21}) and unimodular restriction.}\label
{FigCas.1}\end{figure}%

For $v$ to be in the dimension group means for all sufficiently large $n$,
$A^{n}v$ has integer entries. Any prime $p$ which does not divide $\lambda$
will not occur in the denominator of the expression $B^{m}J A^{-n}(x+y)$.

Consider those primes $p$ which divide $\lambda$. We claim that they cannot
occur in denominators of $y$. Restricted to vectors $y$, the matrix $A$ is
unimodular, so modulo any powers of those primes it lies in a finite group,
$\operatorname*{GL}(N,\mathbb{Z}_{p^{s}})$. Thus we can choose arbitrarily
large $n$ so that $A^{n}$ is congruent to the identity. But then in
$A^{n}(x+y)$, the denominators in $x$ have vanished, being multiplied by
$\lambda^{n}$ and those in $y$ remain. So $x+y$ is not in the dimension group,
a contradiction.

Therefore in
\begin{align}
B^{m}JA^{-n}(x+y)  &  =B^{m}JA^{-n}x+B^{m}JA^{-n}%
y\label{eqProFred.2.7.99.proof}\\
&  =J\lambda^{\prime\,m}\lambda^{-n}x+B^{m}JA^{-n}y\nonumber
\end{align}
both terms are integer for sufficiently large $m$ (with a symmetrical argument
the other way) which verifies the conditions in Section \ref{Int} for
isomorphism of ordered dimension groups.
\end{proof}

\section{\label{CanineRat}The case of no infinitesimal elements and the case
of rational eigenvalues}

In this section we will consider the $C^{\ast}$-equivalence problem in two
extreme cases. To describe these two cases, let us recall some facts from
\cite{Eff81}, \cite{BJO99}. We define a functional $\tau_{A}$ on $G\left(
A\right)  $ by the formula%
\begin{equation}
\tau_{A}\left(  g\right)  =%
\ip{v\left( A\right) }{g}%
,\qquad g\in G\left(  A\right)  , \label{eqSomOct.6}%
\end{equation}
where $v\left(  A\right)  $ is a left Perron--Frobenius eigenvector for $A$.
This functional $\tau_{A}$ is called ``the'' trace since it defines a trace on
the corresponding $C^{\ast}$-algebra. It follows from the eigenvalue equation
that $v\left(  A\right)  $ can be taken to have components in the field
$\mathbb{Q}\left[  \lambda\right]  =\mathbb{Q}\left[  1/\lambda\right]  $,
where $\lambda$ is the Perron--Frobenius eigenvalue. But multiplying $v\left(
A\right)  $ by a positive integer, we may assume that the components of
$v\left(  A\right)  $ are contained in the ring $\mathbb{Z}\left[
1/\lambda\right]  $. It then follows from (\ref{eqInt.8}), and $v\left(
A\right)  A^{-n}=\lambda^{-n}v\left(  A\right)  $, that%
\begin{equation}
\tau_{A}\left(  G\left(  A\right)  \right)  \subset\mathbb{Z}\left[
1/\lambda\right]  . \label{eqCanineRat.X}%
\end{equation}
Furthermore, $\tau_{A}\left(  G\left(  A\right)  \right)  $ is invariant under
multiplication by elements of $\mathbb{Z}$ and by $1/\lambda$, so it is a
$\mathbb{Z}\left[  1/\lambda\right]  $-module. In particular, $\tau_{A}\left(
G\left(  A\right)  \right)  $ is an ideal in the ring $\mathbb{Z}\left[
1/\lambda\right]  $. We need only verify that $\frac{1}{\lambda}\tau
_{A}\left(  g\right)  $ is in $\operatorname*{ran}\left(  \tau_{A}\right)  $
for all $g\in G\left(  A\right)  $, where $\operatorname*{ran}\left(  \tau
_{A}\right)  $ denotes the range of the trace functional $\tau_{A}$, i.e., the
subgroup $\tau_{A}\left(  G\left(  A\right)  \right)  $ from
(\ref{eqCanineRat.X}). Pick $g\in G\left(  A\right)  $, and set $g=A^{-n}m$,
$n\in\mathbb{Z}_{+}$, $m\in\mathbb{Z}^{N}$. Then $\frac{1}{\lambda}\tau
_{A}\left(  g\right)  =%
\ip{v\left( A\right) A^{-1}}{g}%
=%
\ip{v\left( A\right) }{A^{-1}g}%
=\tau_{A}\left(  A^{-\left(  n+1\right)  }m\right)  \in\operatorname*{ran}%
\left(  \tau_{A}\right)  $ as claimed. This is a very special feature of the
constant-incidence-matrix situation which is not shared by the range of a
trace on a general dimension group of general AF-algebras. This range is not
even closed under multiplication in the general case when the incidence matrix
is not assumed constant. We have the natural short exact sequence of groups%
\begin{equation}
0\longrightarrow\ker\left(  \tau_{A}\right)  \hooklongrightarrow G\left(
A\right)  \overset{\tau_{A}}{\longrightarrow}\tau_{A}\left(  G\left(
A\right)  \right)  \longrightarrow0 \label{eqSomOct.4}%
\end{equation}
and the order isomorphism
\begin{equation}
G\left(  A\right)  /\ker\left(  \tau_{A}\right)  \underset{\approx}%
{\overset{\tau_{A}}{\longrightarrow}}\operatorname*{ran}\left(  \tau
_{A}\right)  \subset\mathbb{Z}\left[  1/\lambda_{\left(  A\right)  }\right]  ,
\label{eqSomOct.5}%
\end{equation}
where $\operatorname*{ran}\left(  \tau_{A}\right)  $ inherits the natural
order from $\mathbb{Z}\left[  1/\lambda\right]  $. Note that for the
particular matrices we considered in \cite{BJO99}, we had%
\begin{equation}
\operatorname*{ran}\left(  \tau_{A}\right)  =\mathbb{Z}\left[  1/\lambda
\right]  \label{eqCanineRat.Y}%
\end{equation}
(see \cite[(5.21)--(5.22)]{BJO99}), but be warned that this is not a general
feature. This will be discussed further in Remarks \ref{RemSomOct.3} and
\ref{RemSomOct.pound}. Chapter 5 in \cite{BMT87} also has a nice treatment of
$\operatorname*{ran}\left(  \tau_{A}\right)  $ in the general case. Let us
already at this point state and prove the remarkable fact that \emph{any}
subset $I$ of $\mathbb{Q}\left[  \lambda\right]  $ which is an ideal over
$\mathbb{Z}\left[  1/\lambda\right]  $ occurs as the image of the trace for a
suitable primitive nonsingular matrix $A$ (this is a version of
\cite[Corollary 5.15]{BMT87} which is a consequence of results of Handelman,
see \cite{Han81} and \cite{Han87}):

\begin{proposition}
\label{ProCanineRat.1}Let $\lambda$ be a real algebraic integer larger than
the absolute value of any of its conjugates, and let $I\subset\mathbb{Q}%
\left[  \lambda\right]  $ be an ideal over $\mathbb{Z}\left[  1/\lambda
\right]  $. Then $I$ can occur as the image of the trace for some matrix whose
Perron--Frobenius eigenvalue is a power of $\lambda$ \textup{(}the size of the
matrix will be the degree $\mathbb{Q}\left[  \lambda\right]  /\mathbb{Q}$\textup{).}
\end{proposition}

\begin{proof}
Let $I_{1}=I\cap\mathbb{Z}\left[  \lambda\right]  $; it will be a
$\mathbb{Z}\left[  \lambda\right]  $-ideal which spans $I$ over $\mathbb{Z}%
\left[  1/\lambda\right]  $.

Now define an integer matrix $M$ which expresses the action of $\lambda$ on
$I_{2}$, that is, form an additive basis $w_{i}$ for $I_{1}$, let $\lambda
w_{i}=\sum_{j} m_{ij} w_{j}$, $m_{ij}\in\mathbb{Z}$. This matrix will have an
eigenvalue $\lambda$, and we claim that at the corresponding eigenspace, the
image of the trace is isomorphic to $I$. This is because the action of $M$ on
$\mathbb{Z}^{N}$ has been forced to be that of $\lambda$ on $I_{1}$, and
because the trace reflects this module structure, by means of the short
(nearly exact) sequence.

Finally we claim that we can conjugate $M$ over $\operatorname*{GL}%
(N,\mathbb{Z})$ to a matrix whose powers are eventually positive; then those
powers will be nonnegative matrices whose image of trace is the same. To get
eventual positivity, given that $\lambda$ is the largest eigenvalue (the
largest of its Galois conjugates), it is necessary and sufficient that its row
and column eigenvectors for this eigenvalue be positive, by a limit argument
somewhat like that in Proposition \ref{ProInt.1}. Let $v,w$ be row and column
eigenvectors at the eigenvalue $\lambda$, with signs chosen so that their
inner product is positive. Multiply each by a large integer, and then take
relatively prime integers approximating its components. Such a pair of vectors
can be mapped over $\operatorname*{GL}(N,\mathbb{Z})$ to any vectors whose
entries are relatively prime integers having the same inner product, by
\cite[Lemma 17.19]{BJO99}, in particular, to ones which are positive, if
$N>2$. If $N=2$ we use the same result and get a congruence condition, but
that is compatible with positivity.
\end{proof}

\begin{remark}
\label{RemCanineRat.2}The quotient of the ring $\mathbb{Z}\left[
1/\lambda\right]  $ by any of these ideals will be finite. The ideal can be
lifted to an ideal inside the rank-$N$ additive group $\mathbb{Z}\left[
\lambda\right]  $, and the quotient of two rank-$N$ free abelian groups is
finite---its order is given by the determinant of the map expressing the inclusion.
\end{remark}

Let us return to the two special cases of $C^{\ast}$-equivalence we shall
discuss in this section. These are the following.

\begin{enumerate}
\item \label{CanineRat(1)}The kernel $\ker\left(  \tau_{A}\right)  $ is $0$,
i.e., $G\left(  A\right)  $ has no infinitesimal elements, i.e., the
characteristic polynomial of $A$ is irreducible over $\mathbb{Z}$ (equivalent:
over $\mathbb{Q}$).

\item \label{CanineRat(2)}All the eigenvalues of $A$ are rational (thus
integer), each of them is relatively prime to the rest, and none is equal to
$\pm1$.
\end{enumerate}

\noindent In Sections \ref{Som} and \ref{Str} we will apply this to many
examples. See, for example, Example \ref{Exa2} for an application in the
situation (\ref{CanineRat(2)}) above.

\begin{theorem}
\label{Theorem1Oct4}Two primitive $N\times N$ matrices $A$, $B$ over
$\mathbb{Z}_{+}$ with irreducible characteristic polynomials are $C^{\ast}%
$-equivalent if and only if the following three conditions all hold:

\begin{enumerate}
\item \label{Theorem1Oct4(1)}the roots of their characteristic polynomials
generate the same field,

\item \label{Theorem1Oct4(2)}their Perron--Frobenius eigenvalues are divisible
by the same algebraic primes, and

\item \label{Theorem1Oct4(3)}their dimension groups, as modules over
$\mathbb{Z}[1/\lambda]$ \textup{(}or a full-rank subring\/\textup{),} are
isomorphic. These modules are isomorphic to the fractional ideals given by the
image of the trace $\tau$.

\setcounter{enumlink}{\value{enumi}}
\end{enumerate}

\noindent Moreover, these three conditions are equivalent to the one condition:

\begin{enumerate}
\setcounter{enumi}{\value{enumlink}}

\item \label{Theorem1Oct4(4)}the two ordered additive subgroups in
$\mathbb{Z}\left[  1/\lambda\right]  $ defined by the ranges of the respective
traces are isomorphic.
\end{enumerate}

\noindent If in addition the characteristic polynomials of $A$, $B$ are equal,
then $C^{\ast}$-equivalence \textup{(}isomorphism of ordered dimension
groups\/\textup{)} is the same as shift equivalence.

Note that taking powers of the matrix will preserve the $\mathbb{Z}\left[
1/\lambda\right]  $-module mentioned in \textup{(\ref{Theorem1Oct4(3)}),}
i.e., the ideal in $\mathbb{Z}\left[  1/\lambda\right]  $, and not replace it
by its powers.
\end{theorem}

\begin{remark}
\label{RemCanineRat.4}To say that the dimension groups $G\left(  A\right)  $
and $G\left(  B\right)  $ as modules over $\mathbb{Z}\left[  1/\lambda\right]
$ are isomorphic means that there is an isomorphism $\varphi\colon G\left(
A\right)  \rightarrow G\left(  B\right)  $ of abelian groups such that%
\begin{equation}
\varphi\left(  \omega g\right)  =\omega\varphi\left(  g\right)
\label{eqRemCanineRat.4}%
\end{equation}
for all $g\in G\left(  A\right)  $, $\omega\in\mathbb{Z}\left[  1/\lambda
\right]  $. This is not the same as saying that $G\left(  A\right)  $ is
isomorphic to $G\left(  B\right)  $ as ideals in $\mathbb{Z}\left[
1/\lambda\right]  $. The latter concept means that there is an automorphism
$\varphi$ of the ring $\mathbb{Z}\left[  1/\lambda\right]  $ such that
$\varphi\left(  G\left(  A\right)  \right)  =G\left(  B\right)  $. When we
talk about equivalence of ideals it is the \emph{first} concept we are
thinking about, i.e., there is an element of the quotient field $\mathbb{Q}%
\left[  1/\lambda\right]  =\mathbb{Q}\left[  \lambda\right]  $ mapping the one
ideal into the other by multiplication.
\end{remark}

\begin{proof}
[Proof of Theorem \textup{\ref{Theorem1Oct4}}]The first three statements are a
reformulation of \cite[Proposition 10]{BJKR98}, except for the relationship
with the trace, which we next show. The definition (\ref{eqSomOct.6}),
properties (\ref{eqInt.8}), (\ref{eqInt.10}), and
\[
t_{A}\circ A^{-1}=\lambda^{-1}\tau_{A}%
\]
imply that the image of the trace is a module over $\mathbb{Z}[1/\lambda]$ and
a subset of $\mathbb{Z}[1/\lambda]$. Using the standard basis for
$\mathbb{Z}^{N}$, it is generated by $\langle v\mid e_{i}=v_{i}\rangle$ as a
module over $\mathbb{Z}[1/\lambda]$ since $A^{-n}e_{i}$, $n\in\mathbb{Z}_{+}$
generate the dimension group. The trace mapping is an epimorphism if we pass
to rational coefficients (that is, tensor dimension groups with $\mathbb{Q}$),
just because its image is nonzero (consider $v$ as a Perron--Frobenius column
eigenvector) and closed under field operations in $\mathbb{Q}[1/\lambda]$. Its
kernel is zero since the dimension group with rational coefficients is also a
$1$-dimensional vector space over $\mathbb{Q}[1/\lambda]$ (for instance, by
\cite[Chapter 5]{BMT87}). Thus the trace mapping is an isomorphism to its
image as asserted in the second part of (\ref{Theorem1Oct4(3)}).

A theorem of Handelman (stated as Theorem 5.2 in \cite{BMT87}) in fact allows
us to replace the condition in Theorem \ref{Theorem1Oct4} above that $A$ and
$B$ be nonnegative with the condition that instead they are \emph{integral
eventual positive} (IEP), i.e., that they are in $M_{N}\left(  \mathbb{Z}%
\right)  $ and have respective powers with strictly positive entries.

Next we show that dimension group isomorphism in our sense implies shift
equivalence if the irreducible characteristic polynomials of $A$, $B$ are
equal. The only difference with the isomorphisms used in \cite[Theorem
2.8]{BMT87} is that there the actions of $A$, $B$ are the same, that is, the
matrices themselves represent the field element acting on this module. But the
field element $A$, $B$ represent are roots of the same irreducible
characteristic polynomials, and are the unique Perron--Frobenius roots of
these polynomials so they must be the same field element.

Equivalence to (\ref{Theorem1Oct4(4)}): (\ref{Theorem1Oct4(4)}) is, properly
understood, a rephrasing of (\ref{Theorem1Oct4(3)}), given the isomorphism in
the first paragraph of the proof. We will clarify the kind of module structure
which is involved. Assuming (\ref{Theorem1Oct4(4)}), the images of the traces
generate the fields $\mathbb{Q}[\lambda_{(A)}]=\mathbb{Q}[\lambda_{(B)}]$, but
the rings $\mathbb{Z}[\lambda_{(A)}]$ and $\mathbb{Z}[\lambda_{(B)}]$ may be
different, in which case we work with the full-rank subring $\mathbb{Z}%
[\lambda_{(A)}]\cap\mathbb{Z}[\lambda_{(B)}]$. Condition
(\ref{Theorem1Oct4(1)}) is immediate, and condition (\ref{Theorem1Oct4(2)})
follows since algebraic primes dividing $\lambda$ are those primes which can
divide elements of the dimension groups to arbitrary powers.

Conversely, suppose we are given (\ref{Theorem1Oct4(1)}),
(\ref{Theorem1Oct4(2)}), (\ref{Theorem1Oct4(3)}). The equality of fields
asserted in \cite[Proposition 10]{BJKR98} is taken in the sense of ``equality
of $\mathbb{Q}[\lambda_{(A)}]$ and $\mathbb{Q}[\lambda_{(B)}]$ as subfields of
the real numbers'', which gives embeddings of $\mathbb{Z}[\lambda_{(A)}]$ and
$\mathbb{Z}[\lambda_{(B)}]$ into the real numbers. Thus it also embeds the
modules which can be considered as subsets of $\mathbb{Q}[\lambda_{(B)}]$. The
isomorphism of modules as additive groups acted on multiplicatively (i.e., the
action $y\mapsto xy$) by subrings of $\mathbb{Q}[\lambda_{(A)}]$ having full
rank (in this case $N$) means there is some element of the quotient field
mapping one to the other: if the isomorphism of (\ref{Theorem1Oct4(3)}) maps
some element $y$ to $z$ (considered as images in the real numbers) then the
ratio $y/z$ is independent of the choice of $y$ by definition of the
isomorphism in (\ref{Theorem1Oct4(3)}), and we multiply by this ratio to get
the isomorphism in (\ref{Theorem1Oct4(4)}).
\end{proof}

Note that this applies in particular to Example \ref{ExaSom.5} below.

\begin{theorem}
\label{Theorem2Oct4}Let $A$ and $B$ be matrices over $\mathbb{Z}_{+}$, all of
whose eigenvalues are rational, and each of which is divisible by some prime
not dividing the other eigenvalues. Assume further that $A$ and $B$ have the
same characteristic polynomial. Let $E_{a}$ and $E_{b}$ be their matrices of
column eigenvectors normalized to be integer vectors having greatest common
divisor $1$. Let $D$ be a diagonal matrix whose entries involve only powers of
primes in the respective eigenvalues, let $D_{s}$ be a diagonal matrix
consisting of precisely the diagonal eigenvalues. Then the following are equivalent:

\begin{enumerate}
\item \label{Theorem2Oct4(1)}$A$ and $B$ are $C^{\ast}$-equivalent;

\item \label{Theorem2Oct4(2)}$A$ and $B$ are shift equivalent, as follows: for
some choice of signs in $E_{a}$, $E_{b}$, and some choice of $D$, and for all
sufficiently large $n$, $E_{a}^{{}}DD_{s}^{n}E_{b}^{-1}$ and $E_{b}^{{}}%
D^{-1}D_{s}^{n}E_{a}^{-1}$ are integer matrices.
\end{enumerate}
\end{theorem}

\begin{proof}
Consider an isomorphism of dimension groups. The eigenvectors generate the
$1$-dimensional spaces of vectors such that some multiples of those vectors
are in the dimension group and are divisible by arbitrary powers of the
respective eigenvalues. Hence any dimension group isomorphism must preserve
those subspaces. Moreover, we claim that a dimension group isomorphism must
send normalized eigenvectors to normalized eigenvectors of the image. The
rational multiples of a normalized eigenvector $v$ with rational eigenvalue
$\eta$ which lie in the dimension group are the elements of $M_{v}%
=\{(n/m)v\mid n\in\mathbb{Z},\;\exists\,k\in\mathbb{Z}_{+},m|\eta^{k}\}$: a
vector $w\in M_{v}$ lies in the dimension group $G(A)$, since $w$ is a
multiple of an integer vector by negative powers of $A,B$. And if $w=(n/m)v\in
G(A)$ then there is a $k\in\mathbb{Z}$ such that $A^{k}v=\eta^{k}%
v\in\mathbb{Z}^{N}$, so that $m|\eta^{k}$ in lowest terms, $w\in M_{v}$. It
follows that dimension group isomorphism implies the existence of an
isomorphism of $\mathbb{Q}^{N}$ which sends each eigenvector to a multiple of
the other eigenvector by a number which divides a power of $\eta$. Such a
mapping must preserve the action of multiplication by $A$, given that the
characteristic polynomials are equal, because this multiplies each eigenvector
by its eigenvalue, and the eigenvalues are the same. So the mapping will be a
shift equivalence. Let $D$ be the diagonal matrix whose main diagonal entries
are the multiples just mentioned. Then the isomorphism $J$ of dimension groups
will be, specifically, $E_{b}(E_{a}D)^{-1}$ if it exists. For if we multiply
$J$ and its inverse on the left by a large enough power of $A$ or $B$,
respectively, as in (\ref{eqInt.18}), we see that the resulting matrix
products must be integer matrices. Moreover, these multiples are the matrices
stated in the theorem.
\end{proof}

\section{\label{Som}The transpose map and $C^{\ast}$-symmetry}

\addtolength{\textheight}{1.0625\baselineadjust}In this section we will study
the behavior of the dimension group $\left(  G\left(  A\right)  ,G\left(
A\right)  _{+}\right)  $ under the transpose map $A\rightarrow
A^{\operatorname*{tr}}$. In particular, we say that $A$ is $C^{\ast}%
$-symmetric if $A$ is $C^{\ast}$-equivalent to $A^{\operatorname*{tr}}$, i.e.,
$G\left(  A\right)  $ and $G\left(  A^{\operatorname*{tr}}\right)  $ are
isomorphic as ordered groups. We give several examples showing that $A$ may be
$C^{\ast}$-symmetric, (\ref{eqSom.1}), Remark \ref{RemSomOct.3}, or not,
Example \ref{ExaSomNew.6} ($2\times2$ matrices with rational eigenvalues),
Example \ref{ExaSom.5} ($2\times2$ matrices with irrational eigenvalues) and
Example \ref{Exa2}. An interesting feature with these particular examples is
that when $A$ is a $2\times2$ matrix, then $C^{\ast}$-symmetry is equivalent
to shift-symmetry (i.e., $A$ and $A^{\operatorname*{tr}}$ are shift
equivalent). For $2\times2$ matrices, symmetry seems to be more common than
non-symmetry. Our first example, while very simple, illustrates both $C^{\ast
}$-symmetry and a nontrivial $\operatorname*{Ext}$-element. It has
$\lambda=\lambda_{\left(  A\right)  }=2$. The $\operatorname*{Ext}$-group
represents another contrast between the two cases, $\lambda$ rational (and
hence integral), and the characteristic polynomial irreducible. In the first
case, we generally have $\ker\left(  \tau_{A}\right)  \neq0$, and as we note
in Remark \ref{RemSomOct.3}, $\operatorname*{ran}\left(  \tau_{A}\right)
=\mathbb{Z}\left[  1/\lambda\right]  $. Hence this extra extension structure
for $G\left(  A\right)  $ arises only in the reducible case: The corresponding
short exact sequence%
\begin{equation}
0\longrightarrow\ker\left(  \tau_{A}\right)  \longrightarrow G\left(
A\right)  \overset{\tau_{A}}{\longrightarrow}\mathbb{Z}\left[  1/\lambda
\right]  \longrightarrow0 \label{eqSomNew.J1}%
\end{equation}
may be non-split, which means that $G\left(  A\right)  $ is then not the
direct sum of the two groups $\ker\left(  \tau_{A}\right)  $ and
$\mathbb{Z}\left[  1/\lambda\right]  $.

Recall that for groups $R$ and $S$, $\operatorname*{Ext}\left(  R,S\right)  $
is again a group. Elements in $\operatorname*{Ext}\left(  R,S\right)  $ are
equivalence classes of short exact sequences%
\begin{equation}
0\longrightarrow S\longrightarrow E\overset{\psi}{\longrightarrow
}R\longrightarrow0, \label{eqSomNew.2}%
\end{equation}
see \cite{CaEi56}, \cite[p.~62]{BJO99}, and there are natural operations%
\[
E\longmapsto-E
\]
and%
\[
E,E^{\prime}\longmapsto E+E^{\prime}%
\]
on $\operatorname*{Ext}\left(  R,S\right)  $ which turn it into a group. We
say that (\ref{eqSomNew.2}) splits if there is some $\varphi\in
\operatorname*{Hom}\left(  R,E\right)  $ such that $\psi\circ\varphi
=\operatorname*{id}_{R}$. Then (\ref{eqSomNew.2}) splits if and only if it
represents the zero element in $\operatorname*{Ext}\left(  R,S\right)  $. For
example,%
\[
\operatorname*{Ext}\left(  \mathbb{Z}_{2},\mathbb{Z}\right)  \cong
\mathbb{Z}_{2},
\]
where the non-split $\operatorname*{Ext}$-element is represented by%
\[
0\longrightarrow2\mathbb{Z}\longrightarrow\mathbb{Z}\overset{\text{proj}%
}{\longrightarrow}\mathbb{Z}_{2}\longrightarrow0
\]
and the zero element is represented by%
\[
0\longrightarrow\mathbb{Z}\longrightarrow\mathbb{Z}\oplus\mathbb{Z}%
_{2}\overset{0\oplus\operatorname*{id}}{\longrightarrow}\mathbb{Z}%
_{2}\longrightarrow0.
\]

The next example illustrates how nontrivial $\operatorname*{Ext}$-elements are
part of the invariant structure for the dimension groups.

\begin{example}
\label{Exa1} The dimension group defined by $A$ may be order isomorphic to
that defined by its transpose $B=A^{\operatorname*{tr}}$. Hence an
AF-$C^{\ast}$-algebra built on such a matrix $A$ \textup{(}i.e., from the
corresponding stationary Bratteli diagram\/\textup{)} has a nontrivial
period-two symmetry corresponding to $A\mapsto A^{\operatorname*{tr}}$. An
example here is%
\begin{equation}
A=%
\begin{pmatrix}
1 & 1\\
2 & 0
\end{pmatrix}
,\qquad A^{\operatorname*{tr}}=%
\begin{pmatrix}
1 & 2\\
1 & 0
\end{pmatrix}
. \label{eqSom.1}%
\end{equation}
In this case $A$ and $A^{\operatorname*{tr}}$ have eigenvalues $2$ and $-1$,
and both of the dimension groups $G$ and $G^{\operatorname*{tr}}$ are in
$\operatorname*{Ext}\left(  \mathbb{Z}\left[  1/2\right]  ,\mathbb{Z}\right)
$. It can be checked \textup{(}by use of \cite[Corollary 11.28]{BJO99}%
\textup{)} that this $\operatorname*{Ext}$-element is not zero. Here
$\ker\left(  \tau\right)  =\mathbb{Z}$, $\operatorname*{ran}\left(
\tau\right)  =\mathbb{Z}\left[  1/2\right]  $, and the corresponding short
exact sequence%
\begin{equation}
0\longrightarrow\mathbb{Z}\longrightarrow G\left(  A\right)  \overset{\tau
}{\longrightarrow}\mathbb{Z}\left[  1/2\right]  \longrightarrow0
\label{eqSomNew.J2}%
\end{equation}
does not split, i.e., it is not the zero element in $\operatorname*{Ext}$.
Equivalently, $G\left(  A\right)  $ is not $\mathbb{Z}\oplus\mathbb{Z}\left[
1/2\right]  $ as a group. If it were, we would get $\tau\left(  w\right)
^{-1}\in\mathbb{Z}\left[  1/2\right]  $ by \cite[Corollary 11.28]{BJO99}. But
we computed $\tau\left(  w\right)  =3$, and $1/3$ is not in $\mathbb{Z}\left[
1/2\right]  $. Since $\lambda_{\left(  A\right)  }=2=\left|  \det A\right|  $,
it is tempting to apply Theorem \textup{\ref{ProFred.2.7.99}.} In fact the
inner-product invariants are $%
\ip{v}{w}%
=%
\begin{pmatrix}
2 & 1
\end{pmatrix}
\smash[b]{
\begin{pmatrix}
1\\ 1
\end{pmatrix}
}=3$, and $%
\ip{v^{\prime}}{w^{\prime}}%
=%
\begin{pmatrix}
1 & 1
\end{pmatrix}%
\begin{pmatrix}
2\\
1
\end{pmatrix}
=3$. But since the dimension is $2$ \textup{(}$<3$\textup{),} Theorem
\textup{\ref{ProFred.2.7.99}} does not apply directly, and instead we will
verify directly that $A$ and $A^{\operatorname*{tr}}$ are $C^{\ast}%
$-equivalent. Define matrices $J$, $K$ by
\begin{equation}
J=%
\begin{pmatrix}
1 & 1\\
1 & 0
\end{pmatrix}
,\qquad K=%
\begin{pmatrix}
1 & 0\\
0 & 2
\end{pmatrix}
. \label{eqSomOct.2}%
\end{equation}
One verifies that
\begin{equation}
A=KJ,\qquad A^{\operatorname*{tr}}=JK. \label{eqSomOct.3}%
\end{equation}
Thus $A$ and $A^{\operatorname*{tr}}$ are elementary shift equivalent, and it
follows that they are shift equivalent and $C^{\ast}$-equivalent \textup{(}see
the discussion in \cite{BJKR98}\textup{).}

However, we will see in Examples \textup{\ref{Exa2}} and
\textup{\ref{ExamplebisJul28}} that this is not a general feature of the
transpose map.
\end{example}

We may analyze the $C^{\ast}$-symmetry question by dimension-group analysis:
If we show that the ordered group $G\left(  A\right)  $ is order isomorphic to
$G\left(  A^{\operatorname*{tr}}\right)  $, then $A$ is $C^{\ast}$-equivalent
to $A^{\operatorname*{tr}}$, i.e., $A$ is $C^{\ast}$-symmetric. Clearly then
the two groups $G\left(  A\right)  $ and $\operatorname*{ran}\left(  \tau
_{A}\right)  $ are order isomorphic whenever $\ker\left(  \tau_{A}\right)
=0$, and we have the result:

\begin{proposition}
\label{ProSomOct.2}Let $A\in M_{N}\left(  \mathbb{Z}\right)  $ be nonsingular
and primitive, and suppose its characteristic polynomial $p_{\left(  A\right)
}\left(  x\right)  $ is irreducible, and $\operatorname*{ran}\left(  \tau
_{A}\right)  =\operatorname*{ran}\left(  \tau_{A^{\operatorname*{tr}}}\right)
$: then $A$ is $C^{\ast}$-equivalent to $A^{\operatorname*{tr}}$. Note in
particular that this holds if:

\begin{enumerate}
\item \label{ProSomOct.2(1)}$N=2$,

\item \label{ProSomOct.2(2)}the Perron--Frobenius eigenvalue $\lambda_{\left(
A\right)  }$ is irrational, and

\item \label{ProSomOct.2(3)}$\operatorname*{ran}\left(  \tau_{A}\right)
=\operatorname*{ran}\left(  \tau_{A^{\operatorname*{tr}}}\right)  $.
\end{enumerate}
\end{proposition}

\begin{proof}
This follows directly from Theorem \ref{Theorem1Oct4}.
\end{proof}

\begin{remark}
\label{RemSom.pound}We saw that by scaling out denominators in the entries
$v_{i}$ of the left \textup{(}row\/\textup{)} Perron--Frobenius eigenvector
$v\left(  A\right)  =\left(  v_{1},\dots,v_{N}\right)  $ we can arrange that
$v_{i}\in\mathbb{Z}\left[  1/\lambda\right]  $ for all $i$. But then a further
scaling with a power of $\lambda$ we can get each $v_{i}$ in the subring
$\mathbb{Z}\left[  \lambda\right]  \subset\mathbb{Z}\left[  1/\lambda\right]
$. Suppose that the characteristic polynomial of $A$ is irreducible. Note
that, as a group, $\mathbb{Z}\left[  \lambda\right]  $ is then a copy of the
lattice $\mathbb{Z}^{N}$ so the entries $v_{i}$ may therefore be viewed as
vectors in $\mathbb{Z}^{N}$. Then pick $v\left(  A\right)  $ such that
$\gcd\left(  v_{i}\right)  =1$ for each $i$. In this case the matrix $V$ with
the $v_{i}$'s as rows is in $M_{N}\left(  \mathbb{Z}\right)  $ and is
nonsingular. If we could define greatest common divisors in the ring
$\mathbb{Z}\left[  \lambda\right]  $ then we could divide $v$ by this greatest
common divisor and obtain some new $v$ defined over $\mathbb{Z}\left[
\lambda\right]  $ which has g.c.d.\ $1$. Then the image of its trace would
contain the span of its coordinates $v_{i}$ over $\mathbb{Z}\left[
\lambda\right]  $, that is, the entire ring $\mathbb{Z}\left[  \lambda\right]
$. Moreover the image of the trace will be contained in this ring, so they are
equal. In general, however, this ring will not be a principal ideal domain, so
that the class of the ideal generated by the trace becomes an obstruction. In
fact, the subgroup in $\mathbb{Z}\left[  \lambda\right]  $ which is generated
by the $v_{i}$'s is also an ideal in $\mathbb{Z}\left[  \lambda\right]  $.
Indeed, for $m\in\mathbb{Z}^{N}$, $\sum_{i}m_{i}v_{i}=\tau\left(  m\right)  =%
\ip{v}{m}%
$, so $\lambda\sum_{i}m_{i}v_{i}=%
\ip{vA}{m}%
=%
\ip{v}{A^{\operatorname*{tr}}m}%
$, and $A^{\operatorname*{tr}}m\in\mathbb{Z}^{N}$. As a consequence, we get
that the special incidence matrices $A$ which we considered in \cite{BJO99}
satisfy the condition $\operatorname*{ran}\left(  \tau_{A}\right)
=\mathbb{Z}\left[  1/\lambda_{\left(  A\right)  }\right]  $. However, this
fails for the matrix $A$ from Example \textup{\ref{ExaSom.5},} and others. The
group $\tau\left(  \mathbb{Z}^{N}\right)  $ is contained in $\tau\left(
G\right)  =\operatorname*{ran}\left(  \tau_{A}\right)  $ and the following
proposition indicates their relationship.
\end{remark}

\begin{proposition}
\label{ProSom.Fred}Assume that the Perron--Frobenius row eigenvector $v$ is
chosen to lie in $\mathbb{Z}^{N}\left[  \lambda\right]  $. Then the map
induced by the inclusion $\mathbb{Z}\left[  \lambda\right]  /\tau\left(
\mathbb{Z}^{N}\right)  \rightarrow\mathbb{Z}\left[  1/\lambda\right]
/\tau\left(  G\right)  $ is an epimorphism with kernel precisely%
\[
\left\{  x\in\mathbb{Z}\left[  \lambda\right]  /\tau\left(  \mathbb{Z}%
^{N}\right)  \mid\exists\,m\in\mathbb{Z},\;\lambda^{m}x=0\right\}  .
\]
\end{proposition}

\begin{proof}
Assuming, as we will show, that the image consists of torsion elements
relatively prime to $\lambda$, the inclusion gives a natural mapping. If we
multiply any element in $\mathbb{Z}\left[  1/\lambda\right]  $ by a power of
power of $\lambda$, we can get an element of $\mathbb{Z}\left[  \lambda
\right]  $, so this mapping is an epimorphism. We also claim that if we
multiply any element of $\tau\left(  G\right)  $, say $vA^{-n}x$,
$x\in\mathbb{Z}^{N}$ by a power of $\lambda$, we will get an element of
$\tau\left(  \mathbb{Z}^{N}\right)  $. This is because $vA^{-n}x=v\lambda
^{-n}x$ using the left two factors.

Note that since $\tau\left(  G\right)  $ and $\mathbb{Z}\left[  1/\lambda
\right]  $ are both torsion-free and $\lambda$-divisible, their quotient has
no $\lambda$-torsion. Hence every element annihilated by a power of $\lambda$
lies in the kernel. (We will show that the image consists of torsion elements
relatively prime to $\lambda$.)

Let $y$ be in the kernel of this mapping. Then $y\in\tau\left(  G\right)  $,
so that for some $n\in\mathbb{Z}_{+}$, $\lambda^{n}y\in\tau\left(
\mathbb{Z}^{N}\right)  $ and is zero in the original group. This identifies
the quotient. The left hand group, the quotient of a free abelian group by a
full-rank subgroup, is finite, so some fixed $n$ works for the whole kernel.
\end{proof}

\begin{remark}
[Rational $\lambda$]\label{RemSomOct.3}Even if $N=2$, the dimension group
$G\left(  A\right)  $ is not yet completely understood \cite{BJO99}
\textup{(}perhaps far from it; see, however, \cite{Han87}\/\textup{).} If
$\lambda=\lambda_{\left(  A\right)  }$ is rational, and therefore an integer,
we can have nonisomorphic $G\left(  A_{1}\right)  $ and $G\left(
A_{2}\right)  $ even when $A_{1}$ and $A_{2}$ have the same characteristic
polynomial and thus the same Perron--Frobenius eigenvalue $\lambda$, as
different extensions, $i=1,2$,
\begin{equation}
0\longrightarrow\mathbb{Z}\left[  1/\mu\right]  \hooklongrightarrow G\left(
A_{i}\right)  \overset{\tau}{\longrightarrow}\mathbb{Z}\left[  1/\lambda
\right]  \longrightarrow0, \label{eqSomOct.secsymb}%
\end{equation}
i.e., as different elements of the group $\operatorname*{Ext}\left(
\mathbb{Z}\left[  1/\lambda\right]  ,\mathbb{Z}\left[  1/\mu\right]  \right)
$. Here $\mu$ is the other root of the characteristic polynomial, so $\mu$ is
a nonzero integer with $\left|  \mu\right|  <\lambda$. See
\textup{(\ref{eqSomOct.4})} and \textup{(\ref{eqSomOct.5}).} This may even
happen when $A_{2}$ is the transpose of the matrix $A_{1}$, by Example
\textup{\ref{ExaSomNew.6}} below. Since here $\lambda$ is rational one may
arrange that $\tau_{A}\left(  G\left(  A\right)  \right)  =\mathbb{Z}\left[
1/\lambda\right]  $ by choosing $v$ with $\gcd\left(  v\right)  =1$, and
$\ker\left(  \tau_{A}\right)  $ is a rank-$1$ nonzero group isomorphic to
$\mathbb{Z}\left[  1/\mu\right]  $; see also below. There are specimens of
$2\times2$ primitive matrices $A$, even with integral Perron--Frobenius
eigenvalue such that $\lambda_{\left(  A\right)  }<\left|  \det A\right|  $,
and yet the two groups $G\left(  A\right)  $ and $G\left(
A^{\operatorname*{tr}}\right)  $ are order isomorphic. For example,
$A=\left(
\begin{smallmatrix}
1 & 5\\
3 & 3
\end{smallmatrix}
\right)  $ has that property. To see this, we may use \textup{(\ref{eqInt.17}%
)--(\ref{eqInt.18}).} Since $J=\left(
\begin{smallmatrix}
1 & 2\\
2 & 3
\end{smallmatrix}
\right)  $ satisfies $JA=A^{\operatorname*{tr}}J$, the two conditions hold,
and hence the matrix $A$ is $C^{\ast}$-symmetric. So for this particular pair
$A$, $A^{\operatorname*{tr}}$, the respective groups $G\left(  A\right)  $ and
$G\left(  A^{\operatorname*{tr}}\right)  $ from the middle term in the diagram
\textup{(\ref{eqSomOct.secsymb})} will then in fact represent the same zero
element of $\operatorname*{Ext}\left(  \mathbb{Z}\left[  1/6\right]
,\mathbb{Z}\left[  1/2\right]  \right)  $. For this particular $A$,%
\begin{equation}
G\left(  A\right)  \cong\mathbb{Z}\left[  1/2\right]  \oplus\mathbb{Z}\left[
1/6\right]  \text{\qquad(}\ker\left(  \tau\right)  \cong\mathbb{Z}\left[
1/2\right]  \text{)} \label{eqSom.star}%
\end{equation}
as direct sum of abelian groups. For this, note that the integral column
eigenvectors for $A$ are $\left(
\begin{smallmatrix}
1\\
1
\end{smallmatrix}
\right)  $ and $\left(
\begin{smallmatrix}
5\\
-3
\end{smallmatrix}
\right)  $. Since $\det\left(
\begin{smallmatrix}
1 & 5\\
1 & -3
\end{smallmatrix}
\right)  =-8=-2^{3}$ and the eigenvalues of $A$ are $6=2\cdot3$ and $-2$, we
have $\mathbb{Z}\left[  1/2\right]  ^{2}\subset G\left(  A\right)  $. Thus
$G\left(  A\right)  =\bigcup_{n=0}^{\infty}A^{-n}(\mathbb{Z}\left[
1/2\right]  ^{2})$, and \textup{(\ref{eqSom.star})} follows. Specifically, the
representation \textup{(\ref{eqSom.star})} may be derived from
\textup{(\ref{eqIntNew.21}), (\ref{eqSomOct.6}),} and the two identities%
\begin{equation}
\ker\left(  \tau_{A}\right)  =V\left(  A\right)  \cap G\left(  A\right)
=\mathbb{Z}\left[  1/2\right]
\begin{pmatrix}
5\\
-3
\end{pmatrix}
\label{eqSom.J1}%
\end{equation}
and%
\begin{equation}
G\left(  A\right)  \cap\mathbb{C}w\left(  A\right)  =\mathbb{Z}\left[
1/6\right]  w\left(  A\right)  , \label{eqSom.J2}%
\end{equation}
where $w\left(  A\right)  =\left(
\begin{smallmatrix}
1\\
1
\end{smallmatrix}
\right)  $. The present computation of $G\left(  A\right)  $ is simplified by
the fact that the orthogonal complement of the trace vector $v\left(
A\right)  =\left(
\begin{smallmatrix}
3, & 5
\end{smallmatrix}
\right)  $ is spanned by the nonmaximal column eigenvector. Here
$\lambda_{\left(  A\right)  }=6$, and so $\mathbb{Z}\left[  \lambda_{\left(
A\right)  }\right]  =\mathbb{Z}$. That $\tau\left(  G\left(  A\right)
\right)  =\mathbb{Z}\left[  1/6\right]  $ in this case follows from Remark
\textup{\ref{RemSom.pound}} and the general observation that with our choice
of $v\left(  A\right)  $, we will have $\tau\left(  G\left(  A\right)
\right)  =\mathbb{Z}\left[  1/\lambda_{\left(  A\right)  }\right]  $ provided
the ideal in $\mathbb{Z}\left[  \lambda_{\left(  A\right)  }\right]  $
generated by the $v_{i}\left(  A\right)  $ entries is principal. Ideals in
$\mathbb{Z}$ are principal, of course. Here in this case the
$\operatorname*{Ext}$-element corresponding to $G\left(  A\right)  $ is
trivial. \textup{(}Looking at prime factors in $\det A$, one could also get a
$G\left(  A\right)  $ which is non-split. For example, taking $A=\left(
\begin{smallmatrix}
1 & 6\\
2 & 2
\end{smallmatrix}
\right)  $, we get the spectrum $\left\{  5,-2\right\}  $ and that the
corresponding dimension group $G\left(  A\right)  $ is here represented by a
nonzero element of $\operatorname*{Ext}\left(  \mathbb{Z}\left[  1/5\right]
,\mathbb{Z}\left[  1/2\right]  \right)  $. The analysis here is analogous to
that presented above: We get $\ker\left(  \tau_{A}\right)  \cong
\mathbb{Z}\left[  1/2\right]  $, $\operatorname*{ran}\left(  \tau_{A}\right)
\cong\mathbb{Z}\left[  1/5\right]  $, and the corresponding short exact
sequence%
\begin{equation}
0\longrightarrow\mathbb{Z}\left[  1/2\right]  \longrightarrow G\left(
A\right)  \longrightarrow\mathbb{Z}\left[  1/5\right]  \longrightarrow0
\label{eqSom.J8}%
\end{equation}
is now non-split. The example $A=\left(
\begin{smallmatrix}
1 & 6\\
2 & 2
\end{smallmatrix}
\right)  $ is $C^{\ast}$-symmetric, as $A$ and $A^{\operatorname*{tr}}$ are in
fact shift equivalent: Take $R=\left(
\begin{smallmatrix}
-1 & 2\\
2 & -2
\end{smallmatrix}
\right)  $ and $S=\left(
\begin{smallmatrix}
7 & 4\\
4 & 3
\end{smallmatrix}
\right)  $. Then $RS=A^{\operatorname*{tr}}$ and $SR=A$. It follows from
\cite{BJO99} that $G\left(  A\right)  $, when represented in the
$\operatorname*{Ext}$-group, is generally not the zero element.\textup{)}
\end{remark}

\begin{example}
\label{ExaSomNew.6}Here we will exhibit a primitive nonsingular $2\times2$
matrix $A$ with rational eigenvalues such that $A$ is not $C^{\ast}%
$-equivalent to $A^{\operatorname*{tr}}$ \textup{(}and thus is not shift
equivalent to $A^{\operatorname*{tr}}$\textup{).} The respective dimension
groups $G\left(  A\right)  $ and $G\left(  A^{\operatorname*{tr}}\right)  $
are not even isomorphic as groups, let alone order isomorphic, and hence this
$A$ in \textup{(\ref{eqSonNew.B1})} is ``more'' nonsymmetric than the
corresponding specimen \textup{(\ref{eq1Oct7.1})} in Example
\textup{\ref{ExaSom.5}.} The example is%
\begin{equation}
A=%
\begin{pmatrix}
65 & 7\\
24 & 67
\end{pmatrix}
. \label{eqSonNew.B1}%
\end{equation}
Putting
\begin{equation}
E_{A}=%
\begin{pmatrix}
-7 & 1\\
12 & 2
\end{pmatrix}
,\qquad E_{B}=%
\begin{pmatrix}
-2 & 12\\
1 & 7
\end{pmatrix}
,\qquad D=%
\begin{pmatrix}
53 & 0\\
0 & 79
\end{pmatrix}
, \label{eqSomNew.B2}%
\end{equation}
we have%
\begin{equation}
A=E_{A}^{{}}DE_{A}^{-1},\qquad B=A^{\operatorname*{tr}}=E_{B}^{{}}DE_{B}^{-1}.
\label{eqSomNew.B3}%
\end{equation}
The eigenvalues of $A$ and $B$ are $53$ and $79$, which are both prime and
congruent to $-1\mod{13}$. Using Theorem \textup{\ref{Theorem2Oct4}} it
follows that if $A$ and $B$ were $C^{\ast}$-equivalent there would exist some
diagonal matrix $D_{0}=\left(
\begin{smallmatrix}
x & 0\\
0 & y
\end{smallmatrix}
\right)  $ where $x$, $y$ are congruent to $\pm1\mod{13}$ such that $E_{A}%
^{{}}D_{0}^{{}}E_{B}^{-1}$ would have integral entries. But the $\left(
1,1\right)  $ entry of this matrix is $\left(  49x+y\right)  /26$. If this is
an integer, and $x=\varepsilon_{1}+n_{1}\cdot13$, $y=\varepsilon_{2}%
+n_{2}\cdot13$, where $n_{1}$, $n_{2}$ are integers and $\varepsilon_{i}=\pm
1$, then $\frac{1}{13}\left(  49x+y\right)  =\frac{1}{13}\left(  \left(
4\cdot13-3\right)  x+y\right)  =-\varepsilon_{1}\frac{3}{13}+\varepsilon
_{2}\frac{1}{13}\mod{1}$, but this can never be an integer. Thus $A$ is not
$C^{\ast}$-symmetric.
\end{example}

\begin{remark}
[Irrational $\lambda$]\label{RemSomOct.pound}The assumption in Proposition
\textup{\ref{ProSomOct.2}} that the range of the respective traces $\tau_{A}$
and $\tau_{A^{\operatorname*{tr}}}$ be the same \textup{(}viewed as subgroups
of $\mathbb{Z}\left[  1/\lambda_{\left(  A\right)  }\right]  $\textup{)}
cannot be omitted. It is true in general that $\operatorname*{ran}\left(
\tau_{A}\right)  $ is an ideal in $\mathbb{Z}\left[  1/\lambda_{\left(
A\right)  }\right]  $, but the ideal may be proper, and it may be different
from one to the other. An example showing this to be the case can be found in
\cite[p.~104]{BMT87}, \cite[pp.~79--83]{PaTu82}. The example is a matrix $A$
such that $A$ and its transpose $B=A^{\operatorname*{tr}}$ are not shift
equivalent. We will give another example of this, and then apply Theorem
\textup{\ref{Theorem1Oct4}} to show that they are not $C^{\ast}$-equivalent either:
\end{remark}

\begin{example}
\label{ExaSom.5}The example is $A=\left(
\begin{smallmatrix}
19 & 5\\
4 & 1
\end{smallmatrix}
\right)  $. Here $\lambda=10+\sqrt{101}$, so the characteristic polynomial is
irreducible and therefore $\ker\left(  \tau_{A}\right)  =0$. Since $\det
A=-1$, the unordered dimension groups $G\left(  A\right)  $ and $G\left(
A^{\operatorname*{tr}}\right)  $ are both $\mathbb{Z}^{2}$. However, we will
show that they are not order isomorphic. We have%
\begin{equation}
A=%
\begin{pmatrix}
19 & 5\\
4 & 1
\end{pmatrix}
,\qquad B=A^{\operatorname*{tr}}=%
\begin{pmatrix}
19 & 4\\
5 & 1
\end{pmatrix}
. \label{eq1Oct7.1}%
\end{equation}
We prove that the two ideals $\operatorname*{ran}\left(  \tau_{A}\right)  $
and $\operatorname*{ran}\left(  \tau_{A^{\operatorname*{tr}}}\right)  $ are
nonisomorphic. The eigenvalues are $10\pm\sqrt{101}$. Let $\omega
=(1+\sqrt{101})/2$ so that $1,\omega$ form a $\mathbb{Z}$-basis for the
algebraic integers in $\mathbb{Q}(\sqrt{101})$. \textup{(}The fact that all
algebraic integers in a quadratic field have this form is \cite[Theorem 6-1-1,
p.~234]{Wei98}. One can check that $1$, $\omega$ are algebraic integers, then
that the trace must be an algebraic integer, and see what happens when we
subtract some $a+b\omega$ to simplify, in terms of the norm being an algebraic
integer.\textup{)} The respective \textup{(}column\/\textup{)} eigenvectors
for $A$, $A^{\operatorname*{tr}}$ are
\begin{equation}%
\begin{pmatrix}
4+\omega\\
2
\end{pmatrix}
,\;%
\begin{pmatrix}
5-\omega\\
2
\end{pmatrix}
\text{\quad and\quad}%
\begin{pmatrix}
2\\
-5+\omega
\end{pmatrix}
,\;%
\begin{pmatrix}
-2\\
4+\omega
\end{pmatrix}
. \label{eq1Oct7.2}%
\end{equation}
By transposing and interchanging the two, we get as Perron--Frobenius row
eigenvectors for $A$, $A^{\operatorname*{tr}}$
\begin{equation}%
\begin{pmatrix}
2, & \omega-5
\end{pmatrix}
,\qquad%
\begin{pmatrix}
\omega+4, & 2
\end{pmatrix}
. \label{eq1Oct7.3}%
\end{equation}
Let $I_{1}$, $I_{2}$ denote the ideals they generate. We note that
$\omega-5=\left(  -9+\sqrt{101}\right)  /2$ and that the norm of this number
is $(81-101)/4=-5$. Hence, over the algebraic number ring, which is
$\mathbb{Z}[\omega]$ and properly contains $\mathbb{Z}[1/\lambda]$, both
ideals are the entire ring $(1)=\mathbb{Z}\left[  \omega\right]  $, since the
two have isomorphic spans over the algebraic number ring. Thus, we need to see
whether some element in the quotient field will multiply one ideal to the
other, as additive groups, or modules over $\mathbb{Z}[1/\lambda
]=\mathbb{Z}[\lambda]=\{a+b\sqrt{101}\mid a,b\in\mathbb{Z}\}$.

Note that the two generators listed in \textup{(\ref{eq1Oct7.3})} will
actually generate each ideal over $\mathbb{Z}$ additively, not just as modules
over $\mathbb{Z}[\lambda]$, since multiplication by $\sqrt{101}=2\omega-1$
sends
\begin{equation}%
\begin{pmatrix}
1, & \omega
\end{pmatrix}
\longrightarrow%
\begin{pmatrix}
2\omega-1, & 2\left(  \omega+25\right)  -\omega
\end{pmatrix}
=%
\begin{pmatrix}
2\omega-1, & \omega+50
\end{pmatrix}
, \label{eq3Oct7.4}%
\end{equation}%
\begin{equation}%
\begin{pmatrix}
2, & \omega-5
\end{pmatrix}
\longrightarrow%
\begin{pmatrix}
4\omega-2, & \left(  \omega+50\right)  -5\left(  2\omega-1\right)
\end{pmatrix}
=%
\begin{pmatrix}
4\omega-2, & 55-9\omega
\end{pmatrix}
, \label{eq3Oct7.5}%
\end{equation}%
\begin{equation}%
\begin{pmatrix}
4+\omega, & 2
\end{pmatrix}
\longrightarrow%
\begin{pmatrix}
\left(  \omega+50\right)  +4(2\omega-1), & 4\omega-2
\end{pmatrix}
=%
\begin{pmatrix}
9\omega+46, & 4\omega-2
\end{pmatrix}
, \label{eq3Oct7.6}%
\end{equation}
which are still in the same additive subgroups.

The additive spans of the two pairs of generators in \textup{(\ref{eq1Oct7.3}%
)} are, respectively
\begin{align}
2\mathbb{Z}+\left(  \omega-5\right)  \mathbb{Z}  &  =\left\{  a+b\omega\mid
a,b\in\mathbb{Z}\text{ such that }a-b\equiv0\pmod{2}\right\}
,\label{eq14Oct7.14}\\
2\mathbb{Z}+\left(  \omega+4\right)  \mathbb{Z}  &  =\left\{  a+b\omega\mid
a,b\in\mathbb{Z}\text{ such that }a-b\equiv1\pmod{2}\right\}  .
\label{eq14Oct7.15}%
\end{align}
These are preserved by multiplication by $\sqrt{101}=2\omega-1\equiv1\pmod{2}%
$, so that each span over $\mathbb{Z}$ is a $\mathbb{Z}\left[  \lambda\right]
$-module.

We will now complete the proof. If the ideals were isomorphic under
multiplication by some $f\in\mathbb{Q}\left[  \lambda\right]  $, then $f$
cannot involve primes of the algebraic number ring, since both ideals span the
complete algebraic number ring as modules over it. Therefore $f$ is a unit.
Thus $f$ is up to a sign a power of $\lambda=9+2\omega\equiv1\pmod{2}$. Hence
multiplication by $f$ preserves the congruence conditions defining the two
additive spans, and thus it preserves each ideal separately. So it is
impossible for a unit to send one ideal to the other.
\end{example}

\begin{example}
\label{Exa2} The following is an example of integer matrices $A$, $B$ which
have isomorphic dimension groups (unordered) but such that the corresponding
transposed matrices $A^{\operatorname*{tr}}$, $B^{\operatorname*{tr}}$ do not
have isomorphic dimension groups. Informally, for matrices over $\mathbb{Z}$
with eigenvalues $1$, $p$, $q$, the isomorphism type of the dimension group is
determined by the way the $p$-divisible and $q$-divisible integer vectors lie,
as subspaces, within all integer vectors.
\[%
\begin{pmatrix}
1 & 0 & 0\\
x & p & 0\\
y & 0 & q
\end{pmatrix}
\]
The $p,q$-divisible spaces split off as spanned by $\left(  0,1,0\right)  $,
$\left(  0,0,1\right)  $, hence the column dimension groups are the same no
matter what $x$, $y$ are. We can use the diagonal as a comparison. The row
dimension groups depend on how the $p$ and $q$ eigenvectors lie in their sum.
Row eigenvectors are in the row kernels of
\[%
\begin{pmatrix}
1-p & 0 & 0\\
x & 0 & 0\\
y & 0 & q-p
\end{pmatrix}
,\qquad%
\begin{pmatrix}
1-q & 0 & 0\\
x & p-q & 0\\
y & 0 & 0
\end{pmatrix}
.
\]
These vectors are spanned by $\left(  x,p-1,0\right)  $ and $\left(
y,0,q-1\right)  $, respectively. If $m$ divides $p-1$, $q-1$ and $x-y$ but not
$x$, $y$ then we have sum vectors which are divisible. This will give
non-isomorphism of row dimension groups.

More formally, the column dimension groups consist of all vectors of the form
$\left(  a,b/p^{n},c/q^{n}\right)  $ (we may consider the form of the inverse
and its columns). Any isomorphism on row dimension groups will be an integer
matrix which preserves the row and column divisible eigenspaces. Take the case
$p=3$, $q=5$, $x=y=1$ as compared with $x=y=0$. Then the row vectors $\left(
0,1,0\right)  $, $\left(  0,0,1\right)  $ for the diagonal case must be mapped
to multiples by powers of $3$, $5$, $-1$ of vectors $\left(  1,2,0\right)  $,
$\left(  1,0,4\right)  $, and the total determinant of the matrix must be odd.
This should be a $2$-adic isomorphism, but it cannot be since $\left(
0,1,0\right)  -\left(  0,0,1\right)  $ is not divisible by $2$, whereas
$\left(  1,2,0\right)  -\left(  1,0,4\right)  $ is divisible by $2$, and the
same is true if they are replaced by any of their odd multiples.
\end{example}

\begin{remark}
\label{Rem3} It follows from Theorem 3.1 of Boyle and Handelman \cite{BoHa93}
that there are nonnegative integer matrices that are shift equivalent to the
pair in Example \textup{\ref{Exa2}} and hence have the same dimension groups.
We will now construct an example where the ordered dimension groups are
isomorphic for two matrices, but for the transpose matrices, the ordered
dimension groups are not isomorphic.

We modify Example \ref{Exa2} a little bit so as to get an example of two
nonnegative integer matrices having identical ordered column dimension groups
but non-isomorphic row dimension groups. We start with the same matrices as
before except we take $3$ and $7$ as the two main diagonal primes. Then we add
a large odd prime eigenvalue $101$ which will enable a conjugate to be
positive. This gives matrices $A$, $B$:%
\[
A=%
\begin{pmatrix}
101 & 0 & 0 & 0\\
0 & 1 & 0 & 0\\
0 & 1 & 3 & 0\\
0 & 1 & 0 & 7
\end{pmatrix}
,\qquad B=%
\begin{pmatrix}
101 & 0 & 0 & 0\\
0 & 1 & 0 & 0\\
0 & 0 & 3 & 0\\
0 & 1 & 0 & 7
\end{pmatrix}
.
\]
We multiply these matrices by $5$ and then conjugate by the matrix $C$:%
\[
C=%
\begin{pmatrix}
2 & -1 & 0 & 0\\
1 & 1 & -1 & 0\\
1 & 0 & 2 & -1\\
1 & 0 & -1 & 1
\end{pmatrix}
,
\]
whose determinant is $5$ and which approximately moves the row and column
eigenvectors to $\left(  1,1,1,1\right)  $ to get nonnegative matrices
\[%
\begin{pmatrix}
205 & 200 & 200 & 200\\
104 & 109 & 94 & 94\\
110 & 110 & 105 & 70\\
86 & 86 & 106 & 141
\end{pmatrix}
,\qquad%
\begin{pmatrix}
205 & 200 & 200 & 200\\
103 & 108 & 93 & 93\\
111 & 111 & 106 & 71\\
86 & 86 & 106 & 141
\end{pmatrix}
.
\]
These have identical dimension groups: the original unordered dimension groups
are the same, by the above, just adding the eigenspace for the eigenvalue
$101$ in both cases. The order structure is determined by the
Perron--Frobenius row eigenvalue, which is $\left(  1,0,0,0\right)  C^{-1}$ in
both cases.

But for the transposes, row dimension groups, we have the direct sum of the
$101$ eigenspace with $\mathbb{Z}[1/5]$ times the previous examples. Making
the prime $5$ invertible will not affect the above argument that the row
dimension groups are not isomorphic because this was a $2$-adic
nonisomorphism, given that the row eigenspaces for $5\cdot3$ and $5\cdot7$
eigenvalues must be preserved by an isomorphism. These are the spaces of
vectors in the dimension group divisible by arbitrarily high powers of $3$ and
$7$.
\end{remark}

\section{\label{GZN}The quotient $G/\mathbb{Z}^{N}$ is an invariant}

Recall that $G=G\left(  A\right)  =\bigcup_{n=0}^{\infty}A^{-n}\mathbb{Z}^{N}%
$. In this section we will consider the quotient group $G/\mathbb{Z}^{N}$.
Here $\mathbb{Z}^{N}$ can be replaced with any free abelian subgroup $L$ of
$G$ such that%
\begin{equation}
G\left(  A\right)  =\bigcup_{n=0}^{\infty}A^{-n}L \label{eqGZN.1}%
\end{equation}
and%
\begin{equation}
AL\subset L. \label{eqGZN.2}%
\end{equation}

We used the quotient group in \cite{BJO99}, but at the time we did not know if
it was an invariant, and what the isomorphism classes were (in the category of
abelian torsion groups). These issues are now resolved in the next
proposition, which implies that the quotient is indeed an isomorphism
invariant, i.e., that a given $C^{\ast}$-isomorphism implies that the
corresponding two quotients are isomorphic groups.

Abelian torsion groups are classified in general by the so-called Ulm
invariant \cite[pp.~26--27]{Kap69}, \cite[and references given there]{KaMa51}.
The Ulm invariant in general is a sequence of natural numbers suitably indexed
by ordinals. These numbers are calculated as dimensions of certain vector
spaces over the field $\mathbb{F}_{p}=\mathbb{Z}/p\mathbb{Z}$. First, any
given torsion group decomposes over its $p$-subgroups, and the Ulm dimensions
are then calculated for each ordinal, when $p$ is fixed. In the present
application, the Ulm invariant is, as we show, very simple and concrete.

\begin{example}
\label{ExaGZNNew.1} \textit{The quotient group }$G/\mathbb{Z}^{N}$\textit{ for
the special case of Section} \ref{Cas}. It is easy to understand concretely
the torsion group quotient for the special case of Section \ref{Cas} when it
is assumed that $\left|  \det A\right|  =\lambda_{\left(  A\right)  }$. Of
course then $\lambda_{\left(  A\right)  }$ is an integer, and we may therefore
form $\mathbb{Z}\left[  1/\lambda_{\left(  A\right)  }\right]  $ the usual way
as an inductive limit $\bigcup_{n=1}^{\infty}\mathbb{Z}_{\lambda_{\left(
A\right)  }^{n}}$ as described in (\ref{eqpad.3}) with natural embeddings
$\mathbb{Z}_{\lambda_{\left(  A\right)  }^{n}}\hooklongrightarrow
\mathbb{Z}_{\lambda_{\left(  A\right)  }^{n+1}}\,$, and it follows from the
discussion in Section \ref{pad} and Section \ref{Cas} (Figure \ref{FigCas.1})
that there is then a natural isomorphism between the two groups $G\left(
A\right)  /\mathbb{Z}^{N}$ and $\mathbb{Z}\left[  1/\lambda_{\left(  A\right)
}\right]  /\mathbb{Z}$. Hence, in this very special case,
$\operatorname*{Prim}\left(  \lambda_{\left(  A\right)  }\right)  $ is a
complete invariant for the corresponding torsion group quotient. See also
\cite{BJO99} for more details. It is the case of dimension groups $G\left(
A\right)  $ more general than that of Section \ref{Cas} which requires a
nontrivial localization. The next two propositions deal with the general case,
and the appropriate localizations.
\end{example}

One method of localizing at a prime $p$ is to take the tensor product of an
abelian group with $\mathbb{Z}\left[  1/2,1/3,\dots,\widehat{1/p}%
,\dots\right]  $, inverting all primes except $p$; another is to tensor with
the $p$-adic integers. Both agree for all torsion groups; the latter
localization factors through the former. These tensor products are exact
functors of abelian groups $A$ which are subgroups of $\mathbb{Q}^{N}$, that
is, they preserve exact sequences; this follows from \cite[Proposition 7.2,
p.~138]{CaEi56}, since the group $D(A)$ is a direct sum of copies of
$\mathbb{R}/\mathbb{Z}$ which has no nontrivial continuous homomorphisms into
the totally disconnected $p$-adics. Thus $\operatorname*{Tor}^{1}(A,C)$ is
zero. Unless otherwise specified we will mean the former, smaller tensor
product when we localize.

\begin{proposition}
\label{Proposition1} Assume that $A$ is a nonsingular $N\times N$ integer
matrix and $G=\bigcup_{n}A^{-n}\mathbb{Z}^{N}$ the associated dimension group.
Suppose $L$ is any rank-$N$ lattice in $G$ such that $\bigcup_{n}A^{-n}L=G$.
Then $G/L$ is isomorphic to the product over primes dividing $\det A$ of a
certain number $n(p)$ copies of $\mathbb{Z}_{p^{\infty}}$. The number $n(p)$
is the largest $j$ such that if we write the characteristic polynomial of $A$
as $x^{N}+c_{1}x^{N-1}+\dots+c_{N-1}x+c_{N}$, $p$ does not divide $c_{j}$. In
addition $G$ is dual to the eventual row space in the sense
\begin{equation}
G\otimes\mathbb{Z}\left[  1/2,\dots,\widehat{1/p},\dots\right]  =\left\{
v\in\mathbb{Q}^{N}\mid%
\ip{w}{v}%
\in\mathbb{Z}_{(p)}\;\forall\,w\in G_{\left(  p\right)  }\left(  A\right)
\right\}  . \label{eqGZNNew.2}%
\end{equation}
\end{proposition}

\begin{proof}
$G$ written as $\bigcup_{n=0}^{\infty}A^{-n}\left(  \mathbb{Z}^{N}\right)  $
will have as denominators only primes dividing $\det A$. If $L$ includes
$\mathbb{Z}^{N}$ then we have a torsion group whose torsion involves only
primes in $\det A$.

We first argue that locally at each prime $p$ in it, $G$ consists of those
vectors dual to the eventual $p$-adic row space $G_{\left(  p\right)  }\left(
A\right)  $ of $A$ (see (\ref{eqpad.11})--(\ref{eqpad.12})). That is,
(\ref{eqGZNNew.2}) holds. The dimension group is the group of vectors $x$ such
that for some $n$, $A^{n}x\in\mathbb{Z}^{N}$. This is the group of vectors
such that $\exists\,n\in\mathbb{Z}_{+}$ such that for $\forall\,w\in
\mathbb{Z}^{N}$, we have $wA^{n}x\in\mathbb{Z}$. This is the group of rational
vectors whose products with the row space of $A^{n}$ is integer. This
construction also goes through if we localize at any prime. To say that a
vector has $p$-integer product with the row space of $A^{n}$ for some $n$ then
implies that it has $p$-integer product with the idempotent $p$-adic limit
$E_{\left(  p\right)  }\left(  A\right)  $ of powers of $A$, defined in
(\ref{eqpad.7}) and mentioned in Theorem 7 of \cite{BJKR98}. Conversely
suppose it has $p$-integer product with the idempotent $p$-adic limit, then by
$p$-adic continuity, it must have $p$-integer product with some finite power.
This gives the claim.

Now to show that the quotient group at the prime $p$ is $p$-divisible, take a
$p$-adic dual basis to $G_{\left(  p\right)  }\left(  A\right)  $, which, like
any $p$-adic torsion-free module, must be a free module (the $p$-adic integers
are a principal ideal domain, and argue as with the ordinary integers).
Approximate these vectors $p$-adically by rational vectors $b_{i}$ using a
$p$-adic approximation theorem such as \cite[1-2-3, p.~8]{Wei98}, choosing
these rational vectors so that they give a $p$-adic basis. Take the free
abelian group $L_{1}$ generated by $b_{i}$. As soon as we have a lattice $L$
including $L_{1}$ and $\mathbb{Z}^{N}$, the $p$-adic dimension group consists
of a sum of copies of the $p$-adic integers corresponding to $L_{1}$ and a sum
of copies of the $p$-adic field corresponding to the remaining vectors (in the
null space of $E_{\left(  p\right)  }\left(  A\right)  $---we can take
additional basis vectors for it). When we divide by $L$, we are dividing out
by all the $L_{1}$ part $p$-adically, and by something isomorphic inside a
$p$-adic field in the rest, and the result will be $p$-divisible.

In fact, for any lattice $L$ such that $\bigcup_{n}A^{-n}L$ is the dimension
group, the quotient will be isomorphic to this, since multiplication by
$A^{-1}$ gives an isomorphism of pairs $\left(  \bigcup_{n}A^{-n}L,L\right)
\rightarrow\left(  \bigcup_{n}A^{-n}L,A^{-1}L\right)  $, and eventually this
lattice must be large enough.

Now consider the $p$-adic rank, in relation to the characteristic polynomial.
By Newton's method \cite[3-1-1, p.~74]{Wei98}, if the characteristic
polynomial has the given form, we can factor it over the $p$-adics as a
product of two polynomials, one of which is $x^{N-j}$ modulo $p$, and the
other of which has invertible constant term over the $p$-adics. We can put the
matrix into corresponding block form. The former part will be $p$-adically
nilpotent, and the null space will be its row space.
\end{proof}

\begin{corollary}
\label{CorGZNNew.3}Let $A$, $B$ be matrices as in Proposition
\textup{\ref{Proposition1}.} For any choice of lattices $L$, $L^{\prime}$ in
their dimension groups satisfying the hypotheses in Proposition
\textup{\ref{Proposition1},} $G\left(  A\right)  \cong G\left(  B\right)
\Rightarrow G\left(  A\right)  /L\cong G\left(  B\right)  /L^{\prime}$.
\end{corollary}

\begin{remark}
\label{RemGZNNew.3}As noted, our groups $G\left(  A\right)  $ are contained in
$\mathbb{R}^{N}$ \textup{(}even in $\mathbb{Q}^{N}$\textup{)} where $N$ is the
rank of $G\left(  A\right)  $. But it is clear that general lattices $L$ in
$\mathbb{R}^{N}$ are given by a choice of basis in $\mathbb{R}^{N}$ as a
vector space. Writing the vectors in a basis, equivalently the generators for
$L$, as column vectors, we note that the lattices $L$ may be viewed as, or
identified with, nonsingular real matrices. Making this identification, and
fixing the rank $N$, we further note that the containment $L\subset L^{\prime
}$, for two given lattices, holds if and only if there is some $C\in
M_{N}\left(  \mathbb{Z}\right)  $ such that we have the following matrix
factorization:
\begin{equation}
L=L^{\prime}C. \label{eqRemGZNNew.3.0}%
\end{equation}
There is a similar version of this for row spaces \textup{(}or lattices
defined from row vectors\/\textup{),} as well as a $p$-adic variation,
\emph{mutatis mutandis;} and we have already seen an instance of the latter in
\textup{(\ref{eqpad.10})--(\ref{eqpad.12}).}
\end{remark}

\begin{remark}
\label{RemGZNNew.4}We now show, using (\ref{eqRemGZNNew.3.0}), that the
conditions on $L$ from Proposition \textup{\ref{Proposition1}} are all
integrality conditions. There are three in all, and we proceed to spell them
out. If $A$ is given as usual, and if $G\left(  A\right)  $ is the
corresponding group, i.e., $\bigcup_{n=0}^{\infty}A^{-n}\left(  \mathbb{Z}%
^{N}\right)  $, then a lattice $L$ is a subgroup, i.e., $L\subset G\left(
A\right)  $, if and only if there is a natural number $n$ such that
\begin{equation}
A^{n}L\in M_{N}\left(  \mathbb{Z}\right)  . \label{eqRemGZNNew.3.1}%
\end{equation}
Some given lattice $L$ will satisfy the invariance property $A\left(
L\right)  \subset L$ if and only if the conjugate matrix $L^{-1}AL$ satisfies
\begin{equation}
L^{-1}AL\in M_{N}\left(  \mathbb{Z}\right)  . \label{eqRemGZNNew.3.2}%
\end{equation}
The further condition on $L$ that it is generating, i.e., that $\bigcup
_{n=0}^{\infty}A^{-n}\left(  L\right)  =G\left(  A\right)  $, holds if and
only if for some natural number $n$ we have
\begin{equation}
L^{-1}A^{n}\in M_{N}\left(  \mathbb{Z}\right)  . \label{eqRemGZNNew.3.3}%
\end{equation}
The three conditions should also be compared with \textup{(\ref{eqInt.18})}
from Section \textup{\ref{Int}.}

Now Proposition \textup{\ref{Proposition1}} applies when $G\left(  A\right)  $
is given and some lattice satisfies all three conditions
\textup{(\ref{eqRemGZNNew.3.1})--(\ref{eqRemGZNNew.3.3}),} and we get as a
corollary that if two lattices $L$ and $L^{\prime}$ both satisfy the
conditions, then the two torsion groups $G\left(  A\right)  /L$ and $G\left(
A\right)  /L^{\prime}$ are isomorphic groups.
\end{remark}

In the study of dimension groups, it is convenient to explicitly compute
certain extensions. Let $\mathbb{Z}_{p^{\infty}}$ denote the union of
$\mathbb{Z}_{p^{n}}$ under inclusion, a divisible $p$-torsion group whose
order $p$ subgroup has rank $1$. By standard theory \cite{CaEi56}, the
extension group $\operatorname*{Ext}(\mathbb{Z}_{p^{\infty}},\mathbb{Z})$ can
be computed using the exact sequence
\begin{equation}
0\rightarrow\mathbb{Z}\rightarrow\mathbb{Q}\rightarrow\mathbb{Q}%
/\mathbb{Z}\rightarrow0 \label{eqGZNNew.3}%
\end{equation}
as the cokernel of the map $\operatorname*{Hom}(\mathbb{Z}_{p^{\infty}%
},\mathbb{Q})\rightarrow\operatorname*{Hom}(\mathbb{Z}_{p^{\infty}}%
,\mathbb{Q}/\mathbb{Z})$; the former group is zero and the latter group is
$\operatorname*{Hom}(\mathbb{Z}_{p^{\infty}},\mathbb{Z}_{p^{\infty}})$. Every
$p$-adic integer gives a mapping in this group; we check this mapping is
one-to-one and onto, so that the $\operatorname*{Ext}$ group is the $p$-adic
integers. (To check ``onto'', note that we get every mapping $\mathbb{Z}%
_{p^{n}}\rightarrow\mathbb{Z}_{p^{n}}$ and take limits.) In general, we are
dealing with a direct sum of copies of these $\operatorname*{Ext}$ groups.

Next we look at the problem of isomorphism of dimension groups in a somewhat
different way, by showing that dimension groups can easily be computed as
extensions. In some cases this leads to a quick decision about whether two
dimension groups are isomorphic. However in the most general case, the problem
of deciding isomorphism given this extension structure seems to still require
the methods of Section \ref{Dec}. In view of Remark \ref{RemGZNNew.3} we need
only state the result for the case when the lattice $L$ is $\mathbb{Z}^{N}$.

\begin{corollary}
\label{Proposition3rev}As in Proposition \textup{\ref{Proposition1}}, consider
an unordered dimension group as an extension of $\mathbb{Z}^{N}$ by a
divisible torsion group $G/\mathbb{Z}^{N}$ whose structure, computed as in
Proposition \textup{\ref{Proposition1},} is a direct sum over $i$ of $n_{i}$
copies of $\mathbb{Z}_{p\left(  i\right)  ^{\infty}}=\mathbb{Z}\left[
1/p\left(  i\right)  \right]  /\mathbb{Z}$. The extension class in
$\operatorname*{Ext}^{1}(G/\mathbb{Z}^{N},\mathbb{Z}^{N})$ is an element of
$\bigoplus_{i}\bigoplus_{j=1}^{n_{i}}\mathbb{Z}_{p(i)}$. We write this as an
$N\times\sum_{i}n_{i}$ matrix whose entries are $p\left(  i\right)  $-adic
integers:
\[%
\setlength{\unitlength}{24pt}
\begin{picture}(10.5,2.5)(0.5,0)
\put(1,0){\line(1,0){10}}
\put(1,2){\line(1,0){10}}
\put(1,0){\line(0,1){2}}
\put(3,0){\line(0,1){2}}
\put(5,0){\line(0,1){2}}
\put(11,0){\line(0,1){2}}
\put(0.9,1){\makebox(0,0)[r]{$N$}}
\put(2,1){\makebox(0,0){$\mathbb{Z}_{\left(p\left(1\right)\right)}$}}
\put(4,1){\makebox(0,0){$\mathbb{Z}_{\left(p\left(2\right)\right)}$}}
\put(8,1){\makebox(0,0){$\cdots$}}
\put(2,2.1){\makebox(0,0)[b]{$n_{1}$}}
\put(4,2.1){\makebox(0,0)[b]{$n_{2}$}}
\end{picture}%
\,,
\]
where $p\left(  i\right)  $ runs over all elements in $\operatorname*{Prim}%
\left(  \det\left(  A\right)  \right)  $. Its columns consist precisely of a
basis for the null-space of the matrices $E\left(  A\right)  $ taken at each
prime $p\left(  i\right)  $. Two such matrices $M_{1}$, $M_{2}$ represent
isomorphic unordered dimension groups if and only if there is a matrix
$C\in\operatorname*{GL}(N,\mathbb{Z}\left[  1/\det(A)\right]  )$ and an
invertible direct sum of $p\left(  i\right)  $-adic integer matrices $D$ such
that $CM_{1}D=M_{2}$.
\end{corollary}

\begin{proof}
The given structure (that is, $M_{1}$ up to its equivalence with any $CM_{1}%
D$) is an isomorphism invariant because the $p$-adic row spaces $G_{\left(
p\right)  }\left(  A\right)  $ defined in (\ref{eqpad.11})--(\ref{eqpad.12})
are invariants. Corollary \ref{Corpad.2} shows that a rational matrix over
$\mathbb{Z}\left[  1/\det(A)\right]  $ giving an isomorphism on dimension
groups must give an isomorphism on the $p$-adic row spaces, hence the dual
$p$-adic null spaces of $E_{\left(  p\right)  }\left(  A\right)  $. In fact
Corollary \ref{Corpad.2} gives, as necessary and sufficient conditions for
unordered dimension group isomorphism, in effect, the existence of $C$ and
$D$: the $p$-adic symmetries just mean we are considering the row spaces up to
isomorphism, and the $\operatorname*{GL}(N,\mathbb{Z}\left[  1/\det A\right]
)$ symmetry means that we have a rational map which is an isomorphism at all
primes other than the ones considered here.

The extension class of any extension of $\mathbb{Z}^{N}$ by a group
$G/\mathbb{Z}^{N}$ may be computed by extending the map $\mathbb{Z}^{N}%
\subset\mathbb{Q}^{N}$ to a mapping $G\rightarrow\mathbb{Q}^{N}$, and letting
this give a map in $\operatorname*{Hom}(G/\mathbb{Z}^{N},(\mathbb{Q}%
/\mathbb{Z})^{N})\cong\operatorname*{Ext}^{1}(G/\mathbb{Z}^{N},\mathbb{Z}%
^{N})$. This is the remark of Cartan--Eilenberg \cite[p.~292]{CaEi56}. To
identify this class it suffices to look at the $p$-torsion subgroup of
$G/\mathbb{Z}^{N}$ for each prime $p$ since the group is the direct sum of its
$p$-torsion subgroups. To identify this class, take the tensor product of $G$
with the $p$-adic integers, getting a localized extension of $\mathbb{Z}%
_{\left(  p\right)  }^{N}$ by the $p$-torsion subgroup of $G/\mathbb{Z}^{N}$,
which is $G\otimes\mathbb{Z}_{\left(  p\right)  }^{N}$. But if we write all
$p$-adic vectors as the direct sum $K\oplus R$ of the $p$-adic null space of
$R_{A}$ and a complementary space $R$, by Proposition
\textup{\ref{Proposition1}},
\begin{equation}
G\otimes\mathbb{Z}_{\left(  p\right)  }^{N}=(K\otimes\mathbb{Q}_{\left(
p\right)  })\oplus(R\otimes\mathbb{Z}_{\left(  p\right)  }). \label{eqGZN.4}%
\end{equation}
Thus the extension class is represented taking
\begin{equation}
(K\otimes\mathbb{Q}_{\left(  p\right)  })\oplus(R\otimes\mathbb{Z}_{\left(
p\right)  })\rightarrow(K+R)\otimes\mathbb{Q}_{\left(  p\right)  }
\label{eqGZN.5}%
\end{equation}
and collapsing by $\mathbb{Z}_{\left(  p\right)  }^{N}$ to give the inclusion
\begin{equation}
(K\otimes\mathbb{Q}_{\left(  p\right)  })/(K\otimes\mathbb{Z}_{\left(
p\right)  })\rightarrow(K+R)\otimes(\mathbb{Q}_{\left(  p\right)  }%
/\mathbb{Z}_{\left(  p\right)  })=(\mathbb{Q}_{\left(  p\right)  }%
/\mathbb{Z}_{\left(  p\right)  })^{N}. \label{eqGZN.6}%
\end{equation}
This map is induced by the map
\begin{equation}
(K\otimes\mathbb{Q}_{\left(  p\right)  })\rightarrow(K+R)\otimes
\mathbb{Q}_{\left(  p\right)  }=\mathbb{Z}^{N}\otimes\mathbb{Q}_{\left(
p\right)  } \label{eqGZN.7}%
\end{equation}
which can be taken to send the $i$th unit vector on the left to the $i$th
vector in a basis for $K$ on the right. This means taking basis vectors for
the null space of $E_{A}$ as forming the columns of the matrix giving the extension.
\end{proof}

In Example \ref{ExampleJul28} we will give an example where the groups $G/L$
are the same, but the extensions are different.

\section{\label{Red}Remarks on the singular case}

Except for Sections \ref{SY} and \ref{SX} and Theorem \ref{ThmDec.7} of this
paper, we have considered AF-algebras defined by nonsingular primitive
matrices $A,B,\dots$. Let us comment that the class of dimension groups
arising from primitive matrices does change if one dispenses with the
nonsingularity assumption. George Elliott has given an example of a dimension
group which in his terminology is not ultrasimplicial. For example, the group
arising from the stationary diagram associated with the matrix
\[%
\begin{pmatrix}
1 & 1 & 1\\
0 & 1 & 2\\
2 & 1 & 0
\end{pmatrix}
\]
cannot arise from a diagram associated with a nonsingular matrix \cite{Ell79}.
Elliott actually proves the following: If $G$ is the dimension group defined
by the matrix via Bratteli diagrams, then it is not possible to span a given
finite set of positive elements of $G$ by another finite set of positive
linearly independent (over $\mathbb{Z}$) elements, where the span uses only
positive coefficients. But it is not hard to see that this is equivalent to
the impossibility of writing the group as a direct limit with injective maps,
that is, with nonsingular matrices. In the case at stake, Elliott gives the
following three positive elements of $G$ which cannot be spanned in the above
manner: $(1,0)$, $(1,-1)$, $(1,1)$. In this case $G$ is isomorphic to the
direct product of $\mathbb{Z}[1/3]$ and $\mathbb{Z}$, with order defined by
the condition that the positive elements are those elements which have
strictly positive first coordinate---the one from $\mathbb{Z}[1/3]$. The
matrix is primitive.

Since $C^{*}$-equivalence is weaker than shift equivalence, the matrix above
is also an example of a matrix which is primitive, but not shift equivalent to
a primitive nonsingular matrix. Other such examples may be found in
\cite{BoHa91}.

We do not expect that the analogue of our Theorem \ref{Theorem3Sep15} is true
in the context of shift equivalence. (For example, it looks like the example
\cite[Appendix 3, p.~310]{BoHa91} can be modified as follows (multiply
approximately by $10$). There is no $4\times4$ nonnegative matrix whose
spectrum is $\{14,10i,-10i,3\}$. However, we can readily construct an integer
$4\times4$ matrix with this spectrum (block diagonal), with the eigenvector at
$14$ being $(1,0,0,0)$. Then we take a unimodular matrix mapping the vector
$(1,0,0,0)$ to $(1,1,1,1)$. Then some power of this conjugated matrix is
nonnegative, giving us the situation of Theorem \ref{Theorem3Sep15}. Moreover
it should follow from \cite[Theorem 3.1]{BoHa93} that this matrix over
$\mathbb{Z}$ is shift equivalent to a primitive nonnegative matrix of larger
dimension.) For our weaker $C^{*}$-equivalence we have the added flexibility
of replacing matrices by powers of themselves.

We will prove in Theorem \ref{Theorem3Sep15} below that in some special
circumstances, the condition of nonsingularity of the matrix $A$ can be
removed, and $A$ merely assumed to be primitive, without changing the class of
$C^{\ast}$-algebras. In the general case when $A$ is not assumed invertible,
we may introduce the eventual range of $A$,%
\begin{equation}
\mathcal{W}\left(  A\right)  :=\bigcap_{i=0}^{\infty}A^{i}\mathbb{Q}^{N}%
=A^{N}\mathbb{Q}^{N}. \label{eqRed.J1}%
\end{equation}
Note that $A$ is bijective as a map $\mathcal{W}\left(  A\right)
\rightarrow\mathcal{W}\left(  A\right)  $. We may now introduce an additive
group $G\left(  A\right)  $ by%
\begin{equation}
G\left(  A\right)  :=\left\{  g\in\mathcal{W}\left(  A\right)  \mid A^{k}%
g\in\mathbb{Z}^{N}\text{ for some }k\in\mathbb{Z}_{+}\right\}  ,
\label{eqRedNew.1}%
\end{equation}
and one notes that this group $G\left(  A\right)  $ identifies with the
inductive limit of the sequence (\ref{eqInt.6}), i.e. $G\left(  A\right)  $ is
the dimension group when it is equipped with the obvious order. (This version
of $G\left(  A\right)  $ was used, but not defined, already in Section
\ref{SY} above.) Let us give some details. An element of the inductive limit
(\ref{eqInt.6}) can be represented by a sequence $\left\{  g_{m},g_{m+1}%
,\dots\right\}  $ in $\mathbb{Z}^{N}$ with $Ag_{n}=g_{n+1}$ for $n=m,m+1,\dots
$. Two such sequences represent the same element if they coincide from a
certain step $n$ onward. Given such a sequence, there is a unique sequence
$\left\{  h_{1},h_{2},\dots\right\}  $ in $\mathcal{W}\left(  A\right)  $ such
that $Ah_{n}=h_{n+1}$ for $n=1,2,\dots$ and such that $h_{n}=g_{n}$ for all
large $n$. Then $h_{1}$ is the element of $G\left(  A\right)  \subset
\mathcal{W}\left(  A\right)  $ representing the dimension group element in
(\ref{eqRedNew.1}), so this shows the equivalence between the two definitions
(\ref{eqRedNew.1}) and (\ref{eqInt.6}) of $G\left(  A\right)  $. The
definition (\ref{eqRedNew.1}) is the definition used in \cite[p.~49]{BMT87}.
If $A$ is nonsingular and, as in (\ref{eqRedNew.1}), $A^{k}g=m\in
\mathbb{Z}^{N}$, then $g=A^{-k}m$ is a typical element of the (\ref{eqInt.8}%
)--(\ref{eqInt.10}) version of $G\left(  A\right)  $, and \emph{vice versa.}
If $A$ is primitive, we still have the Perron--Frobenius data, and the order
can be defined as before, \emph{mutatis mutandis.}

\begin{lemma}
\label{Lemma1Sep15}Given a vector $u\in\mathbb{R}^{r}$ there exist $r$ vectors
$w_{i}\in\mathbb{Z}^{r}$ such that the convex cone generated by the $w_{i}$
contains an open neighborhood of $u$, and the determinant of the matrix the
$w_{i}$ form is $\pm1$.
\end{lemma}

\begin{proof}
The standard unit vectors do this for any vector in the subsemigroup of
strictly positive integer vectors. We claim transforms of these by integer row
and column operations, permutations, and reversals of sign, take any vector to
the interior of this subsemigroup---then just reverse those operations on the
standard basis vectors. In fact, we get all coordinates nonzero by certain
linear combinations, then reverse their signs.
\end{proof}

\begin{remark}
\label{Remark2Sep15}It is not in general possible to get a determinant-$1$
system of matrices which approximate multiples by some positive constant $C$
of a given set of nonnegative vectors $w_{i}$ in Lemma
\textup{\ref{Lemma1Sep15}}. This is easiest to see when the vectors $w_{i}$
are chosen diagonally dominant. But Lemma \textup{\ref{Lemma1Sep15}} can
probably be strengthened a little.
\end{remark}

\begin{theorem}
\label{Theorem3Sep15} Let $A$ be an integer primitive matrix. Suppose that
when the vector $u$ in Lemma \textup{\ref{Lemma1Sep15}} is the
Perron-Frobenius eigenvector of $A$, the vectors $w_{i}\in\mathbb{Z}^{r}$ can
be chosen to be positive in terms of the order structure of the dimension
group of $A$. Then the ordered dimension group arising from the primitive
integer matrix $A$ is order isomorphic to one arising from a nonsingular
primitive integer matrix $B$.
\end{theorem}

\begin{proof}
Let the dimension of $A$ be $d$ and the rank of all sufficiently large powers
$A^{s}$ be $r$. By Lemma \ref{Lemma1Sep15}, we find a set of $r$ vectors
$w_{i}$ in the eventual row space $R$, that is, the row space of $A^{N}$, or
some specific higher power, a rank-$r$ subspace of $\mathbb{Z}^{d}$ such that
the cone over $\mathbb{Q}_{+}$ generated by this set includes a neighborhood
of the maximum eigenvector $v$ within $R$. This is sufficient to establish
that all sufficiently large powers of $A$ have their rows expressed as
(unique) nonnegative linear combinations of $w_{i}$, since all rows of $A^{s}$
divided by their lengths converge to fixed multiples of $v$ and hence are
eventually in the convex cone; but to be in the convex cone means that we have
these convex combinations.

However, we also need that it can be chosen that these convex combinations are
eventually integer. For that, it suffices that the determinant of the $w_{i}$
expressed as combinations of a basis for the integral vectors in the eventual
row space, i.e., $R\cap\mathbb{Z}^{d}$, a rank-$r$ free abelian group, is $1$
or $-1$. This follows from the lemma and the extra assumption.

Now let $B$ be the matrix of $A^{s}$ expressed as acting on the vectors
$w_{i}$, which will be nonnegative, and positive. Then $B$ is shift equivalent
to $A^{s}$ over the integers (maybe with negative entries), just by the
inclusion mapping given by the vectors $w_{i}$. By a theorem of Parry and
Williams \cite{PaWi77} (reproved in our 1979 paper \cite{KiRo79}), any shift
equivalence over $\mathbb{Z}$ of primitive matrices can be realized by a shift
equivalence over $\mathbb{Z}_{+}$. This shift equivalence will induce an
isomorphism of ordered dimension groups.
\end{proof}

\section{\label{Str}Strong local isomorphism}

\begin{definition}
\label{Def4Sep15}We will say that two dimension groups $G$, $G^{\prime}$ are
strongly locally isomorphic at the prime $p$ if and only if there is an
isomorphism $G\otimes\mathbb{Z}_{(p)}\rightarrow G^{\prime}\otimes
\mathbb{Z}_{(p)}$ induced by a matrix of integers. \textup{(}The first
paragraph of the proof below shows that this is equivalent to requiring that
the isomorphism be induced by a matrix of rational numbers.\textup{) (}Recall
that $G$ and $G^{\prime}$ are locally isomorphic at prime $p$ if there merely
is an isomorphism $G\otimes\mathbb{Z}_{\left(  p\right)  }\rightarrow
G^{\prime}\otimes\mathbb{Z}_{\left(  p\right)  }$.\textup{)}
\end{definition}

In the next theorem we show that strong local isomorphism is described a
condition similar to that in Corollary \ref{Proposition3rev} if we just take
the submatrix corresponding to the prime in question. This condition is rather
strong and can be decided by a simpler algorithm than the general algorithm in
Section \ref{Dec}. When we speak of realizing some $p$-adic construct over an
algebraic number field $K$, we mean in terms of an inclusion $K\subset
\mathbb{Z}_{(p)}$ corresponding to a non-archimedean completion of $K$. Note
in connection with the following theorem that if the ranks of the $p$-adic
eventual row spaces of $A$ and $B$ are different, local isomorphism cannot
hold. Also note that the equivalent conditions imply that the smallest fields
over which the eventual $p$-adic row spaces can be realized are the same for
$A$ and $B$.

\begin{theorem}
\label{Theorem5Sep15}Given two nonnegative matrices $A$, $B$, form a matrix
whose rows are a basis for the $p$-adic eventual row spaces of $A$, $B$ whose
ranks are $n_{p}$. Their dimension groups are strongly locally isomorphic at
the prime $p$ if and only if the corresponding two matrices $A$, $B$ for each
$p$ admit some matrices $C\in\operatorname*{GL}(n,\mathbb{Q})$, $D\in
\operatorname*{GL}(n_{p},\mathbb{Z}_{(p)})$ such that $AC=DB$. This condition
is decidable.
\end{theorem}

\begin{proof}
The proof of Corollary \ref{Corpad.2}, or alternatively, the proof of
Proposition \ref{Proposition1}, shows that having a rational mapping which
induces isomorphism of $p$-adic dimension groups is equivalent to having a
rational mapping which induces an isomorphism of $p$-adic eventual row spaces.
To make such a mapping integer, multiply by all denominators relatively prime
to $p$. At $p$ we must have integrality on the eventual row space $G_{\left(
p\right)  }\left(  A\right)  $ which is a summand of the space of all $p$-adic
integer vectors. The use of a projection $E_{\left(  p\right)  }\left(
A\right)  $ to this subspace enables us to get a map equal to the given
integer mapping defined over some algebraic number field, $p$-adic integer,
and equal to the rational mapping on the eventual row space. The irrational
part of this map will also consist of algebraic integers, since $\mathbb{Z}$
will be an additive summand of the algebraic number ring, and can be
discarded, since it must be zero on the summand. This gives a matrix of
integers inducing the $p$-adic isomorphism, and thus the local isomorphism is strong.

A rational mapping inducing an isomorphism of $p$-adic row spaces is
equivalent to having a matrix $C$ (giving the rational mapping) and a matrix
$D$ such that $DAC=B$ (where $D$ expresses the image of the row basis vectors
of $A$ as linear combinations of a basis for $B$); it must be $p$-adic integer
and invertible since we can also map backwards by the isomorphism of row spaces.

Next we show this criterion is decidable. The field generated by all
eigenvalues of $A$ or $B$ will be sufficient to realize the $p$-adic eventual
row spaces: we take this field and some prime $\pi$ giving an embedding in an
extension of $\mathbb{Z}_{(p)}$. Then the $p$-adic eventual row space is
spanned by the generalized eigenspaces for eigenvalues which are relatively
prime to $\pi$, since multiplication by powers of $A$ will not send them to
zero, but will annihilate all other generalized eigenspaces, modulo any power
of $\pi$. Some linear combinations of basis vectors for these generalized
eigenspaces must give a $p$-adic basis.

The required $p$-adic matrix $D$ then must lie in this field $K$, since $C$,
$A$, $B$ do. Being $p$-adic integer means that its denominators are relatively
prime to $p$. It can then be expanded as a larger matrix over $\mathbb{Z}%
_{ip}=\mathbb{Z}[1/2,1/3,\ldots,\widehat{1/p},\ldots]$, using a basis for the
$p$-adic integers of $K$ over $\mathbb{Z}_{ip}$. Given the matrix $D$, the
condition that a corresponding $C$ exists over $\mathbb{Q}$ can be stated by
linear equations in $D$. Existence of rational $C$ means that all columns of
$BD$ are rational combinations of the columns of $A$. For the expanded
matrices over $\mathbb{Z}_{ip}$, that means that we have any linear column
combinations of the columns of $A$ yielding the desired columns of $BD^{-1}$.
In turn that means that the columns of $BD^{-1}$ have zero inner product with
all vectors which have zero inner product with the columns of $A$, which is a
linear condition.

We can now determine a $p$-adic basis for this linear space of matrices $D$,
write the determinant of $D$ as a polynomial in the coefficients of a general
linear combination of basis elements, and determine whether or not it is
possible for the determinant to be nonzero modulo $p$ by testing each
congruence class of entries modulo $p$.
\end{proof}

\begin{example}
\label{ExampleJul28}Let
\begin{equation}
A=%
\begin{pmatrix}
4 & 1\\
1 & 2
\end{pmatrix}
,\qquad B=%
\begin{pmatrix}
7 & 0\\
0 & 1
\end{pmatrix}
. \label{eqExampleJul28}%
\end{equation}
The respective characteristic polynomials are $x^{2}-6x+7$ and $x^{2}-8x+7$,
with determinant $7$, and we consider the local dimension groups at $7$. Since
$7$ does not divide $8$, only one root of the former polynomial is divisible
by $7$. Thus only the identity element of the $\mathbb{Z}_{2}$ Galois group
fixes the eigenvalue not divisible by $7$. This implies that the $7$-adic row
space is irrational, and the minimal fields over which eventual row spaces are
defined are respectively $\mathbb{Q}\left[  \sqrt{2}\right]  $, $\mathbb{Q}$,
so the dimension groups are not locally isomorphic. Note that, even so, the
two quotient groups $G\left(  A\right)  /\mathbb{Z}^{2}$ and $G\left(
B\right)  /\mathbb{Z}^{2}$ are isomorphic. This follows from Proposition
\textup{\ref{Proposition1}:} Recall, to verify this we need only compute the
respective Ulm numbers from the characteristic polynomials, and there is only
the prime $p=7$ to check. So the $7$-reduced rank is $1$ for each of the two
quotient torsion groups calculated from $A$ and $B$.
\end{example}

\begin{example}
\label{ExamplebisJul28}Our next example illustrates $C^{\ast}$-symmetry, as
well as the calculation of the $p$-adic eventual row spaces. Consider
\begin{equation}
A=%
\begin{pmatrix}
3 & 1\\
2 & 3
\end{pmatrix}
\label{eqExamplebisJul28.1}%
\end{equation}
and its transpose $B=A^{\operatorname*{tr}}$; for this particular choice of
the pair $A$, $B$, the minimal fields over which eventual row spaces are
defined are isomorphic. We can arbitrarily choose which root of the
characteristic polynomial $x^{2}-6x+7$, the same as for $A$ in Example
\textup{\ref{ExampleJul28},} is divisible by $7$ \textup{(}representing the
unique $p$-adic root which is divisible by $7$\textup{),} say $3-\rho$ where
$\rho$ is a square root of $2$ in $\mathbb{Z}_{\left(  7\right)  }$
\textup{(}see next paragraph\textup{).} The eventual row eigenspaces are
spanned by the other row eigenvectors, which are $\left(  1,\rho\right)  $ for
$B$ and $\left(  \rho,1\right)  $ for $A$. A mapping of eigenspaces must map
one to a multiple $c$ times the other. If it commutes with the Galois action,
then it must do the same for their conjugates, so that it has the form
\begin{equation}%
\begin{pmatrix}
0 & c\\
c & 0
\end{pmatrix}
. \label{eqExamplebisJul28.2}%
\end{equation}
The determinant restricted to the eventual $7$-adic row space is $c$, so the
congruences are $\,c\equiv0\pmod{7}$, which are solvable. The dimension groups
of this matrix and its transpose are locally isomorphic at the prime $7$.
Since $7$ is the only prime involved, this implies global isomorphism of the
unordered dimension groups.

Here we think of the field $\mathbb{Q}\left[  \sqrt{2}\right]  $ as embedded
in the field of $7$-adics $\mathbb{Q}_{\left(  7\right)  }$ via $1\mapsto1$
and $\sqrt{2}\mapsto\rho$, where $\rho\in\mathbb{Z}_{\left(  7\right)
}\subset\mathbb{Q}_{\left(  7\right)  }$. The polynomial $x^{2}-2$ may be
considered in $\mathbb{Z}_{\left(  7\right)  }\left[  x\right]  $, and it is
reducible there. To find the two roots $\pm\rho$ in $\mathbb{Z}_{\left(
7\right)  }$, calculate the terms $t_{0},t_{1},\dots\in\left\{  0,1,2,\dots
,6\right\}  $ in $\rho=t_{0}+t_{1}\cdot7+t_{2}\cdot7^{2}+\cdots\in
\mathbb{Z}_{\left(  7\right)  }$ recursively, starting with $t_{0}=3$ or
$t_{0}=4$ \textup{(}see, e.g., \cite[p.~18]{BoSh66}\textup{). They can be
found with Maple or the PARI program. The results are }$\rho=3+1\cdot
7+2\cdot7^{2}+6\cdot7^{3}+1\cdot7^{4}+2\cdot7^{5}+\cdots$ and $-\rho
=4+5\cdot7+4\cdot7^{2}+0\cdot7^{3}+5\cdot7^{4}+4\cdot7^{5}+\cdots$. Since the
root $3-\rho$ of $x^{2}-6x+7$ is divisible by $7$, the $7$-adic eventual row
spaces in $\left(  \mathbb{Z}_{\left(  7\right)  }\right)  ^{2}$ are the
respective $\mathbb{Z}_{\left(  7\right)  }$-modules%
\[
G_{\left(  7\right)  }\left(  A\right)  =\mathbb{Z}_{\left(  7\right)
}\left(  \rho,1\right)
\]
and%
\[
G_{\left(  7\right)  }\left(  B\right)  =\mathbb{Z}_{\left(  7\right)
}\left(  1,\rho\right)  ,
\]
by Remark \textup{\ref{RempadNew.3},} i.e., they are generated over
$\mathbb{Z}_{\left(  7\right)  }$ by the respective row eigenvectors
corresponding to the second eigenvalue $3+\rho$, the one not divisible by $7$.

Note also in this case that $A$ and $A^{\operatorname*{tr}}$ are conjugate by
the unimodular matrix $\left(
\begin{smallmatrix}
0 & 1\\
1 & 0
\end{smallmatrix}
\right)  $; and hence it even follows directly as in Example
\textup{\ref{Exa1}} that $A$ and $A^{\operatorname*{tr}}$ are \textup{(}%
elementary\/\textup{)} shift equivalent.

Note finally that if one denotes the $A$ matrices in
\textup{(\ref{eqExampleJul28})} and \textup{(\ref{eqExamplebisJul28.1})} by
$A_{1}$, $A_{2}$, respectively, and one defines%
\begin{equation}
J=
\begin{pmatrix}
1 & 0\\
1 & 1
\end{pmatrix}
, \label{eqStrNew.4}%
\end{equation}
then $A_{2}J=JA_{1}$, and thus if $K=A_{1}J^{-1}$ we have the system%
\begin{equation}
A_{2}=JK,\qquad A_{1}=KJ. \label{eqStrNew.5}%
\end{equation}
But $K$ does not have positive matrix entries, so this does not imply
elementary shift equivalence. However, if we redefine%
\begin{equation}
K=A_{1}^{2}J_{{}}^{-1}=
\begin{pmatrix}
11 & 6\\
1 & 3
\end{pmatrix}
, \label{eqStrNew.6}%
\end{equation}
then we have the pair of shift relations for the squares,%
\begin{equation}
A_{1}^{2}=KJ,\qquad A_{2}^{2}=JK, \label{eqStrNew.7}%
\end{equation}
which is the assertion that $A_{1}^{2}$ and $A_{2}^{2}$ are elementary shift
equivalent. In particular, $A_{1}$ and $A_{2}$ are $C^{\ast}$-equivalent. This
latter conclusion and the one in Example \textup{\ref{Exa1}} also follow the
next general observation:
\end{example}

\begin{observation}
\label{ObsStrNew.5}If $A$, $B$ are nonsingular primitive $N\times N$ matrices
and there exists a unimodular matrix $J$ in $M_{N}\left(  \mathbb{Z}\right)  $
such that%
\begin{equation}
v\left(  B\right)  J=\mu v\left(  A\right)  \label{eqStrNew.8}%
\end{equation}
for a positive number $\mu$, and%
\begin{equation}
BJ=JA, \label{eqStrNew.9}%
\end{equation}
then $A$ and $B$ are $C^{\ast}$-equivalent.
\end{observation}

\begin{proof}
Since $J$ is unimodular, we have%
\begin{equation}%
\begin{cases}
B^{n}JA^{-n}=JA^{n}A^{-n}=J\in M_{N}\left( \mathbb{Z}\right) , \\
A^{n}J^{-1}B^{-n}=A^{n}A^{-n}J^{-1}=J^{-1}\in M_{N}\left( \mathbb{Z}\right) ,
\end{cases}
\label{eqStrNew.10}%
\end{equation}
and the observation follows from (\ref{eqInt.17})--(\ref{eqInt.18}). (The
condition (\ref{eqStrNew.8}) may be replaced by the strictly stronger
requirement that $J$ and $J^{-1}$ have only nonnegative matrix entries.)
\end{proof}

\begin{remark}
\label{RemStrOct.6}In fact we have the ``partial'' implication
\textup{(\ref{eqStrNew.9})} $\Rightarrow$ \textup{(\ref{eqStrNew.8}),} but
\textup{(\ref{eqStrNew.8})} for some real scalar $\mu$, while the positivity
restriction on $\mu$ is not a consequence of \textup{(\ref{eqStrNew.9})}
alone. We further stress that \textup{(\ref{eqStrNew.9})--(\ref{eqStrNew.8})}
are more restrictive than $C^{\ast}$-equivalence, even more restrictive than
shift equivalence: take, for example, $A=\left(
\begin{smallmatrix}
2 & 1\\
4 & 4
\end{smallmatrix}
\right)  $ and $B=\left(
\begin{smallmatrix}
2 & 2\\
2 & 4
\end{smallmatrix}
\right)  $, which are shift equivalent by \cite{Bak83}, but do not satisfy
\textup{(\ref{eqStrNew.9}).}
\end{remark}

To summarize, the two examples have four matrices in all, and the first one in
Example \ref{ExampleJul28} is $C^{\ast}$-equivalent to the two in Example
\ref{ExamplebisJul28}, but $\left(
\begin{smallmatrix}
7 & 0\\
0 & 1
\end{smallmatrix}
\right)  $ from Example \ref{ExampleJul28} is not $C^{\ast}$-equivalent to the
other three. The first one, $\left(
\begin{smallmatrix}
4 & 1\\
1 & 2
\end{smallmatrix}
\right)  $ in Example \ref{ExampleJul28}, is symmetric, and $A$ from Example
\ref{ExamplebisJul28} is $C^{\ast}$-symmetric in that it is $C^{\ast}%
$-equivalent to its own transpose.

\begin{remark}
\label{RemStrOct.7}Note that the two matrices $A_{1}$, $A_{2}$ in
(\ref{eqExampleJul28}), (\ref{eqExamplebisJul28.1}) considered above are
elementary shift equivalent over $\mathbb{Z}$ since they are conjugate over
$\mathbb{Z}$. But while $A_{1}^{2}$, $A_{2}^{2}$ are elementary shift
equivalent over $\mathbb{Z}_{+}$, $A_{1}$ and $A_{2}$ are not! (These types of
$2\times2$ examples have been considered earlier by Kirby Baker
\cite{Bak83,Bak87}.) This is seen as follows: Suppose $A_{1}=CD$ where $C$,
$D$ are nonnegative integer $2\times2$ matrices. Then $C$ expresses the rows
of $A_{1}$ as nonnegative integer combinations of the rows of $D$. The entries
$1$ in the rows of $A_{1}$ can only come from entries $1$ in the rows of $D$.
Moreover these $1$'s can only be in the same row. Furthermore, in the linear
combinations these $1$'s can only be multiplied by $1$'s. So the product $CD$
looks like, up to symmetry,%
\begin{equation}%
\begin{pmatrix}
1\cdot c_{12} & d_{11}\cdot1\\
c_{21}\cdot1 & 1\cdot d_{22}%
\end{pmatrix}
=%
\begin{pmatrix}
4 & 1\\
1 & 2
\end{pmatrix}
. \label{eqStrNew.11}%
\end{equation}
But if we write out the equations, there are no solutions unless one of $C$,
$D$ is a permutation matrix, and thus $DC$ cannot be equal to $A_{2}$.
\end{remark}

\section{\label{Con}Concluding remarks}

In the paper we addressed the interplay between the local and the global
versions of the isomorphism problem. There are different, but related,
decidability results in the literature. Ax and Kochen
\cite{AxKo65a,AxKo65b,AxKo66} and Grunewald and Segal \cite{GrSe82} address
decidability in a $p$-adic setting.

\begin{acknowledgements}
The co-authors are very grateful to Brian Treadway for his excellent work in
typesetting, and in coordinating the many pieces of manuscript, and sequences
of revisions, which arrived by fax and e-mail. We also thank Daniele Mundici
and David Stewart for helpful conversations and references on computation and
algorithms. We are especially indebted to Vincenzo Marra for pointing out a
serious mistake in the preprint version of Theorem \ref{Theorem3Sep15}, and
reminding us about the reference \cite{Ell79}. The referee of the paper sent a
very extensive report making many constructive suggestions both for improving
the exposition and for making the reduction to nonsingular matrices explicit.
P.E.T.J. benefited by a Norwegian-funded visit to the University of Oslo in
the winter 1998--99 and in the summer 2000 where part of the work was done,
and he is grateful for the support and hospitality.\bigskip
\end{acknowledgements}

\begin{center}
\textit{Richtiges Auffassen einer Sache und Mi{\ss}verstehen der gleichen
Sache}\linebreak \textit{schlie{\ss}en einander nicht vollst\"{a}ndig
aus.\bigskip}\linebreak \settowidth{\qedskip}{\textit{Richtiges Auffassen
einer Sache und Missverstehen der gleichen Sache}}\makebox[\qedskip
]{\hfill\textsc{---Franz Kafka}, \textit{Der Proze\ss}}
\end{center}

\bibliographystyle{bftalpha}
\bibliography{jorgen}
\end{document}